\documentclass[a4paper,10pt]{article}
\usepackage{latexsym}

\usepackage{amsmath,amscd,amssymb}
\usepackage{graphicx}
\addtocounter{section}{-1}

\usepackage{theorem}
\newtheorem{lemma}{Lemma}[section]
\newtheorem{proposition}{Proposition}[section]
\newtheorem{theorem}{Theorem}[section]
\newtheorem{corollary}{Corollary}[section]
\newtheorem{example}{Example}[section]

\def\BC{\begin{center}}
\def\EC{\end{center}}
\def\BP{\noindent {\it Proof.} }
\def\EP{\hspace*{\fill}$\Box$

\vspace{1ex}
}

\def\BT{\begin{theorem}}
\def\ET{\end{theorem}}
\def\BPR{\begin{proposition}}
\def\EPR{\end{proposition}}
\def\BL{\begin{lemma}}
\def\EL{\end{lemma}}
\def\BCO{\begin{corollary}}
\def\ECO{\end{corollary}}
\def\BPI{\BC\begin{picture}}
\def\EPI{\end{picture}\EC}
\def\BTA{\begin{tabular}}
\def\ETA{\end{tabular}}
\def\BX{\begin{example}}
\def\EX{\end{example}}
\def\BIB{\vspace{-1ex}\bibitem}
\def\BIT{\begin{itemize}}
\def\EIT{\vspace{-1ex}\end{itemize}}

\def\BM{
\vspace{1ex}

\ \ \ $}

\def\EM{$

\vspace{1ex}
\noindent}

\begin{document}
\title{\bf Exponential functions of finite posets and   
the number of extensions with a fixed set of minimal points}
\author{\sc Frank a Campo {\rm and} Marcel Ern\'e}  
\date{\small Seilerwall 33, D 41747 Viersen, Germany\\
{\sf acampo.frank@gmail.com}\\[3pt]
Faculty for Mathematics and Physics,
Leibniz University, Welfengarten 1,\\ D 30167 Hannover, Germany\\ 
{\sf erne@math.uni-hannover.de}\\[3pt]
}

\maketitle

\begin{abstract}
\noindent  We establish formulas for the number of all downsets (or equivalently, of all antichains) of a finite poset $P$. 
Then, using these numbers, we determine recursively and explicitly the number of all posets 
having a fixed set of minimal points and inducing the poset $P$ on the non-minimal points. 
It turns out that these counting functions are closely related to a collection of downset numbers of certain subposets. Since any function of that kind is an exponential sum 
(with the number of minimal points as exponent), we call it the exponential function of the poset.
Some linear equations, divisibility relations, upper and lower bounds, and asymptotical equalities for the counting functions are deduced. 
A list of all such exponential functions for posets with up to five points concludes the paper.

\vspace{1ex}


\end{abstract}

\section{\large Introduction}
\label{intro}

The determination of $p(k)$, the number of all labeled partially ordered sets (posets), and $\overline{p}(k)$, the number of all unlabeled
(i.e.\ isomorphism types of) posets with $k$ points are two old order-theoretical problems that are still far from being completely solved. 
For early attempts see, e.g., Sharp \cite{Sh1}, Butler and Markowsky \cite{BM1,BM2}, Ern\'e \cite{E1,E2}. 
Numerical and structural improvements until the last decade of the past century are reported by Ern\'e and Stege \cite{ESt1}.
The explicit computation of the numbers $p(k)$ was pushed forward by Heitzig and Reinhold \cite{HR1} 
until $k \!=\!17$, and of $\overline{p}(k)$ until $k \!=\!14$. Slightly later Brinkmann and McKay \cite{BM1,BM2} 
calculated $p(k)$ for $ k\leq 19$, $\overline{p}(k)$ for $k\leq 16$, and many tables of related numbers. 
The reason why in the labeled case one can go three steps further than in the unlabeled case is a result due to Ern\'e \cite{E1, E3}, saying that 
$p(k)$ equals the number of posets on $k\!+\!1$ points with a fixed set of exactly two minimal (or maximal) elements.

It seems that a way to a general solution (if any) can only be successful if one considers posets with further additional parameters, 
like the number of related pairs (Culberson and Rawlins \cite{CR}, Ern\'e and Stege \cite{E4,ESt2}), or of downsets, upsets or antichains (Ern\'e \cite{E3}, Stanley \cite{St2}). 
{\em Downsets} ({\em alias lower sets} or {\em decreasing sets}) of a poset $P$ contain with any element all elements below it, 
and {\em upsets} ({\em alias upper sets} or {\em increasing sets}) are defined dually.
The downsets and the upsets of $P$ form two dually isomorphic superalgebraic, hence completetly distributive lattices $\mathcal{D}(P)$ and $\mathcal{U}(P)$. 
Functorial properties of $\mathcal{D}$ and $\mathcal{U}$ with some combinatorial consequences will be discussed briefly in Sections \ref{basic} and \ref{sums}. 
For example, the tensor product of $\mathcal{U}(O)$ and $\mathcal{D}(P)$ turns out to be isomorphic to the lattice of all generalized vertical sums of $O$ and $P$.  
The numbers $d(P) = \sharp\, \mathcal{D}(P) = \sharp\, \mathcal{U}(P)$ of all up- resp.\ downsets of $P$ will play a crucial role in the present study. 

Given natural numbers $j$ and $k$, let $p(j,k)$ denote the number of all labeled posets and $\overline{p}(j,k)$ that of all unlabeled posets 
on $k$ points having exactly $j$ downsets (or upsets, or antichains). 
By the familiar one-to-one correspondence between finite posets and finite $T_0$ spaces (Alexandroff \cite{Al1,Al2}, Birkhoff \cite{BLat}), 
$p(j,k)$ resp.\ $\overline{p}(j,k)$ is also the number of all $T_0$ topologies on $k$ points with exactly $j$ open (or closed) sets.  
Some first progress towards the determination of these numbers 
(at least at the lower or upper end of the $j$-scale for fixed $k$) was made by Sharp \cite{Sh2}, Ern\'e \cite{E1,E3} and Stanley \cite{St2}. 
More comprehensive results were obtained by Ern\'e and Stege \cite{ESt1,ESt2}, and later by Benoumhani and Kolli \cite{Be, BeKo,Ko1,Ko2}. 
Note that the number $\overline{p}(j,k)$ is, by a theorem due to Birkhoff \cite{BLat}, 
equal to the number of all unlabeled distributive lattices with $j$ elements of which  $k$ are join (or meet) irreducible; 
these numbers have been calculated by Ern\'e, Heitzig and Reinhold \cite{EHR} for all $j \!<\! 50$ and arbitrary $k$.

Many of the (often tedious and time-consuming) calculations in that area are simplified essentially 
by the following systematical approach, which we shall pursue in the present paper. 
Given a poset or partial order $P$ on a finite set $K$ and a finite set $M$ disjoint from $K$, let $\mathcal{E}(M,P)$ 
denote the set of all posets with ground set $M\cup K$, inducing $P$ on $K$ and having $M$ as the set of all minimal elements. 
The cardinality $e(m,P)$ of $\mathcal{E}(M,P)$ depends only on the cardinality $m = \sharp\, M$ but not on the specific choice of $M$; 
$e(m,P)$ is called the {\em exponential function} of $P$, because it is always an exponential sum with exponent $m$. 
Using some further combinatorial interpretations of the numbers $e(m,P)$ by means of certain sets of order-preserving functions or relations, 
we shall obtain diverse explicit and recursive formulas for $e(m,P)$, all involving the downset numbers $d(U)$ of certain subposets $U$ as the relevant ingredients. 
Horizontal and vertical decompositions of posets are of particular use for that purpose.
A few number-theoretical divisibility properties of the numbers $e(m,P)$ are easy consequences of their exponential sum representations.
We present exact upper bounds for the downset numbers $d(P)$ where $P$ belongs to specific classes of posets with a given height or number of minimal elements. 
These bounds are helpful for a more precise analysis of the exponential functions $e(m,P)$. 

Given a suitable enumeration of unlabeled finite posets, one may express the relevant formulas in terms of certain succinct matrix equations. 
Furthermore, the knowledge of the leading terms $d(P)^m$ of the functions $e(m,P)$ provides good numerical estimates and information about their asymptotical behavior.
The paper concludes with a list of all unlabeled posets with up to five points and their exponential functions.

\section{\large Basic definitions and concepts}
\label{basic}

By $\mathbb{N}$ we denote the set of positive integers, and $\underline{k}$ will stand for the set of all $n\in \mathbb{N}$ with $n\leq k$\,; 
in particular, $\underline{0}$ is the empty set $\emptyset$. While $\subseteq$ denotes set inclusion, $\subset$ denotes {\em proper} inclusion.

A {\em relation} between (elements of) two sets $X$ and $Y$ is merely a subset $R$ of the cartesian product $X\times Y$. 
Instead of $(x,y) \in R$ one writes $x\,R\,y$. The {\em dual} or {\em converse} relation of $R$ is $R^d = R^{-1} = \{ (y,x)\mid x \, R\, y\}$. For $a\in X$ and $b\in Y$, we put 
\[
aR = \{ y \in Y \mid a\,R\,y\},  \ Rb = \{ x\in X \mid x\,R \,b\},
\]
and for $A\subseteq X$, $B\subseteq Y$, we define
\[
\textstyle{A R = \bigcup\, \{ aR \mid a\in A\}, \ R B = \bigcup\,\{ Rb \mid b\in B\}.}
\]
In case $X = Y$, a subset $R$ of $X\times X$ is a {\em relation on} $X$, and for $A \subseteq X$, the {\em relation induced on} $A$ is
\mbox{$R|_A = R \mathop{\cap}\, (A\times A)$.} Instead of $R|_{X\smallsetminus A}$ we simply write $R\!\mathop{-}\!A$. The {\em identity relation} on $X$ or {\em diagonal\,} is 
\mbox{$\Delta_X \!= \{ (x,x) \mid x\in X\}.$}

Given relations $R\subseteq X\times Y$ and $S\subseteq Y\times Z$, the {\em composite relation} $RS = S\circ R$ consists of all pairs $(x,z)$ for which there exists a $y$ with $x \, R \, y \, S\, z$. 
In particular, this applies to functions from $X$ to $Y$, regarded as special relations $F\subseteq X\times Y$ such that each $xF$ is a singleton, consisting of the single element $F(x)$;
in that case one writes $F[A]$ instead of $AF$.

A {\em quasiorder} (or {\em preorder}) on $X$ is a reflexive and transitive relation on\,\,$X$, in other terms, 
a relation $Q$ on $X$ with $\Delta_X \cup QQ\subseteq Q$; a quasiorder $Q$ that is antisymmetric (i.e.\
$Q\cap Q^{d} \subseteq \Delta_X$) is referred to as a {\em (partial) order} on\,\,$X$. 
The pair $(X,Q)$ is then a {\em quasiordered set} ({\em qoset}) or a {\em partially ordered set} ({\em poset}), respectively.

For any quasiorder $Q$ on $X$ and arbitrary subsets $A \!\subseteq\! X$, the {\em down-closure} is $QA$, the downset generated by $A$ (that is, the least downset containing $A$). 
{\em Upsets} or {\em upper sets} and the {\em up-closure} $AQ$ are defined dually. 
The collection of all downsets of $Q$ is denoted by $\mathcal{D} (Q)$, and that of all upsets by $\mathcal{U}(Q)$. 
It is known and easy to see that $\mathcal{D}(Q)$ as well as $\mathcal{U}(Q)$ is always an {\em Alexandroff topology} (closed under arbitrary unions and intersections),
called the {\em lower} and the {\em upper Alexandroff topology} associated with $Q$ (see \cite{Al1,Al2}). Of course, on a finite set, all topologies are Alexandroff topologies. 
Notice that sending each downset to its complement gives a dual isomorphism between the complete lattices $\mathcal{D}(Q)$ and $\mathcal{U}(Q)$.
By passing from a qoset to the associated (upper or lower) Alexandroff space, one obtains a concrete isomorphism between the category of qosets with isotone, 
i.e.\ order-preserving maps as morphisms and the category of Alexandroff spaces with continuous maps as morphisms. 
Under that isomorphism, the posets correspond to the $T_0$ Alexandroff spaces (in which distinct points have different closures). 
In particular, finite qosets (resp.\ posets) are in one-to-one correspondence with finite topological spaces (resp.\ $T_0$ spaces).

Two fundamental constructions allow us to build from given posets larger ones (cf.\ Stanley \cite{St5}). Suppose $O$ and $P$ are partial orders on
disjoint ground sets $X$ and $Y$, respectively (of course, the disjointness condition is
not an essential restriction, because one always may pass to isomorphic copies having disjoint ground sets).  
The {\em cardinal sum} or {\em disjoint sum} is 
\[O+P = O\mathop{\cup} P ,\vspace{-1ex}\]
and the {\em ordinal sum} is 
\[O\oplus P = O \cup P \cup (X\!\times\! Y).\]
Accordingly, one writes for the corresponding posets
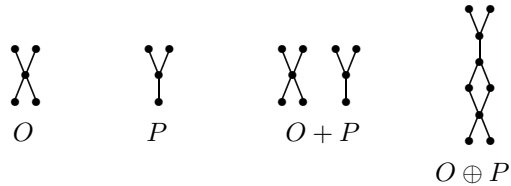
\[ (X,O)+(Y,P) = (X\cup Y, O+P),\]
\[ (X,O)\oplus (Y,P) = (X\cup Y, O\oplus P).\]

\BC
\begin{picture}(200,60)
\put(0,0){
\begin{picture}(10,20)(0,-20)
    \put(.8,20){\circle*{3}}
    \put(8.7,20){\circle*{3}}
    \put(4.7,10){\circle*{3}}
    \put(.8,0){\circle*{3}}
    \put(8.7,0){\circle*{3}}
    \put(.5,20.5){\line(2,-5){4}}
    \put(9,20.5){\line(-2,-5){4}}
    \put(5,10){\line(2,-5){4}}
    \put(4.5,10){\line(-2,-5){4}}
    \put(0,-15){$O$}
\end{picture}
}
\put(50,0){
\begin{picture}(10,20)(0,-20)
    \put(.8,20){\circle*{3}}
    \put(8.7,20){\circle*{3}}
    \put(4.7,10){\circle*{3}}
    \put(4.7,0){\circle*{3}}
    \put(.5,20.5){\line(2,-5){4}}
    \put(9,20.5){\line(-2,-5){4}}
    \put(4.7,0){\line(0,1){10}}
    \put(0,-15){$P$}
\end{picture}
}
\put(100,0){
\begin{picture}(10,20)(0,-20)
    \put(.8,20){\circle*{3}}
    \put(8.7,20){\circle*{3}}
    \put(4.7,10){\circle*{3}}
    \put(.8,0){\circle*{3}}
    \put(8.7,0){\circle*{3}}
    \put(.5,20.5){\line(2,-5){4}}
    \put(9,20.5){\line(-2,-5){4}}
    \put(5,10){\line(2,-5){4}}
    \put(4.5,10){\line(-2,-5){4}}
    \put(2,-15){$O+P$}
\end{picture}
}
\put(120,0){
\begin{picture}(10,20)(0,-20)
    \put(.8,20){\circle*{3}}
    \put(8.7,20){\circle*{3}}
    \put(4.7,10){\circle*{3}}
    \put(4.7,0){\circle*{3}}
    \put(.5,20.5){\line(2,-5){4}}
    \put(9,20.5){\line(-2,-5){4}}
    \put(4.7,0){\line(0,1){10}}
\end{picture}
}
\put(170,15){
\begin{picture}(10,20)(0,-20)
    \put(.8,20){\circle*{3}}
    \put(8.7,20){\circle*{3}}
    \put(4.7,10){\circle*{3}}
    \put(4.7,0){\circle*{3}}
    \put(.5,20.5){\line(2,-5){4}}
    \put(9,20.5){\line(-2,-5){4}}
    \put(4.7,0){\line(0,1){10}}
\end{picture}
}
\put(170,-15){
\begin{picture}(10,20)(0,-20)
    \put(.8,20){\circle*{3}}
    \put(8.7,20){\circle*{3}}
    \put(4.7,10){\circle*{3}}
    \put(.8,0){\circle*{3}}
    \put(8.7,0){\circle*{3}}
    \put(.5,20.5){\line(2,-5){4}}
    \put(9,20.5){\line(-2,-5){4}}
    \put(5,10){\line(2,-5){4}}
    \put(4.5,10){\line(-2,-5){4}}
    \put(-12,-15){$O\oplus P$}
    \put(5,30){\line(2,-5){4}}
    \put(4.5,30){\line(-2,-5){4}}
\end{picture}
}
\end{picture}
\EC

\vspace{2ex}

Up to isomorphism, there is only one {\em antichain} (totally unordered set) with $k$ elements, namely $A_k = (\underline{k}, =)$,
and only one {\em chain} (totally ordered set) with $k$ elements, namely $C_k = (\underline{k},\leq)$. Such extremal posets are obtained by iterated formation of cardinal resp.\ ordinal sums,
starting with singletons. Allowing a mixed iterated construction involving both kinds of sums, again starting with singletons, 
yields exactly all {\em N-free posets}, i.e., all finite posets having no subposet with diagram 
$ $ \begin{picture}(10,8)(0,-1)
    \put(-.2,5.3){\circle*{3}}
    \put(6.2,5.3){\circle*{3}}
    \put(-.2,-3.5){\circle*{3}}
    \put(6.2,-3.5){\circle*{3}}
    \put(6,-5){\line(-1,2){6}}
    \put(-.2,-4.9){\line(0,1){10}}
    \put(6.2,-4.9){\line(0,1){10}}
    \end{picture}
$\!\!$. By obvious reasons, such order relations have also been baptized {\em series-parallel orders} (see, e.g., M\"ohring \cite{Mo}). 

Sometimes orders are identified with the corresponding ordered sets; set-theoretically, this convention is not entirely consistent, but
there is no danger of confusion since the ground set is determined by the order relation. 
If $P$ is a poset with underlying set $X$, the order relation is denoted by $\leq_P$. 
For subsets $A$ of $X$, we denote by $P|_A$ the induced poset on $A$ and by $P\!\mathop{-}\!A$ the induced poset on $X\!\smallsetminus\! A$;  
and we say $P$ is an {\em extension} of a poset $O$ if $O = P|_A$ for some $A$. 
It will then be convenient to write $xP$ for $x\!\leq_P$, $AP$ for $A\!\!\leq_P$\,\,etc.; hence, the symbols $x\leq_P y$, $x\in Py$, and $y\in xP$ all have the same meaning, namely $(x,y)\in \ \leq_P$. 
For any poset (or qoset) $P$, the downset lattice $\mathcal{D}(P) = \mathcal{D}(\leq_P)$ and the upset lattice $\mathcal{U} (P) = \mathcal{U}(\leq_P)$ are {\em superalgebraic}, i.e.\ completely distributive algebraic lattices. 
in fact, $\mathcal{D} (P)$ is the free completely distributive lattice over $P$: 
every isotone (i.e\ order-preserving) map $f$ from $P$ to a completely distributive lattice $L$ uniquely extends to a join-preserving map $f^{\vee}$ from $\mathcal{D}(P)$ to $L$, 
given by  $f^{\vee}(D) = \bigvee f[D]$. Moreover, $\mathcal{D}$ gives rise to a functorial dual equivalence between the category of posets with isotone maps 
and the category of superalgebraic lattices with complete homomorphisms, i.e.\ maps preserving arbitrary joins and meets (see, e.g., Ern\'e \cite{ELR,EAO}). 
As nonempty finite distributive lattices are already superalgebraic, the category of finite posets resp.\ $T_0$ spaces 
is equivalent to the category of finite distributive lattices (see Birkhoff \cite{BLat}).

These remarks enable us to translate all results about (finite) ordered sets into the language of topology (or, if desired, of lattice theory).  
For instance, the down-closure $PA$ of a subset $A$ of a poset or qoset $P$ is the closure w.r.t.\ the associated {\em upper} Alexandroff topology $\mathcal{U}(P)$, 
while the up-closure $AP$ is the closure w.r.t.\ the associated {\em lower} Alexandroff topology $\mathcal{D}(P)$. 
And, consequently, $X\smallsetminus (X\!\smallsetminus \!A)P$ is the topological interior $A^{\circ}$ in $\mathcal{D}(P)$,
whereas $X\smallsetminus P(X\!\smallsetminus \!A)$ is the interior in $\mathcal{U}(P)$. See \cite{E2, EABC, EO, ESt1} for more background.

Although for concrete computations in combinatorial order theory incidence matrices provide a very useful tool  
(see, e.g., \cite{Bu, BuMa,E1,EO,Pa,Sh1,Sh2} for details), we shall not need them in the subsequent considerations.

\newpage

\section{\large Generalized vertical sums of posets}
\label{sums}
 
In this section it will be convenient to denote by $O$ and $P$ two {\em partial orders} on sets $X$ and $Y$, respectively (and not the corresponding posets). 
Recall that a function $f$ between $(X,O)$ and $(Y,P)$ is said to be {\em isotone} (or {\em order-preserving}) 
if $x\, O \,x'$ implies $f(x) \, P \, f(x')$, and {\em antitone} (or {\em order-inverting})
if $x\, O \,x'$ implies $f(x') \, P \, f(x)$. Let $\mathcal{Q}(X)$ denote the set of all quasiorders and $\mathcal{Q}_0(X)$ that of all partial orders on\,$X$.
 
Specifically, given partial orders $O$ and $P$ on {\em disjoint} ground sets $X$ and $Y$, respectively, consider the following sets: 
\begin{eqnarray*}
\mathcal{Q}\hspace{.1ex}(O,P) & \!\!=\!\! &\{ Q\in \mathcal{Q}(X\cup Y) \mid Q|_X =O, \ Q|_Y = P, \ Q\mathop{\cap}\, (Y\!\times\!X) = \emptyset\},\\
\mathcal{D}(O,P) & \!\!=\!\!  & \{ \hspace{.2ex}f:  Y\rightarrow \mathcal{D}(O) \mid f \mbox{ is isotone from } (\hspace{.3ex}Y,P) \mbox{ to } (\mathcal{D}(O),\subseteq) \},\\
\mathcal{U}(O,P) & \!\!=\!\!  & \{ \hspace{.2ex}g: X \rightarrow \mathcal{U}(P) \mid \,g \mbox{ is antitone from } (X,O) \mbox{ to } (\mathcal{U}(P),\subseteq) \},\\
\mathcal{J}(O,P) & \!\!=\!\!  & \{ \hspace{.2ex} h: \mathcal{D}(P) \rightarrow \mathcal{D}(O) \mid \,h \mbox{ preserves arbitrary joins (unions)} \},\\ 
\mathcal{R}(O,P) & \!\!=\!\!  & \{ R\subseteq X\!\times\!Y \mid OR\cup RP \subseteq R\,\}.
\end{eqnarray*}
Note that $\mathcal{Q}(O,P)$ is the collection of all {\em partial} orders extending $O+P$ and contained in $O \mathop{\oplus} P$. 
The posets with order relations in $\mathcal{Q}\hspace{.1ex}(O,P)$ are called {\em generalized vertical sums} of $(X,O)$ and $(Y,P)$. 

Being closed under arbitrary unions and intersections, the sets $\mathcal{R}(O,P)$ are Alexandroff topologies on $X\times Y$, hence superalgebraic lattices.

The {\em tensor product} $C\otimes D$ of two complete lattices $C$ and $D$ (Shmuely \cite{Sh}) may be represented as the complete lattice of all join-preseving maps from $C$ to the dual of $D$.
Notice that $C\otimes D$ is always isomorphic to $D\otimes C$. As usual, we order sets of functions into posets pointwise. 

\BT
\label{biext}
The following six superalgebraic lattices are isomorphic: 
\[\mathcal{Q}\hspace{.1ex}(O,P), \ \mathcal{D}(O,P), \ \mathcal{U}(O,P), \ \mathcal{J}(O,P), \ \mathcal{R}(O,P), \ \mathcal{U}(O)\mathop{\otimes} \mathcal{D}(P).\] 
\ET

\BP
Define $\iota : \mathcal{R}(O,P) \rightarrow \mathcal{Q}(O,P)$  by $\iota (R) = O\cup P \cup R$. This is a well-defined map: 
for $R \in \mathcal{R}(O,P)$ and $Q = \iota(R),$ we have $Q|_X = O$, $Q|_Y = P$, and $R\subseteq X\!\times\!Y$, hence
$Q \cap\,(Y\!\times\!X) = \emptyset$; furthermore, $Q$ is reflexive since $\Delta_{X \cup Y} = \Delta_X \cup \Delta_Y \subseteq Q$, and transitive since
$QQ = OO \mathop{\cup} PP \mathop{\cup} OR \mathop{\cup} RP$ is contained in $O\cup P \cup R = Q$; this proves $Q\in \mathcal{Q}(X\cup Y)$. 
By the equations $\iota(R)|_X \!=\! O$ and $\iota(R)|_Y \!=\! P$, the map $\iota$ is one-to-one. 

Now, let $Q$ be an arbitrary partial order in $\mathcal{Q}(O,P)$, and put $R = Q\mathop{\cap}\,(X\!\times\!Y)$.
Then $RP\subseteq QQ \cap (X\!\times\!Y) = Q \mathop{\cap}\,(X\!\times\!Y) = R$, and analogously $OR \subseteq R$. Thus, $R\in \mathcal{R}(O,P)$.
In order to show that $\iota$ is onto, it remains to verify that $Q = \iota (R)$. 
In fact, $Q = Q|_X\cup Q|_Y \cup (Q\cap (X\!\times\!Y)) = O\cup P \cup R = \iota (R)$, because $Q\cap (Y\!\times\!X)$ is empty. 
By definition, $R\subseteq S$ is equivalent to $\iota(R)\subseteq \iota(S)$, and consequently, $\iota$ is an isomorphism between the superalgebraic lattices $\mathcal{R}(O,P)$ and $\mathcal{Q}(O,P)$.

The map 
\BM
\alpha : \mathcal{U}(O,P) \rightarrow \mathcal{R}(O,P), \ g \mapsto \{ (x,y)\in X\!\times\!Y\mid y\in g(x)\}
\EM
is a well-defined isomorphism, too: for $R =\! \alpha (g)$ and $(x,z)\in RP$ there is a $y \in Y$ with $x \, R\, y \, P z$. 
Thus $y \in g(x) \in \mathcal{U}(P)$ and so $z \in g(x)$, i.e.\ $(x,z) \in R$. 
Hence, we have $RP \subseteq R$. For $OR \subseteq R$, note that $x\, O\, y\, R\, z$ implies $z\in g(y)\subseteq g(x)$, hence $x\,R\,z$. 
Surjectivity is straightforward: given $R\in \mathcal{R}(O,P)$, define $g: X \rightarrow \mathcal{U}(P)$ by $g(x) = xR \in \mathcal{U}(P)$ (since $xRP$ is contained in $xR$ by definition of $\mathcal{R}(O,P)$)
and verify that $g$ is an antitone map from $(X,O)$ to $(\mathcal{U}(P),\subseteq)$ with $\alpha(g) = R$. 
The proof of the isomorphism properties is completed by the equivalences 

$g\leq h \Leftrightarrow \forall \,x\in X\ (g(x)\subseteq h(x)) \Leftrightarrow \alpha (g) \subseteq \alpha (h)$.

\noindent Similarly, one sees that the map
\BM
\beta : \mathcal{D}(O,P) \rightarrow \mathcal{R}(O,P), \ f \mapsto \{ (x,y)\in X\!\times\!Y\mid x\in f(y)\}
\EM
is an isomorphism. 

By the universal property of $\mathcal{D}$ (see Section \ref{basic}), 
the lattice $\mathcal{D}(O,P)$ is isomorphic to the lattice $\mathcal{J}(O,P)$ of all join-preserving maps from $\mathcal{D}(P)$ to $\mathcal{D}(O)$, 
which in turn is  isomorphic to the tensor product $\mathcal{D}(P) \mathop{\otimes} \mathcal{U}(O)$ 
(because $\mathcal{D}(O)$ and $\mathcal{U}(O)$ are duals of each other), hence to $\mathcal{U}(O) \mathop{\otimes} \mathcal{D}(P)$.  
\EP

In the special case that $O$ is the antichain on $X$, we immediately deduce from Theorem \ref{biext} (denoting, as usual, by $\mathcal{P} (X)$ the power set of $X$): 

\BCO
\label{exta}
For disjoint sets $X,Y$ and any partial order $P$ on $Y$, the following five lattices are isomorphic to $\mathcal{P}(X)\otimes\mathcal{D}(P)$: 
\begin{eqnarray*}
\mathcal{Q}\hspace{.1ex}(X,P)\!\!& = &\!\!\{ Q\in \mathcal{Q}(X\cup Y) \mid Q|_X = \Delta_X, \, Q|_Y = P, \ Q \mathop{\cap}\, (Y\!\times\!X)=\emptyset\},\\
\mathcal{D}(X,P)\!\!& = &\!\!\{ \hspace{.2ex}f: Y \rightarrow \mathcal{P}(X) \mid x \,P\, y \, \Rightarrow \, f(x)\subseteq f(y)\},\\
\mathcal{U}(X,P)\!\!& = &\!\!\{ \hspace{.2ex}g \mid g : X \rightarrow \mathcal{U}(P) \} = \mathcal{U}(P)^{\,X},\\
\mathcal{J}(X,P)\!\!& = &\!\!\{ \hspace{.2ex}h: \mathcal{D}(P) \rightarrow \mathcal{P}(X) \mid \,h \mbox{ \rm preserves arbitrary unions} \},\\
\mathcal{R}(X,P)\!\!& = &\!\!\{ R\subseteq X\!\times\!Y \mid RP\subseteq R\}.
\end{eqnarray*}
\ECO

Observe that $\mathcal{Q}(X,P)$ consists of certain partial orders on $X\cup Y$ whose set of minimal elements contains but is not necessarily equal to the set\,\,$X$.
By definition, $\mathcal{D}(X,P)$ consists of all {\em isotone} maps from $(Y,P)$ to \mbox{$(\mathcal{P}(X),\subseteq)$,} while $\mathcal{U} (X,P)$ consists of {\em all} maps from $X$ to $\mathcal{U} (P)$. 

\newpage

Although we shall not need the following categorical aspect in our combinatorial computations, it is worthwhile to observe that one may form a new category $\mathcal{R}$ of posets by 
taking the sets $\mathcal{R}(O,P)$ as hom-sets, with the relation product $\circ$ as composition 
(this works because $OR\cup RP \subseteq R$ and $PS\cup SQ \subseteq S$ imply $ORS\cup RSQ \subseteq RS$, and $O\Delta_X \cup \Delta_X O = O$). 
In view of Theorem \ref{biext}, one might take alternatively the sets $\mathcal{Q}\hspace{.1ex}(O,P)$, $\mathcal{U}(O,P)$, $\mathcal{D} (O,P)$, or $\mathcal{J}(O,P)$ as hom-sets, with suitable compositions.

It is well-known that the superalgebraic lattices are, up to isomorphism, exactly the downset lattices of arbitrary posets, 
and that the category of posets and isotone maps is dual to the category of superalgebraic lattices and maps preserving arbitrary joins and meets (see Ern\'e \cite{EABC,EAO}). 
Now, let $\mathcal{S}_{\bigvee}$ denote the category of superalgebraic lattices with maps preserving arbitrary joins (but not necessarily meets) as morphisms. 
Within the extended framework of relations instead of functions, one shows by a series of straightforward verifications, involving Theorem \ref{biext}: 

\BT
\label{category}
The following facts hold without any finiteness restrictions:

{\rm (1)} The category $\mathcal{R}$ is self-dual under the contravariant functor sending any poset to its dual and each relation $R\in \mathcal{R}(O,P)$ to the dual $R^d$. 

{\rm (2)} The category $\mathcal{S}_{\bigvee}$ is self-dual under the contravariant functor sending any superalgebraic lattice to its dual 
and any join-preserving map \mbox{$h : B \rightarrow C$} between superalgebraic lattices to its upper adjoint $h^* : C \rightarrow B$,
given by 
\[\textstyle{h^*(c) = \bigvee\, \{ b \in B \mid h(b)\leq c\}.}\]

{\rm (3)} The upset operator $\mathcal{U}$, defined on morphisms $R\in \mathcal{R}(O,P)$ by
\[ \mathcal{U}(R): \mathcal{U}(O) \rightarrow  \mathcal{U}(P),\ \mathcal{U}(R)(U) = UR\]
gives rise to a functorial equivalence between the categories $\mathcal{R}$ and $\mathcal{S}_{\bigvee}$. 

{\rm (4)} Similarly, the downset operator $\mathcal {D}$ gives rise to a functorial duality between the categories $\mathcal{R}$ and $\mathcal{S}_{\bigvee}$. 
On morphisms $R\in \mathcal{R}(O,P)$, the contravariant functor $\mathcal{D}$ is defined by
\[\mathcal{D}(R): \mathcal{D}(P) \rightarrow  \mathcal{D}(O),\ \mathcal{D}(R)(D) = RD.\]  

\ET

\vspace{1ex}

Of course, the category of posets and isotone maps may be regarded as a subcategory of $\mathcal{R}$ (viewing functions as relations), 
and the category of superalgebraic lattices and complete homomorphisms as a subcategory of $\mathcal{S}_{\bigvee}$. 
However, these two dually equivalent subcategories fail to be self-dual!

\vspace{1ex}

In the finite case, Theorems \ref{biext} and \ref{category} give rise to several combinatorial `bijection equalities'. We shall return to that topic in Section \ref{exponential}.

\newpage

\section{\large Counting antichains, upsets and downsets}
\label{antichains}

In what follows, we deal with {\em finite posets}, denoted by $O$, $P$ etc.\ and represented in the usual way by diagrams (see Section \ref{basic} and the Appendix).
 
For a finite poset $P$, we denote by $d(P)$ the cardinality of $\mathcal{D}(P)$, the downset lattice or free distributive lattice over $P$. 
Observe that $d(P)$ is not only the number of all downsets and the number of all upsets, but also the number of all {\em antichains},
i.e., subsets $A$ of pairwise incomparable elements (such that the induced order on $A$ is the identity relation). 
This coincidence is caused by the bijection between antichains and downsets (or upsets) sending $A$ to $PA$ (or to $AP$).

A very flexible formula for the recursive enumeration of all downsets (hence of all anti\-chains) of a single finite poset is given by

\BT
\label{recdown}
For any subset $A$ of a finite poset $P$, the number $d(P)$ of downsets satisfies the inequality 
\[\textstyle{d(P) \leq \sum_{B\subseteq A} d(P - A \mathop{\updownarrow}_P\! B) \ \mbox{ with } A \mathop{\updownarrow}_P\! B =  (A \!\mathop{\smallsetminus}\! B)P\, \cup PB.}
\]
Equality holds whenever $A$ is an antichain. 
\ET

\BP
Clearly, $\mathcal{D}(P)$ is the disjoint union of the sets 
\vspace{-1ex}
\[
\mathcal{D}_{AB}(P) = \{ D\in \mathcal{D}(P)\mid D\cap A = B\},
\]
where $B$ runs through the subsets of $A$. Hence it suffices to show that 
\[
\textstyle{\varphi_{AB} : \mathcal{D}_{AB}(P) \rightarrow \mathcal{D}(P \!\mathop{-}\! A \mathop{\updownarrow}_P\! B), \ D \mapsto D \mathop{\smallsetminus} PB}
\]
is a well-defined injection, and a bijection in case $A$ is an antichain.

Any downset $D\in \mathcal{D}(P)$ satisfies the implications
\[
\textstyle{D \mathop{\cap} A = B \,\Rightarrow \, D \,\mathop{\cap}\, (A\!\smallsetminus\!B)P = \emptyset \, \Rightarrow \, (D\smallsetminus PB) \,\mathop{\cap}\, (A \mathop{\updownarrow}_P\! B) = \emptyset.}
\]
Consequently, $D\in \mathcal{D}_{AB}(P)$ implies $\varphi_{AB}(D)\in \mathcal{D}(P \!\mathop{-}\! A \mathop{\updownarrow}_P\! B)$, so that $\varphi_{AB}$ is well-defined,
and $D = \varphi_{AB}(D)\mathop{\cup}\,(PB\mathop{\smallsetminus}(A\!\smallsetminus\!B)P)$; indeed, from $B\subseteq D$ it follows that $PB\subseteq D$, and assuming $y\in D\mathop{\cap}\,(A\!\smallsetminus\!B)P$ 
we would find an $x\in (A\!\mathop{\smallsetminus}\! B)\mathop{\cap} Py$, which would lead to the contradiction $x\in D\cap A = B$.
This shows that $\varphi_{AB}$ is one-to-one. (In fact, it is an embedding of the distributive lattice $\mathcal{D}_{AB}(P)$ in the distributive lattice $\mathcal{D}(P \!\mathop{-}\! A \mathop{\updownarrow}_P\! B)$.) 

Now suppose that $A$ is an antichain of $P$. Given $D'\!\in \mathcal{D}(P \!\mathop{-}\! A \mathop{\updownarrow}_P\! B)$, put $D = D'\cup PB$. 
Then $D$ is a downset of $P$; indeed, for \mbox{$y\in D$ and $x \in Py$,} we either have $y \in PB$ and so $x \in PB \subseteq D$, 
or we have the implications 

\pagebreak

\noindent \mbox{$y \in D \smallsetminus PB \,\Rightarrow\, y\in D' \,\Rightarrow\, y \notin A \mathop{\updownarrow}_P\! B \,\Rightarrow\, y \notin (A \smallsetminus B )P\Rightarrow\, x \notin (A \smallsetminus B)P$}

\hspace{9.5ex} $\Rightarrow x \notin A \mathop{\updownarrow}_P\! B \mbox{ or } x\in PB  \ \Rightarrow \ x\in D'\cup PB = D$,

\vspace{.5ex}

\noindent since $D'\in \mathcal{D}(P - A \mathop{\updownarrow}_P\! B)$. Furthermore, $D$ cannot meet $A\!\smallsetminus\!B$:

\vspace{.5ex}

$D \cap ( A \smallsetminus B)= (D' \mathop\cup PB) \cap ( A \smallsetminus B )$ $ = PB \cap (A \smallsetminus B) = \emptyset$, 

\vspace{.5ex}

\noindent as $A$ is an antichain and $B \subseteq A$. On the other hand, \mbox{$B\subseteq PB\cap A \subseteq D\cap A$.} This proves the equation  $B = D\mathop{\cap} A$. 
By definition, \mbox{$\varphi_{AB}(D) = D \!\smallsetminus\! PB =\! D'$} (because $D'$ is disjoint from $PB$). 
Thus, in that case, $\varphi_{AB}$ is onto, hence bijective (and an isomorphism between $\mathcal{D}_{AB}(P)$
and $\mathcal{D}(P \!- \! A \mathop{\updownarrow}_P\! B)$). 
\EP

Since the minimal elements of a poset always form an antichain, we have the following consequence:

\BCO
\label{downpos}
For the set $A$ of all minimal elements of a finite poset $P$,
\[
\textstyle{d(P) = \sum_{B\subseteq A} d(P - ( B \cup (A\!\smallsetminus\!B)P)) = \sum_{C\subseteq A} d(P - ( A \cup CP)).} 
\]
\ECO

Already the case of a singleton $A$ is useful (and well-known, see \cite{E1,EO}):

\BCO
\label{downqos}
For any element $x$ of a finite poset $P$,
\[
d(P) = d(P- P\{ x\}) + d(P- \{ x\}P).
\]
In particular, if $x$ is a minimal element of $P$ then
\[
d(P) = d(P-\{x\}) + d(P- \{x\}P).
\]
\ECO

\vspace{-1ex}

A further counting formula for antichains is obtained by applying a result on clique numbers in graphs (Erd\"os and Ern\'e \cite{EE}, Vollert \cite{V}) 
to the incomparability graphs of posets. Let $\mathcal{A}(P)$ denote the set of antichains of the poset $P$.

\BT
\label{clique}
For any antichain $A$ of a finite poset $P$,
\[
\textstyle{d(P) = \sum_{B\in \mathcal{A}(P\!-\!A)} 2^{\,\sharp (A\smallsetminus (PB\mathop{\cup} BP))}.} 
\]
\ET

The formulas for downsets in sums of posets are almost obvious:

\BL
\label{sumdown}
For any two disjoint finite posets $O$ and $P$,
\[
d(O+P) = d(O)\cdot d(P), \ \ d(O\oplus P) = d(O)+d(P)-1.
\]
\EL

\noindent These equations allow to determine recursively in a straightforward manner the number of downsets in $N$-free posets. 
In particular, 
$d(A_k) \!=\! 2^k$ and $d(C_k) \!=\! k\!+\!1$.

\BX
\label{Ex21}
{\rm Typical orders that are not $N$-free are the zigzag orders:
 a {\em zigzag poset} with $m$ minimal elements has one of the following diagrams:
\vspace{-1ex}
\BC
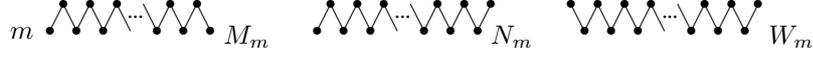
\begin{picture}(80,20)
\put(-15,-3){$m$}
\put(0,0){\circle*{3}}
\put(0,0){\line(1,2){5}}
\put(5,10){\circle*{3}}
\put(10,0){\line(-1,2){5}}
\put(10,0){\circle*{3}}
\put(10,0){\line(1,2){5}}
\put(15,10){\circle*{3}}
\put(20,0){\line(-1,2){5}}
\put(20,0){\circle*{3}}
\put(20,0){\line(1,2){5}}
\put(25,10){\circle*{3}}
\put(25,10){\line(1,-2){5}}
\put(30,5){\circle*{1}}
\put(32,5){\circle*{1}}
\put(34,5){\circle*{1}}

\put(40,0){\line(-1,2){5}}
\put(40,0){\circle*{3}}
\put(40,0){\line(1,2){5}}
\put(45,10){\circle*{3}}
\put(50,0){\line(-1,2){5}}
\put(50,0){\circle*{3}}
\put(50,0){\line(1,2){5}}
\put(55,10){\circle*{3}}
\put(60,0){\line(-1,2){5}}
\put(60,0){\circle*{3}}
\put(65,-5){$M_m$}
\end{picture} \ \ \ \ \
\begin{picture}(80,20)
\put(0,0){\circle*{3}}
\put(0,0){\line(1,2){5}}
\put(5,10){\circle*{3}}
\put(10,0){\line(-1,2){5}}
\put(10,0){\circle*{3}}
\put(10,0){\line(1,2){5}}
\put(15,10){\circle*{3}}
\put(20,0){\line(-1,2){5}}
\put(20,0){\circle*{3}}
\put(20,0){\line(1,2){5}}
\put(25,10){\circle*{3}}
\put(25,10){\line(1,-2){5}}
\put(30,5){\circle*{1}}
\put(32,5){\circle*{1}}
\put(34,5){\circle*{1}}

\put(40,0){\line(-1,2){5}}
\put(40,0){\circle*{3}}
\put(40,0){\line(1,2){5}}
\put(45,10){\circle*{3}}
\put(50,0){\line(-1,2){5}}
\put(50,0){\circle*{3}}
\put(50,0){\line(1,2){5}}
\put(55,10){\circle*{3}}
\put(60,0){\line(-1,2){5}}
\put(60,0){\circle*{3}}
\put(60,0){\line(1,2){5}}
\put(65,10){\circle*{3}}
\put(65,-5){$N_m$}
\end{picture} \ \ \ \ \
\begin{picture}(80,20)
\put(-5,10){\circle*{3}}
\put(0,0){\line(-1,2){5}}
\put(0,0){\circle*{3}}
\put(0,0){\line(1,2){5}}
\put(5,10){\circle*{3}}
\put(10,0){\line(-1,2){5}}
\put(10,0){\circle*{3}}
\put(10,0){\line(1,2){5}}
\put(15,10){\circle*{3}}
\put(20,0){\line(-1,2){5}}
\put(20,0){\circle*{3}}
\put(20,0){\line(1,2){5}}
\put(25,10){\circle*{3}}
\put(25,10){\line(1,-2){5}}
\put(30,5){\circle*{1}}
\put(32,5){\circle*{1}}
\put(34,5){\circle*{1}}

\put(40,0){\line(-1,2){5}}
\put(40,0){\circle*{3}}
\put(40,0){\line(1,2){5}}
\put(45,10){\circle*{3}}
\put(50,0){\line(-1,2){5}}
\put(50,0){\circle*{3}}
\put(50,0){\line(1,2){5}}
\put(55,10){\circle*{3}}
\put(60,0){\line(-1,2){5}}
\put(60,0){\circle*{3}}
\put(60,0){\line(1,2){5}}
\put(65,10){\circle*{3}}
\put(69,-5){$W_m$}
\end{picture}
\EC
\vspace{1ex}

For a zigzag poset with $k$ points, the number of downsets is the $k$-th Fibonacci number $f_k$ with $f_0 = 1$, $f_1= 2$ and $f_{k+2} = f_{k+1} + f_k$.
}
\EX

\vspace{1ex}

The formulas in Lemma \ref{sumdown} are extremal instances of one formula for the number of downsets in generalized vertical sums, as introduced in Section\,\ref{sums}. 
In slightly changed notation, such a vertical sum of two posets $O = (X,\leq_O)$ and $P = (Y,\leq_P)$ is of the form
\[O\oplus_R P = (X\cup Y, Q) \mbox{ with } Q = \ \leq_O \cup \,R \ \cup  \leq_P\vspace{-.5ex}\]
and $R\in \mathcal{R}(O,P)$, i.e., $R\subseteq X\times Y$ and $\leq_O\!\!R \,\mathop{\cup}\, R\!\leq_P \ \subseteq R.$ 

\BT
\label{vertdown}
The collection of all downsets in a generalized vertical sum is 
\[\mathcal{D}(O\oplus_R P) \ = \bigcup_{D\in \mathcal{D}(P)} \{ D\cup RD \cup E \mid E \in \mathcal{D}(O-RD)\}\vspace{-.5ex}\]
and this is a disjoint union. Hence, in the finite case,
\[\hspace{-20.5ex}d(O\oplus_R P) \ = \sum_{D\in \mathcal{D}(P)} d(O-RD).\]
\ET

\BP
By the definition of the quasiorder $Q$, it follows from $D\in \mathcal{D}(P)$ that $QD = D \mathop{\cup} RD$, hence $X \cap\, QD = RD $ and $X \smallsetminus QD = X \smallsetminus RD$. 
For any $D'$ in $\mathcal{D}(O \mathop{\oplus}_R P)$, we have $D = D' \cap Y \in \mathcal{D}(P)$,
$E = D' \cap (X\smallsetminus RD)\in \mathcal{D}(O-RD)$, and $QD \cup E = D'$, since $D'$ is  a downset containing $D$.

On the other hand, given any $D\in \mathcal{D}(P)$ and any $E\in \mathcal{D}(O-RD)$, we obtain $Q(QD \cup E) = QD \cup QE = QD \cup E$, since
$QE = OE$ is contained in $RD \cup E \subseteq QD \cup E$. Thus, $D \cup RD \cup E = QD \cup E \in \mathcal{D}(O\oplus_R P)$. 

The formula for $d(O\oplus_R P)$ follows from that for $\mathcal{D}(O\oplus_R P)$, because the union is disjoint (as $D = (QD \cup E)\cap Y$ for $D\in \mathcal{D}(P)$),
and the map $E \mapsto QD \cup E$ is one-to-one (as $E = (QD\cup E)\smallsetminus QD$ for $E\in \mathcal{D}(O-QD)$). 
\EP

 \vspace{-2ex}
 
\BCO
\label{downA}
For any finite poset $P$ disjoint from an antichain $A_m$, \vspace{-.5ex}
\[\hspace{-4ex}d(A_m \oplus_R P) = \sum_{D\in \mathcal{D}(P)} 2^{\,m-\sharp RD}.\vspace{-.5ex}\]
\ECO
This may also be regarded as a special instance of Theorem \ref{clique}.
Another application of Theorem \ref{vertdown} is the following `associative law':

\BCO
\label{OPP}
For partially ordered sets $O$, $P_1$ and $P_2$ on disjoint ground sets $X$, $Y_1$ and $Y_2$, respectively, one has \vspace{-.5ex}
\[R \in \mathcal{R}(O,P_1+P_2) \ \Leftrightarrow \ R_i \in \mathcal{R}(O,P_i) \mbox{ for } R_i =  R\cap (X\times Y_i), \ i = 1,2,\]
and \vspace{-1ex}
\[\hspace{-4ex}O \oplus_R (P_1 + P_2) \ = \ (O \oplus_{R_1}\! P_1) \oplus_{R_2}\! P_2.\]
Hence, in the finite case, \vspace{-.5ex}
\[\hspace{2.3ex}d(O\oplus_R\! (P_1 \!+\! P_2)) = \sum_{D\in \mathcal{D}(P_2)} d((O \oplus_{R_1}\!\! P_1) - R_2 D).\]
\ECO

\noindent Finally, a few additional arguments, similar to those for Corollary \ref{antitone} in the next section,  lead to the following extremality results (cf.\ Beinling \cite{Bei}):  

\BCO
\label{kn}
The maximal downset number of posets $P$ with $k$ points of which exactly $m$ are minimal (or maximal) is
\[d(P) = 2^{m-1} + 2^{k-1}. \] 
All posets satisfying that equation are isomorphic to $A_{m-1} + (A_1 \oplus A_{k-m})$.
There are exactly $\binom{k}{m}m$ such posets on $\underline{k}$ if $m < k$, and only  one if $m = k$.
\ECO

\BCO
\label{kh}
The maximal downset number of posets with $k$ points and height $h$ (the maximal size of chains in the poset) is 
\[d(P) = 2^{k-h}(h+1).\] 
All posets satisfying that equation are isomorphic to $A_{k-h} + C_h$.
There are exactly $(k)_h = \prod_{\,i = 0}^{h-1}(k \!-\! i)$ such posets on $\underline{k}$ if $h > 1$, and only  one if $h = 1$.
\ECO

\section{\large The exponential function of a finite poset}
\label{exponential}

From now on, let $M$ and $K$ always be disjoint finite sets of cardinality $m$ and\,\,$k$, respectively. For example, one may take $M = \underline{m}$, \mbox{$K = \underline{m\!+\!k}\smallsetminus \underline{m}$.} 
If $P$ is a poset with underlying set $K$, we write $\mathcal{Q} (M,P)$ for $\mathcal{Q}(M, \leq_P)$ etc. 
From Corollary \ref{exta}, we immediately deduce a few combinatorial identities:

\BCO
\label{ordmin}
For a finite poset $P$ on $K$, the following six distributive lattices have all the same cardinality $d(P)^m$:
$$\mathcal{Q} (M,P), \ \mathcal{D} (M,P), \ \mathcal{U}(M,P), \ \mathcal{R}(M,P), \ \mathcal{J} (M,P), \ \mathcal{P}(M)\otimes \mathcal{D}(P).$$
\ECO

We denote the set of minimal elements of a poset $P$ by ${\rm Min}\,P$ and its cardinality by $m(P)$. Furthermore, from now on, $\mathcal{Q}_0(K)$ will denote the set of all {\em posets} on $K$,
while $\mathcal{P}_0 (M)$ will denote the collection of all nonempty subsets of $M$.
Restricting the sets in Corollary \ref{ordmin}, we define:\\[-6ex]

\begin{eqnarray*}
\mathcal{E}\hspace{.1ex}(M,P)\!\!& = &\!\!\{ Q \in \mathcal{Q}_0(M\cup K) \mid \ {\rm Min}\, Q = M, \,Q|_K = P\},\\
\mathcal{F}(M,P)\!\!& = &\!\!\{ \hspace{.2ex}f: K \rightarrow \mathcal{P}_0(M) \mid x \leq_P y \Rightarrow \, f(x)\subseteq f(y)\},\\
\mathcal{G}(M,P)\!\!& = &\!\!\{ \hspace{.2ex} g : M \rightarrow \mathcal{U}(P) \mid \textstyle{\bigcup}\, g[M] = K \},\\
\mathcal{H}(M,P)\!\!& = &\!\!\{ \hspace{.2ex}h\!: \mathcal{D}(P) \rightarrow \mathcal{P}(M) \mid \,h \mbox{ \rm preserves unions}, h(D)\!=\! \emptyset \Rightarrow D \!=\! \emptyset \},\\ 
\mathcal{I}(M,P)\!\!& = &\!\!\{ R\subseteq M\!\times\!K \mid R\!\leq_P \ \subseteq R,\ MR = K\}.
\end{eqnarray*}

Thus, $\mathcal{E}\hspace{.1ex}(M,P)$ is the set of all posets on $M\cup K$ extending $P$ such that $M$ is exactly the set of all minimal elements. The numbers 
\[e(m,P) = \sharp \,\mathcal{E}\hspace{.1ex}(M,P)\] 
do not depend on the special choice of the set $M$ of cardinality $m$ with $M\cap K \!=\! \emptyset$. 
Now, analogous reasonings as for Theorem \ref{biext} and Corollary \ref{exta} lead to the next crucial result.

\BT
\label{EmP}
For finite disjoint sets $M, K$ and any poset $P$ on $K$,
\[
e(m,P) = \sharp \,\mathcal{E}(M,P) = \sharp \,\mathcal{F}(M,P) = \sharp \,\mathcal{G}(M,P) = \sharp \,\mathcal{H}(M,P) = \sharp \,\mathcal{I}(M,P). 
\]
Moreover, $e(m,P)$ is also the number of all $T_0$ topologies on $M\cup K$ inducing the topology $\mathcal{D}(P)$ on $K$ and making exactly the points in $M$ isolated.
\ET

\BCO
\label{antitone}
Regarded as functions, $\mathcal{G}$ and $e$ are strictly antitone in the second argument: for $M \neq \emptyset$ and posets $O,P$ on $K$, 

\noindent $\leq_O\, \subset\,\leq_P \ \Rightarrow \, \mathcal{U}(O) \supset\mathcal{U}(P) \, \Rightarrow \, \mathcal{G}(M,O) \supset\mathcal{G}(M,P) \, \Rightarrow \, e(m,O) > e(m,P).$
\ECO

\BP
For posets $O, P$ on $K$ such that $\leq_O$ is properly contained in $\leq_P$, it follows from \mbox{$U \!\in \mathcal{U}(P)$} that \mbox{$UO \ \subseteq UP \, \subseteq U$,} hence $U \in \mathcal{U}(O)$; 
and if $x \leq_P y$ but not $x \leq_O y$ then $xO$ lies in $\mathcal{U}(O)$ but not in $\mathcal{U}(P)$. 

\noindent The second implication is clear from the definition, and the third from the proven equation $e(m,P) = \sharp\,\mathcal{G}(M,P)$ in Theorem \ref{EmP}. 
\EP

We call $e(m,P)$ the {\em exponential function} of the poset $P$. 
This terminology is justified by the subsequent results, which provide some economic representations of such functions.

\BT
\label{expo}
\mbox{For the antichain\,$A$ of minimal elements of a finite poset\,$P$,}
\vspace{-2ex}
\begin{eqnarray} 
e(m,P) & \!\!=\!\! & \sum_{B\subseteq A} (-1)^{\sharp B}d(P\!-\! B)^m
\end{eqnarray}
\ET

\vspace{-1ex}

\BP
For the ground set $K$ of $P$ and any $g\in \mathcal{U}(P)^M$, the equation $K = \bigcup g[M]$ is tantamount to the inclusion 
$A\subseteq \bigcup g[M]$, since $K$ is the upset generated by $A$ and $\bigcup g[M]$ is an upset of $P$.
For each subset $B$ of\,\,$A$, the set \vspace{-1ex}
\[
\textstyle{\mathcal{G}_B = \{ g\in \mathcal{U}(P)^M \mid \bigcup g[M] \cap B = \emptyset\} = \mathcal{U}(P\!-\!B)^M }
\]
has cardinality $d(P\!-\!B)^m$.
Using the equations 
$$\textstyle{\mathcal{G}_B = \bigcap \,\{ \mathcal{G}_{\{ a\}} \mid a\in B\}, \ \mathcal{G}(M,P) = \mathcal{U}(P)^M\! \smallsetminus \bigcup \,\{ \mathcal{G}_{\{ a\}}\! \mid a\!\in\! A\},}$$

\noindent and a well-known combinatorial sieve inversion formula (see, e.g., \cite[II.1.17]{Ai}), 
one obtains from Theorem \ref{EmP}:

\vspace{1ex}

$e(m,P) = \sharp \,\mathcal{G}(M,P) = \sum_{B\subseteq A} (-1)^{\sharp B}d(P\!-\! B)^m.$
\EP

\vspace{1ex}

\noindent An informal algorithmic description of (1) was mentioned by Beinling \cite{Bei}.

Now, defining
\begin{eqnarray*}
b(i,j,P) & \!\!=\!\! & \sharp\, \{ B \subseteq {\rm Min}\,P \mid \sharp\, B = i, \, d(P\!-\!B) = j\},\\
c(j,P)  & \!\!=\!\! & \sum_{i = 0}^{m(P)} (-1)^i \,b(i,j,P),
\end{eqnarray*}
we arrive at a crucial formula:

\BT
\label{expo2}
For any finite poset $P$ and all natural numbers $m$,
\vspace{-1ex}
\begin{eqnarray}
e(m,P) & \!\! = \!\! & \sum_{j = 1}^{d(P)} c(j,P)\cdot j^m \ \mbox{ with } \ c(d(P),P) = 1 \vspace{-3ex}
\end{eqnarray}
In particular, 
\[\sum_{j = 1}^{d(P)} c(j,P) = e(0,P) = 0 \mbox{ for } P\neq \emptyset, \ \ \sum_{j = 1}^{d(P)} c(j,P)\cdot j = e(1,P) = 1.\]
\ET

A list of all exponential functions for posets with at most five points, correcting a few misprints in \cite{Bei}, is presented in the Appendix of this paper.
We pick two posets from that list and form their ordinal sum:

\BX
\label{Ex41}
 $P_4 = \ \begin{picture}(10,8)(0,0) 
    \put(0.5,4){\circle*{3}}
    \put(6.5,4){\circle*{3}}
    \end{picture}$,     
 $P_7 = \ \begin{picture}(10,8)(0,0) 
    \put(0.5,-1.7){\circle*{3}}
    \put(6.5,-1.7){\circle*{3}}
    \put(3.5,6.7){\circle*{3}}
    \put(.5,-2.7){\line(1,3){3.5}}
    \put(6.5,-2.7){\line(-1,3){3.5}}
    \end{picture}$,     
  $P_{12} = \ \begin{picture}(10,8)(0,0) 
    \put(-.2,3.1){\circle*{3}}
    \put(6.2,3.1){\circle*{3}}
    \put(3,10.8){\circle*{3}}
    \put(3,-4.8){\circle*{3}}
    \put(-.5,4.4){\line(2,-5){4}}
    \put(6.5,4.4){\line(-2,-5){4}}
    \put(-.5,2){\line(2,5){4}}
    \put(6.5,2){\line(-2,5){4}}
    \end{picture}$,
  $P_{47} = \ \begin{picture}(10,8)(0,-8)
    \put(-.2,-13.7){\circle*{3}}
    \put(6.2,-13.7){\circle*{3}}
    \put(-.2,-3.5){\circle*{3}}
    \put(6.2,-3.5){\circle*{3}}
    \put(0,-2.7){\line(1,-2){6}}
    \put(6,-2.7){\line(-1,-2){6}}
    \put(-.3,-3){\line(0,-1){10}}
    \put(6.2,-3){\line(0,-1){10}}
    \put(-.5,-2.7){\line(2,5){4}}
    \put(6.5,-2.7){\line(-2,5){4}}
    \put(3,6){\circle*{3}}
    \end{picture} \ = P_4 \oplus P_7\,,$\\   
    
$e(m,P_{47}) = d(P_{47})^m - 2\cdot d(P_{12})^m + d(P_7)^m = 8^m - 2\cdot 6^m + 5^m.$    
\EX

\noindent More involved is the computation of exponential functions for so-called fences:

\BX
\label{Ex42}
{\rm The {\em fence} on $\underline{2t} \!=\! \{ 1,...,2t\}$ is defined by the special zigzag order \vspace{-1ex}
\BM
x \,N_t\, y \ \Leftrightarrow x = y \mbox{ or } (x \mbox{ odd} \mbox{ and } |x\!-\!y| = 1).
\EM
A calculation based on Theorem \ref{expo} yields the following explicit formula:
\begin{eqnarray*}
e(m, N_{t}) & \!\!=\!\! &\sum_{s = t}^{2t} (-1)^s \sum_{i_1 +...+i_{2t-s+1} = s, \, 2|i_r \Leftrightarrow r =1}(\!\prod_{r=1}^{2t-s+1} f_{i_r})^m\\ 
            & \!\!=\!\! &\sum_{s = 0}^{t} (-1)^s \sum_{\!\!j_0 = \frac{1}{2},\, j_1 < \ldots < j_{s+1} = t+1,\, j_r\in \mathbb{N}} (\ \prod_{r=1}^{s+1} f_{2(j_r - j_{r-1}) - 1})^m,  
\end{eqnarray*}
where $f_i$ is the $i$-th Fibonacci number:

\vspace{1ex}

$\begin{array}{r|rrrrrrrr}
  i & 0 & 1 & 2 & 3 & 4 &  5 &  6 &  7\\
\hline
f_i & 1 & 2 & 3 & 5 & 8 & 13 & 21 & 34 
\end{array}
$
\vspace{2ex}
   
\noindent $N_1 =  
\begin{picture}(10,4)(-2,-2) 
    \put(-.2,6.3){\circle*{3}}
    \put(-.2,-4.5){\circle*{3}}
    \put(-.3,-5){\line(0,1){10}}
    \end{picture}
   , \hspace{2.5ex} e(m,N_1) = 3^m\! -2^m\! = {f_2}^m\! - (f_0 f_1)^m,$
    
\vspace{1ex}
    
\noindent $N_2 =  
\begin{picture}(10,8)(-2,-2) 
    \put(-.2,6.3){\circle*{3}}
    \put(6.2,6.3){\circle*{3}}
    \put(-.2,-4.5){\circle*{3}}
    \put(6.2,-4.5){\circle*{3}}
    \put(6,-5){\line(-1,2){6}}
    \put(-.2,-5){\line(0,1){10}}
    \put(6.2,-5){\line(0,1){10}}
    \end{picture}
    \ , \hspace{1.7ex} e(m,N_2)= 8^m\! - 6^m\! - 5^m\! + 4^m\!$ \\[1ex]
\indent $\hspace{18ex}= {f_4}^m\! - (f_2 f_1)^m\! -(f_0 f_3)^m\! + (f_0 f_1 f_1)^m,$

\vspace{1ex}

\noindent $N_3 =  
\begin{picture}(10,8)(-2,-2) 
    \put(-.2,6.3){\circle*{3}}
    \put(6.2,6.3){\circle*{3}}
    \put(-.2,-4.5){\circle*{3}}
    \put(6.2,-4.5){\circle*{3}}
    \put(6,-5){\line(-1,2){6}}
    \put(-.2,-5){\line(0,1){10}}
    \put(6.2,-5){\line(0,1){10}}
    \put(12.2,6.3){\circle*{3}}
    \put(12.2,-4.5){\circle*{3}}
    \put(12.2,-5){\line(0,1){10}}
    \put(12,-5){\line(-1,2){6}}
    \end{picture}
    \ \ , \,\ e(m,N_3) = 21^m\! - 16^m\! - 15^m\! - 13^m\! + 12^m\! + 2\cdot 10^m\! - 8^m$
    
    \vspace{1ex}

\noindent $ =  {f_6}^m\! - (f_4 f_1)^m - (f_2 f_3)^m - (f_0 f_5)^m + (f_2 f_1 f_1)^m +2 (f_0 f_1 f_3)^m - (f_0 f_1 f_1 f_1)^{m}.$     
    
}   
\EX

The exponential sum representation entails some divisibility properties of the numbers $e(m,P)-1$. 

\BCO
\label{Fermat}
For any finite poset $P$, natural number $n$, prime number $p$, and $m = n(p\!-\!1)+1$, the number $e(m,P)-1$ is divisible by $p$. 
In particular, $e(m,P)-1$ is always even, divisible by $\,6$ if $m$ is odd, by $\,30$ if $m \equiv 1 \!\!\mod\! 4$, by $\,42$ if $m \equiv 1 \!\!\mod\! 6$, and by $\,210$ if $m \equiv 1 \!\!\mod\! 12$.
\ECO 

\BP
Theorem \ref{expo2} and Fermat's Little Theorem give for $m = n(p\!-\!1)+1$:

$e(m,P) = \sum_{j = 1}^{d(P)} c(j,P)\cdot j \cdot ( j^n)^{p-1} \equiv \sum_{j = 1}^{d(P)} c(j,P)\cdot j = 1 \!\mod p.$ \vspace{1ex}

\noindent The remaining claims are obtained for the special cases $p = 2,3,5,7$.
\EP

For cardinal and ordinal sums of posets, the exponential functions are easily calculated from those of the summands, 
and the exponential behavior is also observed in the second argument:

\BT
\label{exposum}
For disjoint finite posets $O$ and $P$,
\begin{eqnarray}
e(m,O \mathop{+}P) &\!\!=\!\!& e(m,O)\cdot e(m,P)\\
e(m,O \mathop{\oplus} P) & \!\!=\!\! & \sum_{j=1}^{d(O)} c(j,O)\cdot (j+d(P)-1)^m 
\end{eqnarray}
provided $O$ is nonempty; hence, the exponential function of $\,O \mathop{\oplus} P$ is obtained from that of $O$ simply by increasing all base numbers by $d(P)-1$.
\ET

\BP
(3) The assignment $\ f \mapsto (f|_O, f|_P)$ gives a map from $\mathcal{F}(M,O\mathop{+} P)$ to $\mathcal{F}(M,O) \times \mathcal{F}(M,P)\,$ 
that is easily seen to be bijective. Now, apply Theorem \ref{EmP}.

(4) By the equations $\,{\rm Min} (O\mathop{\oplus} P) = {\rm Min}\,O$ \ and, for $B \subseteq {\rm Min}\,O$,

$\,O\mathop{\oplus} P - B = (O\!-\!B)\mathop{\oplus} P$, \ $d((O\!-\!B)\mathop{\oplus} P) = d(O\!-\!B)+d(P)-1$\,,

\noindent we get

$b(i,j,O \mathop{\oplus} P) = \sharp\,\{ B\subseteq {\rm Min}\, O\! \mid \sharp B = i, \ d(O\mathop{\oplus} P \mathop{-} B) =\! j\}$

$\hspace{13.5ex} = b(i,j\mathop{-}d(P)\mathop{+}\!1,O)$,\\[1ex]
and consequently $c(j,O\mathop{\oplus} P) = c(j\mathop{-}d(P)\mathop{+}1,O)$. 

Now, $d(O\oplus P - {\rm Min}\,O) \geq d(P)$ entails $c(j,O\mathop{\oplus} P) = 0$ for $j < d(P)$,\\[1ex]
and Theorem \ref{expo2} yields \vspace{1ex}

$e(m,O \mathop{\oplus} P) = \sum_{j = d(P)}^{d(O) \mathop{+} d(P)-1} c(j\mathop{-}d(P)\mathop{+}1,O) \cdot j^m$

\vspace{1ex}

$\hspace{13ex} = \sum_{j=1}^{d(O)} c(j,O)\cdot (j\mathop{+}d(P)\mathop{-}1)^m .$
\EP

\BCO
\label{exposum2}
For all finite posets $P$ disjoint from the antichain $A_{\ell}$,
\begin{eqnarray}
e(m,A_{\ell} \mathop{+}P) &\!\!=\!\!& (2^m-1)^\ell \cdot e(m,P)\\
e(m,A_{\ell} \mathop{\oplus} P) & \!\!=\!\! & \sum_{i=0}^{\ell} (-1)^{\ell -i} \binom{\ell}{i}\cdot (2^i+d(P)-1)^m\\[-1ex]
e(m,A_{\ell} \mathop{\oplus} A_{k}) & \!\!=\!\! & \sum_{i=0}^{\ell} (-1)^{\ell -i} \binom{\ell}{i}\cdot (2^i+2^k-1)^m
\end{eqnarray}
In particular, for a chain $C_k$ with $k$ elements not in $\underline{\ell}$,
\begin{eqnarray}
e(m,C_k) &\!\!=\!\!& (k+1)^m -k^m\\[2ex]
e(m,A_{\ell} \mathop{+} C_k) &\!\!=\!\!& (2^m-1)^\ell \cdot ((k+1)^m -k^m) \hspace{6ex}\\
e(m,A_{\ell} \mathop{\oplus} C_k) & \!\!=\!\! & \sum_{i=0}^{\ell } (-1)^{\ell -i} \binom{\ell}{i}\cdot (2^i+k)^m 
\end{eqnarray}
\ECO

\vspace{-1ex}

\BP
Equation (5) is obtained by an iterated application of (3). 
Equation (6) results from (4), observing that $c(2^i, A_{\ell}) = (-1)^{\ell - i}\binom{\ell}{i}$ for $i\leq \ell$, and $c(j,A_{\ell}) = 0$ otherwise. 
The formulas (7) and (8) are special cases of (6). Equation (9) is now immediate from (5) and (8), and (10) from (6).  
\EP

With the help of Theorem \ref{exposum} and Corollary \ref{exposum2}, one may determine recursively the exponential functions for all $N$-free posets. 
In particular, this gives all exponential functions for posets with at most four points except the fence $N = (\underline{4},N_2)$ (see Example \ref{Ex42}). 

\BX
\label{Ex43}
{\rm For posets $Q$ with a least element, i.e.\ $Q \simeq A_1\oplus P$, (6) gives 
\[e(m,Q) = d(Q)^m - (d(Q)\!-\!1)^m,\]
and for posets of the form $Q \simeq A_2\oplus P$, (6) amounts to 
\[e(m,Q) = (d(P)\!+\!3)^m -2\cdot (d(P)\!+\!1)^m + d(P)^m.\]
For instance, this confirms the formula for $P_{47} = P_4 \oplus P_7$ in Example \ref{Ex41}.
}
\EX

\BX
\label{Ex44}
{\rm For the poset $P = A_3\oplus C_2$ with diagram \begin{picture}(10,8)(-3,-5) 
    \put(-1,-13){\circle*{3}}
    \put(4,-13){\circle*{3}}
    \put(9,-13){\circle*{3}}
    \put(4,5){\circle*{3}}
    \put(-1,-13){\line(1,2){5}}
    \put(9.5,-14){\line(-1,2){5}}
    \put(4,-14){\line(0,1){20}}
    \put(4,-4){\circle*{3}}
    \end{picture}$ \ \ $, (10) gives

\vspace{2ex}
    
$e(m,P) = \sum_{i =0} ^3 \binom{3}{i} (-1)^{\,3-i} (2^i + 2)^m = 10^m -3\cdot 6^m + 3\cdot 4^m -3^m.$    
}
\EX

Some insight into the downset structure of a poset $P$ and its extensions $Q\in \mathcal{E}(M,P)$ is encoded in a sort of `characteristic polynomial' $p_Q$, defined by
\[\hspace{4ex}p_Q (z) = \sum_{D\in \mathcal{D}(P)} z^{\,\sharp (M\smallsetminus \,Q D)} = \sum_{i = 0}^m p_{\,Q,i} \, z^i, \vspace{-1ex}\]
where
\[\hspace{-3ex}p_{\,Q,i} = \sharp \,\{ D\in \mathcal{D}(P)\mid \sharp (M\smallsetminus QD) = i\},\]
in particular,
\[p_{Q,0} = p_Q(0) = \sharp\,\{ D\in \mathcal{D}(P) \mid M\subseteq \, QD\}.\]

\BL
\label{poly}
For any $Q\in \mathcal{E}(M,P)$, $p_Q$ is a normalized polynomial of degree $m$ with $p_Q(1) = d(P)$ and $p_Q(2) = d(Q)$.
Specifically, 
\begin{eqnarray*}
\mbox{for }Q\in \mathcal{E}(\underline{1},P), \ p_Q(z) & \!\!\! = \!\!\! & z + d(P) -1,\\[1ex]
\mbox{for }Q\in \mathcal{E}(\underline{2},P), \ p_Q(z) & \!\!\! = \!\!\! & z^2 + (d(Q)-d(P)-3)\,z -d(Q)+2d(P)+2,\\[1ex]
\mbox{for }Q\in \mathcal{E}(\underline{3},P), \ p_Q(z) & \!\!\! = \!\!\! & z^3 + (\tfrac{1}{2}d(Q)\, - \hspace{.4ex} d(P)+ \tfrac{1}{2}p_{Q,0}-3)\,z^2\\
&  &  \hspace{3ex}-\,(\tfrac{1}{2}d(Q) -\! 2d(P) +\tfrac{3}{2}p_{Q,0} - 2)\,z + p_{Q,0}.
\end{eqnarray*}
\EL

\noindent The second part is obtained by solving an easy system of linear equations.

\BCO
\label{poly2}
For $Q,Q'\!\in \mathcal{E}(M,P)$, \vspace{1ex}

$p_Q = p_{Q'}  \Leftrightarrow  d(Q) = d(Q') $ \ in case $m \leq 2$,

$p_Q = p_{Q'} \, \Leftrightarrow \, (d(Q)= d(Q')$\,and $p_{Q,0} = p_{Q',0})$ \ in case $m  = 3$.
\ECO

\noindent This is also an immediate consequence of the fact that the difference polynomial $p_Q - p_{Q'}$ has degree at most $m-1$ and the roots $1$ and $2$, whence \vspace{-1ex}
\[p_Q(z) - p_{Q'}(z) = (z-1)(z-2)\, r(z) = (z^2-3z +2)\sum_{i = 0}^{m-3} r_i\, z^i\]
for an integer polynomial $r$ of degree at most $m-3$; in particular, $r$ is an integer constant if $m\leq 3$ and equal to $0$ if $m\leq 2$. 
The above factorization provides the difference between $p_Q$ and $p_{Q'}$ coefficientwise: 
\[p_{\,Q,i} - p_{\,Q',i} = r_{i-2} - 3\, r_{i-1} +2\,r_i \ \mbox{ for } 0\leq i < m \, , \mbox{ and } r_i = 0 \mbox{ for all other } i.\]

\BX
\label{Ex45}
{\rm We sketch the three posets $P_n$ with $k_n\! = \! 5$ points, $m_n\! =\! 3$ minimal points, and $d_n\! =\! 14$ downsets (see the Appendix).
One of them has height $h_n \! = \! 3$ and belongs to $\mathcal{E}(\underline{3},C_2)$, while the other two have height $h_n \! = \! 2$ and belong to $\mathcal{E}(\underline{3},A_2)$.

\begin{picture}(400,54)
\put(30,10)
{   \begin{picture}(10,8)
    \put(-1,15){\circle{3}}
    \put(4,15){\circle{3}}
    \put(9,15){\circle{3}}
    \put(4,36){\circle*{3}}
    \put(9,16){\line(-1,4){5}}
    \put(4,16){\line(0,1){20}}
    \put(4,26){\circle*{3}}
    \put(20,30){$P = C_2$}
    \put(20,16){$Q = P_{73}$}
    \put(-19,-2){$p_Q(z) = z^3 +z^2+z$}
    \end{picture}
}
\put(130,10)
{   \begin{picture}(10,8)
    \put(-1,30){\circle*{3}}
    \put(4,19){\circle{3}}
    \put(9,30){\circle*{3}}
    \put(-1,19){\circle{3}}
    \put(-1,30){\line(1,-1){10}}
    \put(9,30){\line(-1,-1){10}}
    \put(-1,30){\line(0,-1){10}}
    \put(9,30){\line(0,-1){10}}
    \put(9,19){\circle{3}}
    \put(20,30){$P = A_2$}
    \put(20,16){$Q = P_{78}$}
    \put(-8,-2){$p_Q(z) = z^3 + 3z$}
    \end{picture}
}
\put(230,10)
{   \begin{picture}(10,8) 
    \put(-1,30){\circle*{3}}
    \put(4,19){\circle{3}}
    \put(9,30){\circle*{3}}
    \put(-1,19){\circle{3}}
    \put(-1,30){\line(1,-1){10}}
    \put(-1,30){\line(0,-1){10}}
    \put(9,30){\line(0,-1){10}}
    \put(4,20){\line(-1,2){5}}
    \put(9,19){\circle{3}}
    \put(20,30){$P = A_2$}
    \put(20,16){$Q = P_{79}$}
    \put(-8,-2){$p_Q(z) = z^3 + z^2 +2$}
    \end{picture}
}
\end{picture}
The difference between the last two polynomials is just $z^2 -3z +2$.}
\EX

\section{\large Recursions and matrix equations}
\label{rec}

Henceforth, given a subset $U$ of a poset $P$, we write $e(m,U)$ instead of $e(m,P|_U)$, regarding $U$ as a subposet of $P$. 
We now are going to establish two recursion formulas for the exponential functions of posets $P$ by representing $e(m,P)$ (resp.\ $e(m\!+\!1,P)$)
as integer linear combinations of the functions $e(m,U)$, where $U$ runs through all proper (resp.\ all) upsets of $P$.

\BT
\label{recurs1}
For all finite posets $P$ and natural numbers $m$,
\[\hspace{13ex}d(P)^m \, = \sum_{D\in \mathcal{D}(P)} e(m,P\!-\!D) = \ \ \sum_{U\in \,\mathcal{U}(P)} \, e(m,U) \ ; \mbox{ hence }\hspace{8ex}\]
\vspace{-3ex}
\begin{eqnarray}r(m,P) & \!\!\! = \!\!\! & d(P)^m - e(m,P) \, = \!\!\sum_{U\in\, \mathcal{U}(P)\smallsetminus \{ P\}} \hspace{-2ex} e(m,U) \ \
\end{eqnarray}
\ET

\BP
For each $U \in \mathcal{U}(P)$ let $\iota_U$ be the inclusion map from $\mathcal{U}(P|_U)$ into $\mathcal{U}(P)$. 
Since for $g\in \mathcal{U}(P)^M$ the union $\bigcup g[M]$ is an upset of $P$, we have 
\[ \mathcal{U}(P)^M \, \ \subseteq \bigcup_{U \in \,\mathcal{U}(P)} \{ \iota_U \circ g \mid g \in \mathcal{G}(M,P|_U)\}. \]
Since the reverse inclusion is trivial and the union on the right is disjoint, we get the result by invoking Theorem \ref{EmP} once more. 
\EP

\BX
\label{Ex51}
{\rm For the fence $N = P_{19}$ and its eight upsets, isomorphic to $P_1$, $P_2$ (twice), $P_3$, $P_4$, $P_6$, $P_8$, and $P_{19}$ (see the Appendix), one obtains:}

\vspace{1ex}

\noindent $8^m   = e(m,P_1) + 2 e(m,P_2) + e(m,P_3) + e(m,P_4) + e(m,P_6)$

$ + \ e(m,P_8) + e(m,P_{19})$ 

$ = 1  +  (2\cdot 2^m - 2)  +  (3^m  -  2^m)  +  (4^m  -  2\cdot 2^m  + 1)  +  (5^m  -  4^m)$

$ + \ (6^m  -  4^m - 3^m  +  2^m)  +  (8^m  -  6^m  -  5^m  +  4^m).$  
\EX

\noindent In order to proceed from $m$ to $m\!+\!1$, consider the sets
\[\mathcal{F}\,'(M,P) = \{ (D,U,f)\mid D\in \mathcal{D}(P), \, U\in \mathcal{U}(P \!\mathop{-}\! DP), f\in \mathcal{F}(M,P \!-\! D)\}.
\]
A careful distinction of cases (nine in all) shows:

\BL
\label{Fprime}
For all $(D,U,f)\in \mathcal{F}'(M,P)$ and elements $i$, the function  

$f_i : P \rightarrow \mathcal{P}_0 (M\cup\{ i\})$ with
\[ f_i (x) = \left\{ 
\begin{array}{ll}
\{ i\} & \mbox{if } x\in D\\
f(x) \cup \{ i\} & \mbox{if } x\in (U\cup DP)\smallsetminus\! D\\
f(x) & \mbox{if } x\in K\!\smallsetminus\!(U\cup DP)
\end{array}
\right.
\]
is isotone, hence a member of $\mathcal{F}(M \cup \{ i\},P)$ if $i \notin K$.
\EL

\noindent Recall that for every subset $A$ of a poset $P$ on the ground set $K$, the set 
$A^{\circ} = K\!\smallsetminus (K\!\smallsetminus\!A)P$ 
is the interior in the lower Alexandroff topology $\mathcal{D}(P)$.

\BT
\label{recurs2}
For every finite poset $P$ disjoint from the set $\underline{m\! +\! 1}$, the assignment 
$(D,U,f) \mapsto f_{m+1}$ gives a bijection $\varphi$ between $\mathcal{F}\,' (\underline{m},P)$ and $\mathcal{F}(\underline{m\!+\!1},P)$. Hence 
\begin{eqnarray} \label{recurs2_formula}
e(m \!+\!1,P)\! & \!\! \! = \!\!\! & \!\!\!\! \sum_{D\in \mathcal{D}(P)} \!\! d(P\!-DP)\, e(m,P\!-\!D) = \!\!\! \sum_{U\in \,\mathcal{U}(P)} \! d(U^{\circ})\,e(m,U)\ \ \ \
\end{eqnarray}
\ET

\BP
Let us verify that the function below is well-defined and\,inverse\,to\,\,$\varphi$:

\vspace{1ex}

$\psi : \mathcal{F}(\underline{m\!+\!1},P) \rightarrow \mathcal{F}\,' (\underline{m},P), \ f \mapsto (D_f,U_f,f')$ \ with

\vspace{1ex}

$D_f = \{ x\in K \mid f(x) = \{ m\!+\!1\}\},$

$\,U_f = \{ x\in K \mid m\!+\!1 \in f(x)\} \!\smallsetminus\! D_f P,$

$f' : K\!\smallsetminus\!D_f \rightarrow \mathcal{P}_0(\underline{m}),\ x \mapsto f(x)\!\smallsetminus\!\{ m\!+\!1\}.$\\[1ex]
(1) $\psi$ is well-defined: as $f$ is isotone and each $f(x)$ is nonempty, $D_f$ lies in $\mathcal{D}(P)$, $U_f$ in $\mathcal{U}(P\!\mathop{-}\! D_fP)$, and $f'$ is isotone, too.
Furthermore, $f'(x)$ belongs to $\mathcal{P}_0(\underline{m})$ for each $x\in K \!\smallsetminus\! D_f$, since the set $f(x)$ contains at least one element distinct from $m\!+\!1$.\\[1ex]
(2) $(\psi\circ \varphi) (D,U,f) = (D,U,f)$: in fact, $\{ x \in K \mid f_{m+1}(x) = \{ m\!+\!1\}\} = D$, 
$\{ x\in K \mid m\!+\!1 \in f_{m+1}(x)\} \smallsetminus DP = U$, and $f_{m+1}(x)\smallsetminus \{ m\!+\!1\} = f(x)$.\\[1ex]
(3) $(\varphi\circ \psi)(f) = f$ for $f\in \mathcal{F}(\underline{m\!+\!1},P)$: in fact, $\psi (f) = (D_f,U_f, f')$, and by definition,
$f'_{m+1} = f$.\\[1ex]
The recursion formula for $e(m\!+\!1,P)$ is now clear, again by Theorem \ref{EmP}.
\EP

We find it instructive to express the explicit and recursive formulas found for $e(m,P)$ in terms of certain matrix equations. 
Let $(P_n \mid n \in \mathbb{N})$ be an enumeration of all unlabeled (i.e.\ isomorphism types of) finite posets  
such that $m \leq n$ implies $k_m \leq k_n$, where $\underline{k_n}$ is the ground set of $P_n$ (see the Appendix). 
Then we have the following (infinite) matrices:
\begin{eqnarray*}
A = (a_{mn}), \ a_{mn}\! & \!\!\! = \! & \textstyle{\sum_{U\in \,\mathcal{U}(P_n), \, U \simeq P_m} \! d(U^{\circ})},\\[.5ex]
B\hspace{.2ex} = (b_{mn}), \ b_{mn} & \!\!\! = \! & \sharp \,\{ U \! \in \mathcal{U}(P_n) \mid U\simeq P_m\},\\[.5ex]
C = (c_{mn}), \ c_{mn} & \!\!\! = \! & \delta_{mn} -  \textstyle{\sum_{j = m+1}^n b_{mj}c_{jn}},\\[.5ex]
D = (d_{mn}), \ d_{mn} \!& \!\!\! = \! & d(P_n)^m,\\[.5ex]
E\, = (e_{mn}), \ e_{mn} \!& \!\!\! = \! & e(m,P_n).
\end{eqnarray*}
The columns of the matrix $E$ represent the exponential functions $e(m,P_n)$ for fixed $n$, while the rows $E_m$ list the values $e(m,P_n)$ for fixed $m$.
Using the identity matrix $I = (\delta_{mn})$ and its upper triangular shift $\vec{I} = (\delta_{m,n-1})$, we now are able to prove some nice matrix equations.

\BT
\label{matrix}
$A$, $B$, and $C$ are upper triangular matrices with
\vspace{-1ex}
\begin{eqnarray*}
(1) \ \ BC & \!\! = \!\! & \,I\ , \ \mbox{ hence } C = B^{-1},\\
(2) \ \ EB & \!\! = \!\! & D\,, \ \mbox{ hence } DC = E,\\
(3) \ \ EA & \!\! = \!\! & \vec{I}E\hspace{.1ex},  \mbox{ hence } E_m A = E_{m+1}.
\end{eqnarray*}
\ET

\BP
For $m > n$, no $U\in \mathcal{U}(P_n)$ can be isomorphic to $P_m$, because that would imply that $P_m$ is embeddable in $P_n$, and then either $k_m < k_n$ (hence $m < n$)
or $\sharp\, U = k_m = k_n$ (hence $P_m \simeq U = P_n$, and so $m =n$). This shows that $A$ and $B$ are upper triangular, and from the recursive definition of 
the coefficients of $C$ it follows that $C$ is upper triangular, too, and inverse to $B$, proving (1).\\[1ex]
(2) By Theorem \ref{recurs1}, $d_{mn} = d(P_n)^m = \sum_{U\in \,\mathcal{U}(P_n)} e(m,U) = \sum e_{mj} b_{jn}$, whence $D = EB$, and by (1), $DC = E$.\\[1ex]
(3) By Theorem \ref{recurs2}, $(E_m A)_n$ is equal to

\noindent $\sum_{j : U\in \mathcal{U}(P_j),U\simeq P_m} \!e(m,P_j)\,d(P_m\!-(P_m\!-\!U)P_m) = e(m\!+\!1,P_n) = (E_{m+1})_n$.
\EP

\section{\large Estimates and asymptotical results} 

As expected, the largest exponential functions are obtained for antichains.

\BT
\label{est1}
For all posets $P$ with $k$ points, except the antichain $A_k$, and for all natural numbers $m > 1$, one has the inequality 
\[e(m,P) < (\tfrac{3}{4})^m \, (2^m -1)^k = (\tfrac{3}{4})^m \, e(m,A_{k}).\]
\ET

\BP
$P$ contains two comparable elements. Hence, by Corollary \ref{antitone}, we obtain $e(m,P)\leq e(m, A_{k-2} + C_2)$, and so by 
(9),\\[1ex]
$e(m,P) \leq \frac{12^m-8^m}{4^m}\,(2^m-1)^{k-2} \! < \! \frac{12^m -2\cdot 6^m +3^m}{4^m} \,(2^m - 1)^{k-2} = (\frac{3}{4})^m (2^m-1)^k$\\[1ex]
since $2\cdot 6^2 - 3^2 = 63 < 64 = 8^2$ and $2\cdot 6^m \leq 6^m +2m\cdot 6^{m-1} < (6+2)^m = 8^m$ for $m>2$.
\EP

Compare Theorem \ref{est1} with the well-known inequality 
$$\textstyle{d(P) \leq \frac{3}{4}\,2^k = \frac{3}{4}\,d(A_k) \mbox{ for } P \neq A_k.}$$ 

The leading term in an exponential sum is very dominant and makes asymptotical estimates rather easy. 

\BT
\label{est2}
For any poset $P$ on $K$, the function $e(m,P)$ is asymptotically equal to $d(P)^m$. More precisely, for \vspace{1ex}

$d'(P) = \max \{ d(P|_A) \mid A \subset K \}$, one has $d'(P) \leq 2^{k-1}$ and 
\begin{eqnarray}
0 \leq r(m,P) & \!\!\! = \!\!\! & d(P)^m - e(m,P) \leq 2^{\,m(P)-1} d'(P)^m \leq 2^{\,(k-1)(m+1)} \ \ \label{est2_formula} \\[1ex] 
e(m,P) & \!\!\! = \!\!\! & d(P)^m (1 - O(\varepsilon(P)^m)) \ \mbox{ with } \  \varepsilon (P) = \tfrac{d'(P)}{d(P)} < 1 
\end{eqnarray}
\ET

\BP
By (11), we have $0 \leq r(m,P) = d(P)^m - e(m,P)$, and by (1), 
$$
\textstyle{r(m,P) \leq \sum_{B \subseteq {\rm Min}\,P, \, \sharp B\, odd \,} d(P-B)^m \leq 2^{\,m(P)-1} d'(P)^{m}.} 
$$

\noindent The asymptotical equation (14) is an immediate consequence of (13). 
\EP

\vspace{1ex}

\noindent The definition of the coefficients used in Theorem \ref{expo2} gives the inequality
\vspace{-1ex}
\[\sum_{j=1}^{d(P)} |c(j,P)| \leq 2^{m(P)}\]
and one might conjecture that even equality holds; indeed, that is true very frequently, but not always (because certain summands $(-1)^i b(i,j,P)$ may absorb each other). 
For posets with up to five points, there are only two exceptions (see the Appendix):

\BX
\label{Ex61}
{\rm For the posets with the numbers $n \!=\! 71$ and $n \!=\! 81$ one has

\vspace{1ex}
$P_{71} =$ 
\begin{picture}(10,10)(-3,-5) 
    \put(-1,-13){\circle*{3}}
    \put(4,-13){\circle*{3}}
    \put(9,-13){\circle*{3}}
    \put(4,5){\circle*{3}}
    \put(9,-14){\line(-1,4){5}}
    \put(4,-14){\line(1,4){2.5}}
    \put(6.5,-4){\circle*{3}}
    \end{picture} \ , \ $e(m,P_{71}) =12^m - 2\cdot 8^m + 2 \cdot 4^m - 3^m$,
    
\vspace{2ex}    
$P_{81} = $
\begin{picture}(10,8)(-3,-5) 
    \put(4,0){\circle*{3}}
    \put(4,-10){\circle*{3}}
    \put(9,0){\circle*{3}}
    \put(-1,-10){\circle*{3}}
    \put(4,0){\line(1,-2){5}}
    \put(4,0){\line(0,-1){10}}
    \put(9,0){\line(0,-1){10}}
    \put(9,-10){\circle*{3}}
    \end{picture}\ \ , \  $e(m,P_{81}) = 16^m - 12^m - 10^m + 6^m + 5^m - 4^m$,
    
\vspace{2ex}

$m(P_n) = 3\,, \ \sum_{j=1}^{d(P_n)} |c(j,P_n)| = 6 < 2^3.$ 
}
\EX

\vspace{1ex}

\BT
\label{est3}
$e(m \!+\!1,P)$ is asymptotically equal to $d(P)\, e(m,P)$. \\[1ex] 
More precisely, for $\varepsilon (P) = \frac{d'(P)}{d(P)} \leq 1 -\frac{1}{d(P)}$ and $m \geq (m(P)\!+\!1)\, d(P) \,{\rm ln}\,2$: \vspace{-1ex}
\[
1 \leq \frac{e(m\!+\!1,P)}{d(P)\, e(m,P)}  \leq  1+ 2^{\,m(P)}\,\varepsilon (P)^m , \ \frac{e(m\!+\!1,P)}{d(P)\,e(m,P)} = 1+ O(\varepsilon (P)^m). 
\]
\ET

\BP 
For $ m \geq - \dfrac{(m(P)+1) \ln 2}{\ln \,\varepsilon (P)}$, that is, $2^{m(P)}\,\varepsilon (P)^m \!\leq 2^{-1}$, we get
\vspace{-1ex}
\begin{eqnarray*}
\,0 \!\!& \!\!\!\stackrel{(\ref{recurs2_formula})}{\leq}\!\!\! &\!\! \frac{e(m\!+\!1,P)}{d(P)\, e(m,P)} - 1 = \frac{d(P)^{m+1} -r(m\!+\!1,P) -d(P)^{m+1} +d(P)\,r(m,P)}{d(P)\,e(m,P)}\\
& \!\!\!\leq\!\!\! &\frac{r(m,P)}{e(m,P)} = \frac{r(m,P)}{d(P)^m-r(m,P)} \stackrel{ (\ref{est2_formula})}{\leq} 
\frac{2^{\,m(P)-1}\,\varepsilon (P)^m} {1 - 2^{m(P) - 1 }\,\varepsilon (P)^m} \leq 2^{\,m(P)}\,\varepsilon (P)^m, 
\end{eqnarray*}
and using the straightforward inequality 
\vspace{-2ex}
\[ d(P) \geq - \frac{1}{\ln (1-1/d(P))} \geq - \frac{1}{\ln \varepsilon (P)} \] 
we arrive at the claimed estimate in Theorem \ref{est3}.
\EP

\noindent This estimate is quite good for large $m$, but not for small $m$, because then $\varepsilon (P)^m$ may be close to $1$. 
For posets of size $k$ much smaller than $m$, the following estimates may sometimes be more effective:

\BL
\label{est4k}
For any nonempty downset $D$ of a finite poset $P$,
\[m\cdot e(m,P\!\mathop{-}\!D) \leq  2^{\,\sharp (DP \smallsetminus D)}\, e(m,P).
\]
\EL

\BP
For each $f\in \mathcal{F}(M,P\!\mathop{-}\!D)$ and each $i\in M$, the map $f_i : K \rightarrow \mathcal{P}_0 (M)$ with
\[ f_i (x) = \left\{ 
\begin{array}{ll}
\{ i\} & \mbox{if } x\in D\\
f(x) \cup \{ i\} & \mbox{if } x\in DP \!\smallsetminus\! D\\
f(x) & \mbox{if } x\in K\!\smallsetminus\! DP
\end{array}
\right.
\]
is a member of $\mathcal{F}(M,P)$ (Lemma \ref{Fprime} with $U = \emptyset$). 
For any $g\in \mathcal{F}(M,P\!\mathop{-}\!D)$, the equation $f_i = g_j$ forces $i = j$ (because the set $D$ is nonempty) and then 
$f(x) \subseteq f_i(x) = f(x) \cup \{ i\} = g_i(x) = g(x) \cup \{ i\} $ for all $x\in DP\smallsetminus D$. 

In all, there are at most $2^{\,\sharp (DP\smallsetminus D)}$ possibilities for $g$ so that $f_i = g_j$. Now, Theorem \ref{EmP} yields
$e(m,P\!\mathop{-}\!D) = \sharp \, \mathcal{F}(M,P\!-\!D)$ and so \vspace{1ex}

$m\cdot e(m,P\!\mathop{-}\!D) \leq 2^{\,\sharp (DP\smallsetminus D)} \sharp\, \mathcal{F}(M,P) = 2^{\,\sharp (DP\smallsetminus D)} e(m,P).$ 
\EP

\BT
\label{est4}
For all finite posets $P$ on $k$ points and natural numbers $m$, 
\[1 \leq \frac{e(m\!+\!1,P)}{d(P)\,e(m,P)} < \, 1+ m^{-1} \,2^k \,d(P) \,\leq\, 1 + m^{-1}\, 4^k.\] 
\ET

\BP
By Theorem \ref{recurs2}, we have
$$\frac{e(m\!+\!1,P)}{d(P)\,e(m,P)}- 1 \!= \hspace{-2ex}\sum_{\emptyset \neq D\in \mathcal{D}(P)} \hspace{-3ex} \frac{d(P \!\mathop{-}\! DP)\,e(m,P\!-\!D)}{d(P)\,e(m,P)} 
\! < d(P) \max_{_{\emptyset \neq D\in \mathcal{D}(P)}} \hspace{-1ex}\frac{e(m,P\!-\!D)}{e(m,P)}.$$ 
Now, Lemma \ref{est4k} and the inequality $\sharp (DP \smallsetminus\!D) < k$ conclude the proof. 
\EP

\vspace{1ex}

As in Section \ref{rec}, let $(P_n : n \in \mathbb{N})$ be an enumeration of all unlabeled finite posets, 
arranged in such a manner that $P_n$ is a poset on $\underline{k_n}$ and $m \leq n$ implies $k_m \leq k_n$. 
For any finite poset $P$ with $k$ points, let $a(P)$ denote the number of all automorphisms (symmetries) of $P$. Then
$i(P) = k!/a(P)$ is the number of isomorphic copies of $P$ on the same ground set, and 
\begin{equation} \label{ekm_als_summe}
e_k(m) = \sum_{P\in \mathcal{Q}_0 (\underline{k})} e(m,P) = \sum_{n \mbox{ \footnotesize with } k_n = k} i(P_n)\,e(m,P_n) 
\end{equation}
is the number of all ({\em labeled}) partial orders or posets on $\underline{m\!+\!k}$ such that $\underline{m}$
is the set of all minimal elements. From this description it is clear that 
\[e_k(1) = p(k) = \sharp\, \mathcal{Q}_0 (k)\]
but it is a surprising and non-trivial fact, established in \cite{E3}, that also
\[e_{k-1}(2) = p(k).\] 

The asymptotical behavior of the exponential sums $e_k(m)$ for fixed $k$ is easily determined with the help of an old result due to Stanley \cite{St2}: putting 
\[z_2 = 2, z_3 = 6, z_4 =20, \mbox{ and }z_i = 2 i \mbox{ for } i\geq 5,
\] 
there are $z_i\binom{k}{i}$ $T_0$ topologies (or partial orders) on $k$ points with $2^{k-1}+2^{k-i}$ open sets (resp.\ downsets), and no others with more than $2^{k-1}$ open sets.  
Now, Theorems \ref{est1} and \ref{est2} together with $(2^m \! - 1)^k = 2^{km} + O(2^{(k-1)m})$ lead to the following result:

\BT
\label{asym}
$e_k(m)$ is asymptotically equal to $2^{km}$. More precisely,
\vspace{-1ex} 
\[e_k(m) = 2^{km} + \sum_{i=2}^k z_i\binom{k}{i} (2^{k-1} + 2^{k-i})^m + O(2^{(k-1)m}) = 2^{km} (1+ O((\tfrac{3}{4})^m)). 
\]
\ET  

The leading difference has an interesting number-theoretical property:

\BT
\label{Cauchy}
$e_k(m)-2^{mk} - (-1)^k$ is divisible by $k$ if $k$ is prime.
\ET

\BP By (5) and (\ref{ekm_als_summe}),
$e_k(m)-(2^m-1)^k = \sum_{P_n \in \mathcal{Q}_0(\underline{k})\setminus \{ A_k\}} \frac{k!}{a(P_n)}\,e(m,P_n)$. \\[1ex] 
But if $k_n = k$ is prime, it cannot divide $a(P_n)$ unless $P_n = A_k$;  
indeed, if $k$ divides $a(P_n)$ then, by Cauchy's Theorem, ${\rm Aut}\, P_n$ must have a cyclic subgroup of order $k$, 
and consequently all points must have the same height, i.e., $P_n$ must be an antichain. Hence, $k$ must divide the sum on the right.
The congruence $(2^m -1)^k \equiv 2^{mk} + (-1)^k \! \mod k$ completes the proof.
\EP

Finally, note a generalization of Corollary \ref{kn}, based on Corollary \ref{downA}:

\BT
For any $Q\in \mathcal{E}(M,P)$,
\[d(Q) = \sum_{D\in \mathcal{D}(P)} 2^{\sharp(M\smallsetminus QD)} \leq 2^{\,m-1}(d(P)+1),
\]
and the upper\,bound is attained exactly for those posets which are isomorphic to
\[A_{m-1} + (A_1 \oplus P).
\]
Hence, for disjoint finite sets $K$, $M$, $N$ with $\sharp M = m$ and $\sharp N = n$, and for any poset $P$ on $K$, 
the number of posets $Q\in \mathcal{Q}_0(K\cup M\cup N)$ with $Q|_K = P$, ${\rm Min}\, Q|_{K\cup M}\! = M$ and ${\rm Min}\, Q = N$ is
asymptotically equal to 
$$2^{\,(m-1)n}(d(P)+1)^n.$$
\ET

\newpage

\subsection*{Appendix: posets with at most five points}

\begin{footnotesize}

\hspace{.5ex}\begin{tabular}{|l|c|c|c|cc|ccccc|}
  \hline
  & $n$ &$\!$1$\!$&$\!$2$\!$& 3 & \hspace{-2ex}4 & 5 & \hspace{-1ex}6 & \hspace{-1ex}7 & \hspace{-1ex}8 & \hspace{-1ex}9\\
  \hline  
poset 
  & \begin{picture}(10,8)(-2,1)
    \put(-3,-2){$P_n$}
    \end{picture}
  & \hspace{-1ex}
    \begin{picture}(10,8)(-2.5,1) 
    \put(0,-2){$\emptyset$}
    \end{picture}
  & \hspace{-1.8ex}
    \begin{picture}(10,8)(-4,3) 
    \put(3,3){\circle*{3}}
    \end{picture}
  & \hspace{-1ex}
    \begin{picture}(10,8)(-2,3) 
    \put(3,0){\circle*{3}}
    \put(3,6){\circle*{3}}
    \put(3,0){\line(0,1){6}}
    \end{picture}
  & \hspace{-2.5ex}
    \begin{picture}(10,8)(-2,0) 
    \put(0,0){\circle*{3}}
    \put(6,0){\circle*{3}}
    \end{picture}
  & \hspace{-1.5ex}
    \begin{picture}(10,8)(-3,3) 
    \put(3,-2.8){\circle*{3}}
    \put(3,2.2){\circle*{3}}
    \put(3,7.2){\circle*{3}}
    \put(3,-2.8){\line(0,1){10}} 
    \end{picture} 
  & \hspace{-2ex}
    \begin{picture}(10,8)(-2,3) 
    \put(.5,6.7){\circle*{3}}
    \put(6.5,6.7){\circle*{3}}
    \put(3.5,-1.7){\circle*{3}}
    \put(.5,7.7){\line(1,-3){3.5}}
    \put(6.5,7.7){\line(-1,-3){3.5}}
    \end{picture}
  & \hspace{-2ex}
    \begin{picture}(10,8)(-2,3) 
    \put(0.5,-1.7){\circle*{3}}
    \put(6.5,-1.7){\circle*{3}}
    \put(3.5,6.7){\circle*{3}}
    \put(.5,-2.7){\line(1,3){3.5}}
    \put(6.5,-2.7){\line(-1,3){3.5}}
    \end{picture}
  & \hspace{-2ex}
    \begin{picture}(10,8)(-2,3) 
    \put(0,-1.7){\circle*{3}}
    \put(6,7){\circle*{3}}
    \put(6,-1.7){\circle*{3}}
    \put(6,-1.7){\line(0,1){10}}
        \end{picture}
  & \hspace{-2ex}
    \begin{picture}(10,8)(-2,0) 
    \put(-2,-1){\circle*{3}}
    \put(3,-1){\circle*{3}}
    \put(8,-1){\circle*{3}}
    \end{picture}\\[2ex]
  \hline  
number of points & $k_n$ &$\!$0$\!$&$\!$1$\!$&$\!$2$\!$&$\!$&$\!$3$\!$& & & &\\ 
  \hline 
number of minimal points & $m_n$ &$\!$0$\!$&$\!$1$\!$&$\!$1$\!$& \hspace{-2ex}2 & 1 &\hspace{-2ex} 1 & \hspace{-2ex}2 & \hspace{-2ex}2 & \hspace{-2ex}3\\
  \hline 
height\,(maximal\,size\,of\,chains)\hspace{-.7ex} &  $h_n$ &$\!$0$\!$&$\!$1$\!$&$\!$2$\!$& \hspace{-2ex}1 & 3 & \hspace{-1.5ex}2 & \hspace{-2ex}2 & \hspace{-2ex}2 & \hspace{-2ex}1\\
  \hline
number of automorphisms & $a_n$ &$\!$1$\!$&$\!$1$\!$&$\!$1$\!$& \hspace{-2ex}2 & 1 &\hspace{-2ex} 2 & \hspace{-2ex}2 & \hspace{-2ex}1 & \hspace{-2ex}6\\
  \hline
number of isomorphic copies &  $i_n$ &$\!$1$\!$&$\!$1$\!$&$\!$2$\!$& \hspace{-2ex}1 & 6 & \hspace{-2ex}3 & \hspace{-2ex}3 & \hspace{-2ex}6 & \hspace{-2ex}1\\
  \hline
number of downsets & $d_n$ &$\!$1$\!$&$\!$2$\!$&$\!$3$\!$& \hspace{-2ex}4 & 4 & \hspace{-2ex}5 & \hspace{-2ex}5 & \hspace{-2ex}6 & \hspace{-2ex}8\\
  \hline
\end{tabular}


\hspace{-4.7ex}
\begin{tabular}{|c|cccccccccccccccc|}
  \hline
  $n$ &$\!$10$\!$&$\!$11$\!$&$\!$12$\!$&$\!$13$\!$&$\!$14$\!$&$\!$15$\!$&$\!$16$\!$&$\!$17$\!$&$\!$18$\!$&$\!$19$\!$&$\!$20$\!$&$\!$21$\!$&$\!$22$\!$&$\!$23$\!$&$\!$24$\!$&$\!$25$\!$\\
  \hline
  \begin{picture}(10,8)(-2,1)
  \put(-3,-4){$P_n$}
  \end{picture}
  & \hspace{-2ex}
    \begin{picture}(10,8)(-4,3.5) 
    \put(3,-7){\circle*{3}}
    \put(3,-2){\circle*{3}}
    \put(3,3){\circle*{3}}
    \put(3,8){\circle*{3}}
    \put(3,-6){\line(0,1){14}}
    \end{picture} 
  & \hspace{-2ex}
    \begin{picture}(10,8)(-4,4) 
    \put(0,8.5){\circle*{3}}
    \put(6,8.5){\circle*{3}}
    \put(3,1){\circle*{3}}
    \put(3,-6.5){\circle*{3}}
    \put(-.5,9.8){\line(2,-5){4}}
    \put(6.5,9.8){\line(-2,-5){4}}
    \put(3,-7.5){\line(0,1){10}}
    \end{picture}
  & \hspace{-2ex}
    \begin{picture}(10,8)(-5,6) 
    \put(-.2,3.1){\circle*{3}}
    \put(6.2,3.1){\circle*{3}}
    \put(3,10.8){\circle*{3}}
    \put(3,-4.8){\circle*{3}}
    \put(-.5,4.4){\line(2,-5){4}}
    \put(6.5,4.4){\line(-2,-5){4}}
    \put(-.5,2){\line(2,5){4}}
    \put(6.5,2){\line(-2,5){4}}
    \end{picture}
  & \hspace{-2ex}
    \begin{picture}(10,8)(-5,6) 
    \put(-.2,3.4){\circle*{3}}
    \put(6.2,10.8){\circle*{3}}
    \put(-.2,10.8){\circle*{3}}
    \put(3,-4.8){\circle*{3}}
    \put(-.2,4.6){\line(2,-5){4}}
    \put(6.6,10.8){\line(-1,-5){3}}
    \put(-.2,2){\line(0,1){10}}
    \end{picture}
& \hspace{-2ex}
    \begin{picture}(10,8)(-6,4) 
    \put(-2.2,6.3){\circle*{3}}
    \put(8.2,6.3){\circle*{3}}
    \put(3,6.3){\circle*{3}}
    \put(3,-4.5){\circle*{3}}
    \put(-2.5,7.5){\line(2,-5){5}}
    \put(8.5,7.5){\line(-2,-5){5}}
    \put(3,6){\line(0,-1){10}}
    \end{picture}        
  & \hspace{-2ex}
    \begin{picture}(10,8)(-5,4) 
    \put(0,-6.5){\circle*{3}}
    \put(6,-6.5){\circle*{3}}
    \put(3,2){\circle*{3}}
    \put(3,9){\circle*{3}}
    \put(-.5,-7){\line(2,5){4}}
    \put(6.5,-7){\line(-2,5){4}}
    \put(3,.5){\line(0,1){10}}
    \end{picture}
  & \hspace{-2ex}
    \begin{picture}(10,8)(-6,6) 
    \put(6.2,2.6){\circle*{3}}
    \put(6.2,-4.8){\circle*{3}}
    \put(-.2,-4.8){\circle*{3}}
    \put(3.2,10.8){\circle*{3}}
    \put(7,1.4){\line(-2,5){4}}
    \put(-.2,-4.8){\line(1,5){3}}
    \put(6.2,4){\line(0,-1){10}}
    \end{picture}
  & \hspace{-2ex}
    \begin{picture}(10,8)(-4,3.5) 
    \put(6,-7){\circle*{3}}
    \put(0,-7){\circle*{3}}
    \put(6,.5){\circle*{3}}
    \put(6,8){\circle*{3}}
    \put(6,-6){\line(0,1){14}}
    \end{picture} 
  & \hspace{-2ex}
    \begin{picture}(10,8)(-5,4) 
    \put(-.2,6.3){\circle*{3}}
    \put(6.2,6.3){\circle*{3}}
    \put(-.2,-4.5){\circle*{3}}
    \put(6.2,-4.5){\circle*{3}}
    \put(0,-5){\line(1,2){6}}
    \put(6,-5){\line(-1,2){6}}
    \put(-.3,-5){\line(0,1){10}}
    \put(6.2,-5){\line(0,1){10}}
    \end{picture}
  & \hspace{-2ex}
    \begin{picture}(10,8)(-5,4) 
    \put(-.2,6.3){\circle*{3}}
    \put(6.2,6.3){\circle*{3}}
    \put(-.2,-4.5){\circle*{3}}
    \put(6.2,-4.5){\circle*{3}}
    \put(6,-5){\line(-1,2){6}}
    \put(-.3,-5){\line(0,1){10}}
    \put(6.2,-5){\line(0,1){10}}
    \end{picture}
  & \hspace{-2ex}
    \begin{picture}(10,8)(-5,4) 
    \put(-.2,6.3){\circle*{3}}
    \put(6.2,6.3){\circle*{3}}
    \put(-.2,-4.5){\circle*{3}}
    \put(6.2,-4.5){\circle*{3}}
    \put(-.3,-5){\line(0,1){10}}
    \put(6.2,-5){\line(0,1){10}}
    \end{picture}
  & \hspace{-2ex}
    \begin{picture}(10,8)(-5,4) 
    \put(.8,6.3){\circle*{3}}
    \put(8.2,6.3){\circle*{3}}
    \put(-3,-4.5){\circle*{3}}
    \put(4.5,-4.5){\circle*{3}}
    \put(.5,7.5){\line(1,-3){4}}
    \put(8.5,7.5){\line(-1,-3){4}}
    \end{picture}
  & \hspace{-2ex}
    \begin{picture}(10,8)(-4,4) 
    \put(-2.2,-4.5){\circle*{3}}
    \put(8.2,-4.5){\circle*{3}}
    \put(3,-4.5){\circle*{3}}
    \put(3,6.3){\circle*{3}}
    \put(-2.5,-5.5){\line(2,5){5}}
    \put(8.5,-5.5){\line(-2,5){5}}
    \put(3,6){\line(0,-1){10}}
    \end{picture}
  & \hspace{-2ex}
    \begin{picture}(10,8)(-4,4) 
    \put(-2.2,-4.5){\circle*{3}}
    \put(8.2,-4.5){\circle*{3}}
    \put(3,-4.5){\circle*{3}}
    \put(3,6.3){\circle*{3}}
    \put(-2.5,-5.5){\line(2,5){5}}
    \put(8.5,-5.5){\line(-2,5){5}}
    \end{picture}
  & \hspace{-2ex}
    \begin{picture}(10,8)(-5,4) 
    \put(-2.2,-4.5){\circle*{3}}
    \put(8.2,-4.5){\circle*{3}}
    \put(3,-4.5){\circle*{3}}
    \put(3,6.3){\circle*{3}}
    \put(3,6){\line(0,-1){10}}
    \end{picture}
  & \hspace{-1.5ex}
    \begin{picture}(10,8)(-4,-1) 
    \put(-3,-4.5){\circle*{3}}
    \put(1,-4.5){\circle*{3}}
    \put(5,-4.5){\circle*{3}}
    \put(9,-4.5){\circle*{3}}
    \end{picture}
  \\[3.5ex]

  \hline  
  $k_n$ & 4 & & & & & & & & & & & & & & &\\ 
  \hline 
  $m_n$ & 1 & 1 & 1 & 1 & 1 & 2 & 2 & 2 & 2 & 2 & 2 & 2 & 3 & 3 & 3 & 4\\
  \hline 
  $h_n$ & 4 & 3 & 3 & 3 & 2 & 3 & 3 & 3 & 2 & 2 & 2 & 2 & 2 & 2 & 2 & 1\\
  \hline
  $a_n$ & 1 & 2 & 2 & 1 & 6 & 2 & 1 & 1 & 4 & 1 & 2 & 2 & 6 & 2 & 2 &$\!$24$\!$\\
  \hline
  $i_n$ &$\!$24$\!$&$\!$12$\!$&$\!$12$\!$&$\!$24$\!$& 4 &$\!$12$\!$&$\!$24$\!$&$\!$24$\!$& 6 &$\!$24$\!$&$\!$12$\!$&$\!$12$\!$& 4 &$\!$12$\!$&$\!$12$\!$& 1\\
  \hline
  $d_n$ & 5 & 6 & 6 & 7 & 9 & 6 & 7 & 8 & 7 & 8 & 9 &$\!$10$\!$& 9 &$\!$10$\!$&$\!$12$\!$&$\!$16$\!$\\
  \hline
\end{tabular}


\hspace{-11.5ex}
\begin{tabular}{|c|ccccccccccccccccccccc|}
  \hline
$n$ &$\!$26$\!$&$\!$27$\!$&$\!$28$\!$&$\!$29$\!$&$\!$30$\!$&$\!$31$\!$&$\!$32$\!$&$\!$33$\!$&$\!$34$\!$&$\!$35$\!$&$\!$36$\!$&
$\!$37$\!$&$\!$38$\!$&$\!$39$\!$&$\!$40$\!$&$\!$41$\!$&$\!$42$\!$&$\!$43$\!$&$\!$44$\!$&$\!$45$\!$&$\!$46$\!$\\
  \hline
  \begin{picture}(10,8)(-2,1)
  \put(-3,-7){$P_n$}
  \end{picture}
  & \hspace{-2ex}
    \begin{picture}(10,8)(-4,5) 
    \put(3,-12){\circle*{3}}
    \put(3,-6.5){\circle*{3}}
    \put(3,-1){\circle*{3}}
    \put(3,4.5){\circle*{3}}
    \put(3,10){\circle*{3}}
    \put(3,-11){\line(0,1){20}}
    \end{picture}
  & \hspace{-2ex}
    \begin{picture}(10,8)(-4,4)
    \put(0,8.5){\circle*{3}}
    \put(6,8.5){\circle*{3}}
    \put(3,.5){\circle*{3}}
    \put(3,-13.3){\circle*{3}}
    \put(3,-6.3){\circle*{3}}
    \put(-.5,9.3){\line(2,-5){4}}
    \put(6.5,9.3){\line(-2,-5){4}}
    \put(3,-13){\line(0,1){14}}
    \end{picture}
  & \hspace{-2ex}
    \begin{picture}(10,8)(-5,7) 
    \put(-.2,3.9){\circle*{3}}
    \put(6.2,3.9){\circle*{3}}
    \put(3,11.7){\circle*{3}}
    \put(3,-3.7){\circle*{3}}
    \put(3,-10.8){\circle*{3}}
    \put(-.5,5.2){\line(2,-5){4}}
    \put(6.5,5.2){\line(-2,-5){4}}
    \put(-.5,2.8){\line(2,5){4}}
    \put(6.5,2.8){\line(-2,5){4}}
    \put(3,-12){\line(0,1){10}}
    \end{picture}
  & \hspace{-2ex}
    \begin{picture}(10,8)(-5,7) 
    \put(.2,3.9){\circle*{3}}
    \put(6.2,3.9){\circle*{3}}
    \put(.2,11.7){\circle*{3}}
    \put(3,-3.5){\circle*{3}}
    \put(3,-10.8){\circle*{3}}
    \put(-.5,5.2){\line(2,-5){4}}
    \put(6.5,5.2){\line(-2,-5){4}}
    \put(.2,2.8){\line(0,1){10}}
    \put(3,-12){\line(0,1){10}}
    \end{picture}
  & \hspace{-2ex}
    \begin{picture}(10,8)(-6,2) 
    \put(-2.2,6){\circle*{3}}
    \put(8.2,6){\circle*{3}}
    \put(3,6){\circle*{3}}
    \put(3,-4){\circle*{3}}
    \put(3,-14){\circle*{3}}
    \put(-2.5,7.5){\line(2,-5){5}}
    \put(8.5,7.5){\line(-2,-5){5}}
    \put(3,6){\line(0,-1){19}}
    \end{picture}
  & \hspace{-2ex}
    \begin{picture}(10,8)(-5,4) 
    \put(0,-6){\circle*{3}}
    \put(6,-6){\circle*{3}}
    \put(3,2){\circle*{3}}
    \put(3,9){\circle*{3}}
    \put(-.5,-6.5){\line(2,5){4}}
    \put(6.5,-6.5){\line(-2,5){4}}
    \put(3,.5){\line(0,1){10}}
    \put(-.5,-5){\line(2,-5){4}}
    \put(6.5,-5){\line(-2,-5){4}}
    \put(3,-14){\circle*{3}}
    \end{picture}
  & \hspace{-2ex}
    \begin{picture}(10,8)(-6,5) 
    \put(-.2,-4.8){\circle*{3}}
    \put(6.2,2.6){\circle*{3}}
    \put(6.2,-4.8){\circle*{3}}
    \put(3,10){\circle*{3}}
    \put(6.8,1.2){\line(-2,5){4}}
    \put(3,10.8){\line(-1,-5){3}}
    \put(-.5,-4){\line(2,-5){4}}
    \put(2.7,-13.5){\line(2,5){4}}
    \put(3,-13){\circle*{3}}
    \put(6.2,1.4){\line(0,-1){6}}
    \end{picture}
  & \hspace{-2ex}
    \begin{picture}(10,8)(-4,3.5) 
    \put(6,-6){\circle*{3}}
    \put(0,-6){\circle*{3}}
    \put(6,1.5){\circle*{3}}
    \put(6,9){\circle*{3}}
    \put(6,-5){\line(0,1){14}}
    \put(-.5,-5.5){\line(2,-5){4}}
    \put(6.5,-5.5){\line(-2,-5){4}}
    \put(3,-14.5){\circle*{3}}
    \end{picture}
  & \hspace{-2ex}
    \begin{picture}(10,8)(-5,2)
    \put(-.2,6){\circle*{3}}
    \put(6.2,6){\circle*{3}}
    \put(-.2,-4){\circle*{3}}
    \put(6.2,-4){\circle*{3}}
    \put(0,-5){\line(1,2){6}}
    \put(6,-5){\line(-1,2){6}}
    \put(-.3,-5){\line(0,1){10}}
    \put(6.2,-5){\line(0,1){10}}
    \put(-.5,-5){\line(2,-5){4}}
    \put(6.5,-5){\line(-2,-5){4}}
    \put(3,-14){\circle*{3}}
    \end{picture}
  & \hspace{-2ex}
    \begin{picture}(10,8)(-5,2) 
    \put(-.2,6){\circle*{3}}
    \put(6.2,6){\circle*{3}}
    \put(-.2,-4){\circle*{3}}
    \put(6.2,-4){\circle*{3}}
    \put(6,-5){\line(-1,2){6}}
    \put(-.3,-5){\line(0,1){10}}
    \put(6.2,-5){\line(0,1){10}}
    \put(-.5,-5){\line(2,-5){4}}
    \put(6.5,-5){\line(-2,-5){4}}
    \put(3,-14){\circle*{3}}
    \end{picture}
  & \hspace{-2ex}
    \begin{picture}(10,8)(-5,2) 
    \put(-.2,6){\circle*{3}}
    \put(6.2,6){\circle*{3}}
    \put(-.2,-4){\circle*{3}}
    \put(6.2,-4){\circle*{3}}
    \put(-.3,-5){\line(0,1){10}}
    \put(6.2,-5){\line(0,1){10}}
    \put(-.5,-5){\line(2,-5){4}}
    \put(6.5,-5){\line(-2,-5){4}}
    \put(3,-14){\circle*{3}}
    \end{picture}
  & \hspace{-2ex}
    \begin{picture}(10,8)(-6,2) 
    \put(-2.5,6){\circle*{3}}
    \put(4.5,6){\circle*{3}}
    \put(-2.5,-4){\circle*{3}}
    \put(4.5,-4){\circle*{3}}
    \put(-2,6){\line(2,-3){7}}
    \put(4.5,6.5){\line(0,-1){10}}
    \put(-2.5,-5){\line(2,-5){4}}
    \put(4.5,-5){\line(-2,-5){4}}
    \put(1,-14){\circle*{3}}
    \end{picture}
  & \hspace{-2ex}
    \begin{picture}(10,8)(-4,2) 
    \put(-2.2,-4){\circle*{3}}
    \put(8.2,-4){\circle*{3}}
    \put(3,-4){\circle*{3}}
    \put(3,6){\circle*{3}}
    \put(-2.5,-5.5){\line(2,5){5}}
    \put(8.5,-5.5){\line(-2,5){5}}
    \put(3,6){\line(0,-1){20}}
    \put(-1.5,-5){\line(2,-5){4}}
    \put(7.5,-5){\line(-2,-5){4}}
    \put(3,-14){\circle*{3}}
    \end{picture}
  & \hspace{-2ex}
    \begin{picture}(10,8)(-5,2) 
    \put(-2.2,-4){\circle*{3}}
    \put(8.2,-4){\circle*{3}}
    \put(3,-4){\circle*{3}}
    \put(3,6){\circle*{3}}
    \put(3,-14){\line(0,1){10}}
    \put(-1.5,-3){\line(2,5){4}}
    \put(-1.5,-5){\line(2,-5){4}}
    \put(7.5,-5){\line(-2,-5){4}}
    \put(3,-14){\circle*{3}}
    \put(4.2,5){\line(2,-5){4}}
    \end{picture}
  & \hspace{-2ex}
    \begin{picture}(10,8)(-5,2) 
    \put(-2.2,-4){\circle*{3}}
    \put(8.2,-4){\circle*{3}}
    \put(3,-4){\circle*{3}}
    \put(3,6){\circle*{3}}
    \put(3,6){\line(-0,-1){20}}
    \put(-1.5,-5){\line(2,-5){4}}
    \put(7.5,-5){\line(-2,-5){4}}
    \put(3,-14){\circle*{3}}
    \end{picture}
  & \hspace{-1.5ex}
    \begin{picture}(10,8)(-4,-2) 
    \put(-3,-4){\circle*{3}}
    \put(1,-4){\circle*{3}}
    \put(5,-4){\circle*{3}}
    \put(9,-4){\circle*{3}}
    \put(2,-13.5){\line(-3,5){6}}
    \put(2.5,-13.5){\line(-1,5){2}}
    \put(3.5,-13.5){\line(1,5){2}}
    \put(2.5,-14.5){\circle*{3}}
    \put(3.5,-13.5){\line(3,5){6}}
    \end{picture}
  & \hspace{-2ex}
    \begin{picture}(10,8)(-4,4)
    \put(0,-13){\circle*{3}}
    \put(6,-13){\circle*{3}}
    \put(3,8.5){\circle*{3}}
    \put(3,1.5){\circle*{3}}
    \put(3,-5.3){\circle*{3}}
    \put(-.5,-14){\line(2,5){4}}
    \put(6.5,-14){\line(-2,5){4}}
    \put(3,7){\line(0,-1){14}}
    \end{picture}
  & \hspace{-2ex}
    \begin{picture}(10,8)(-5,5) 
    \put(-.2,-4.9){\circle*{3}}
    \put(6.2,-12.3){\circle*{3}}
    \put(-.2,-12.3){\circle*{3}}
    \put(3,2.8){\circle*{3}}
    \put(-.8,-6.1){\line(2,5){4}}
    \put(6.8,-12.5){\line(-1,4){4}}
    \put(-.2,-3.5){\line(0,-1){10}}
    \put(3,11){\line(0,-1){10}}
    \put(3,10){\circle*{3}}
    \end{picture}
  & \hspace{-2ex}
    \begin{picture}(10,8)(-6,5) 
    \put(-.2,2.1){\circle*{3}}
    \put(6.2,-12.5){\circle*{3}}
    \put(-.2,-12.5){\circle*{3}}
    \put(-.2,-5){\circle*{3}}
    \put(-.8,0.9){\line(2,5){4}}
    \put(6.8,-12.5){\line(-1,6){4}}
    \put(-.2,3.5){\line(0,-1){15}}
    \put(3,10){\circle*{3}}
    \end{picture}
  & \hspace{-2ex}
    \begin{picture}(10,8)(-1,3.5) 
    \put(9,-13){\circle*{3}}
    \put(9,-6){\circle*{3}}
    \put(9,1){\circle*{3}}
    \put(9,8){\circle*{3}}
    \put(9,-13){\line(0,1){20}}
    \put(3,-13){\circle*{3}}
    \end{picture}
  & \hspace{-2ex}
    \begin{picture}(10,8)(-2,2) 
    \put(.8,6){\circle*{3}}
    \put(8.7,6){\circle*{3}}
    \put(4.7,-4){\circle*{3}}
    \put(.8,-14){\circle*{3}}
    \put(.5,6.5){\line(2,-5){4}}
    \put(9,6.5){\line(-2,-5){4}}
    \put(5,-5){\line(2,-5){4}}
    \put(4.5,-5){\line(-2,-5){4}}
    \put(8.7,-14){\circle*{3}}
    \end{picture}
  \\[5ex]

  \hline  
  $k_n$ & 5 & & & & & & & & & & & & & & & & & & & &\\ 
  \hline 
  $m_n$ & 1 & 1 & 1 & 1 & 1 & 1 & 1 & 1 & 1 & 1 & 1 & 1 & 1 & 1 & 1 & 1 & 2 & 2 & 2 & 2 & 2\\
  \hline 
  $h_n$ & 5 & 4 & 4 & 4 & 3 & 4 & 4 & 4 & 3 & 3 & 3 & 3 & 3 & 3 & 3 & 2 & 4 & 4 & 4 & 4 & 3\\
  \hline
  $a_n$ & 1 & 2 & 2 & 1 & 6 & 2 &  1 &  1 & 4 & 1 & 2 & 2 & 6 & 2 & 2 &$\!$24$\!$& 2 & 1 & 1 & 1 & 4\\
  \hline
  $i_n$ &$\!\!\!$120$\!\!\!$&$\!$60$\!$&$\!$60$\!$&$\!\!\!$120$\!\!\!$&$\!$20$\!$&$\!$60$\!$&$\!\!\!$120$\!\!\!$&$\!\!\!$120$\!\!\!$&$\!$30$\!$&$\!\!\!$120$\!\!\!$&$\!$60$\!$&
  $\!$60$\!$&$\!$20$\!$&$\!$60$\!$&$\!$60$\!$&$\!$5$\!$&$\!$60$\!$&$\!\!\!$120$\!\!\!$&$\!\!\!$120$\!\!\!$&$\!\!\!$120$\!\!\!$&$\!\!$30$\!\!$\\
  \hline
  $d_n$ & 6 & 7 & 7 & 8 &$\!$10$\!$& 7 & 8 & 9 & 8 & 9 &$\!$10$\!$&$\!$11$\!$&$\!$10$\!$&$\!$11$\!$&$\!$13$\!$&$\!$17$\!$& 7 & 8 & 9 &$\!$10$\!$ & 8\\
  \hline
\end{tabular}


\hspace{-11.5ex}
\begin{tabular}{|c|ccccccccccccccccccccc|}
  \hline
$n$ &$\!$47$\!$&$\!$48$\!$&$\!$49$\!$&$\!$50$\!$&$\!$51$\!$&$\!$52$\!$&$\!$53$\!$&$\!$54$\!$&$\!$55$\!$&$\!$56$\!$&$\!$57$\!$&
$\!$58$\!$&$\!$59$\!$&$\!$60$\!$&$\!$61$\!$&$\!$62$\!$&$\!$63$\!$&$\!$64$\!$&$\!$65$\!$&$\!$66$\!$&$\!$67$\!$\\
  \hline
  \begin{picture}(10,8)(-2,1)
  \put(-3,-7){$P_n$}
  \end{picture}
  & \hspace{-2ex}
    \begin{picture}(10,8)(-5,2)
    \put(-.2,-13.7){\circle*{3}}
    \put(6.2,-13.7){\circle*{3}}
    \put(-.2,-3.5){\circle*{3}}
    \put(6.2,-3.5){\circle*{3}}
    \put(0,-2.7){\line(1,-2){6}}
    \put(6,-2.7){\line(-1,-2){6}}
    \put(-.3,-3){\line(0,-1){10}}
    \put(6.2,-3){\line(0,-1){10}}
    \put(-.5,-2.7){\line(2,5){4}}
    \put(6.5,-2.7){\line(-2,5){4}}
    \put(3,6){\circle*{3}}
    \end{picture}
  & \hspace{-2ex}
    \begin{picture}(10,8)(-5,2)
    \put(-.2,-13.7){\circle*{3}}
    \put(6.2,-13.7){\circle*{3}}
    \put(-.2,-3.5){\circle*{3}}
    \put(6.2,-3.5){\circle*{3}}
    \put(0,-2.7){\line(1,-2){6}}
    \put(6,-2.7){\line(-1,-2){6}}
    \put(-.3,-3){\line(0,-1){10}}
    \put(6.2,-3){\line(0,-1){10}}
    \put(6.5,-2.7){\line(-2,5){4}}
    \put(3,6){\circle*{3}}
    \end{picture}
  & \hspace{-2ex}
    \begin{picture}(10,8)(-5,2)
    \put(-.2,-13.7){\circle*{3}}
    \put(6.2,-13.7){\circle*{3}}
    \put(-.2,-3.5){\circle*{3}}
    \put(6.2,-3.5){\circle*{3}}
    \put(0,-2.7){\line(1,-2){6}}
    \put(-.3,-3){\line(0,-1){10}}
    \put(6.2,-3){\line(0,-1){10}}
    \put(-.5,-2.7){\line(2,5){4}}
    \put(6.5,-2.7){\line(-2,5){4}}
    \put(3,6){\circle*{3}}
    \end{picture}
  & \hspace{-2ex}
    \begin{picture}(10,8)(-5,2)
    \put(-.2,-13.7){\circle*{3}}
    \put(6.2,-13.7){\circle*{3}}
    \put(-.2,-3.5){\circle*{3}}
    \put(6.2,-3.5){\circle*{3}}
    \put(-.3,-3){\line(0,-1){10}}
    \put(6.2,-3){\line(0,-1){10}}
    \put(-.5,-2.7){\line(2,5){4}}
    \put(6.5,-2.7){\line(-2,5){4}}
    \put(3,6){\circle*{3}}
    \end{picture}
  & \hspace{-2ex}
    \begin{picture}(10,8)(-5,2)
    \put(-.2,-13.7){\circle*{3}}
    \put(6.2,-13.7){\circle*{3}}
    \put(-.2,-3.5){\circle*{3}}
    \put(6.2,-3.5){\circle*{3}}
    \put(0,-2.7){\line(1,-2){6}}
    \put(0,-13.7){\line(1,6){3}}
    \put(6.2,-3){\line(0,-1){10}}
    \put(-.5,-2.7){\line(2,5){4}}
    \put(6.5,-2.7){\line(-2,5){4}}
    \put(3,6){\circle*{3}}
    \end{picture}
  & \hspace{-2ex}
    \begin{picture}(10,8)(-5,2)
    \put(-.2,6){\circle*{3}}
    \put(6.2,6){\circle*{3}}
    \put(-.2,-4){\circle*{3}}
    \put(6.2,-14){\circle*{3}}
    \put(0,-5){\line(1,2){6}}
    \put(5.7,-15){\line(-1,4){5}}
    \put(-.3,-5){\line(0,1){10}}
    \put(6.2,5){\line(0,-1){20}}
    \put(-.5,-5){\line(0,-1){10}}
    \put(-.2,-14){\circle*{3}}
    \end{picture}
  & \hspace{-2ex}
    \begin{picture}(10,8)(-5,2)
    \put(-.2,6){\circle*{3}}
    \put(6.2,6){\circle*{3}}
    \put(-.2,-4){\circle*{3}}
    \put(6.2,-14){\circle*{3}}
    \put(0,-15){\line(1,4){5}}
    \put(5.7,-15){\line(-1,4){5}}
    \put(-.3,-5){\line(0,1){10}}
    \put(6.2,5){\line(0,-1){20}}
    \put(-.5,-5){\line(0,-1){10}}
    \put(-.2,-14){\circle*{3}}
    \end{picture}
  & \hspace{-2ex}
    \begin{picture}(10,8)(-5,2)
    \put(-.2,6){\circle*{3}}
    \put(6.2,6){\circle*{3}}
    \put(6.2,-4){\circle*{3}}
    \put(6.2,-14){\circle*{3}}
    \put(5.7,-5){\line(-1,2){5}}
    \put(-.3,5){\line(0,-1){20}}
    \put(6.2,5){\line(0,-1){20}}
    \put(-.2,-14){\circle*{3}}
    \end{picture}
  & \hspace{-2ex}
    \begin{picture}(10,8)(-5,2)
    \put(-.2,-13.7){\circle*{3}}
    \put(6.2,-13.7){\circle*{3}}
    \put(-.2,-4){\circle*{3}}
    \put(6.2,6){\circle*{3}}
    \put(0,-2.7){\line(1,-2){6}}
    \put(-.3,-3){\line(0,-1){10}}
    \put(6.2,-13){\line(0,1){20}}
    \put(-.3,-2.7){\line(0,1){10}}
    \put(-.2,6){\circle*{3}}
    \end{picture}
  & \hspace{-2ex}
    \begin{picture}(10,8)(-5,2)
    \put(-.2,6){\circle*{3}}
    \put(6.2,6){\circle*{3}}
    \put(-.2,-4){\circle*{3}}
    \put(6.2,-14){\circle*{3}}
    \put(5.7,-15){\line(-1,4){5}}
    \put(-.3,5){\line(0,-1){20}}
    \put(6.2,5){\line(0,-1){20}}
    \put(-.2,-14){\circle*{3}}
    \end{picture}
  & \hspace{-2ex}
    \begin{picture}(10,8)(-5,2)
    \put(-.2,6){\circle*{3}}
    \put(6.2,6){\circle*{3}}
    \put(6.2,-4){\circle*{3}}
    \put(6.2,-14){\circle*{3}}
    \put(5.7,-15){\line(-1,4){5}}
    \put(-.3,5){\line(0,-1){20}}
    \put(6.2,5){\line(0,-1){20}}
    \put(-.2,-14){\circle*{3}}
    \end{picture}
  & \hspace{-2ex}
    \begin{picture}(10,8)(-5,2)
    \put(-.2,6){\circle*{3}}
    \put(6.2,6){\circle*{3}}
    \put(3,-4){\circle*{3}}
    \put(6.2,-14){\circle*{3}}
    \put(5.7,-15){\line(-1,4){5}}
    \put(-.3,5){\line(0,-1){20}}
    \put(6.2,5){\line(0,-1){20}}
    \put(-.2,-14){\circle*{3}}
    \end{picture}
  & \hspace{-2ex}
    \begin{picture}(10,8)(-5,2)
    \put(-.2,6){\circle*{3}}
    \put(6.2,6){\circle*{3}}
    \put(6.2,-4){\circle*{3}}
    \put(6.2,-14){\circle*{3}}
    \put(5.7,-5){\line(-1,2){5}}
    \put(6.2,5){\line(0,-1){20}}
    \put(-.2,-14){\circle*{3}}
    \end{picture}
  & \hspace{-2ex}
    \begin{picture}(10,8)(-5,2)
    \put(-.2,-14){\circle*{3}}
    \put(6.2,-14){\circle*{3}}
    \put(-.2,-4){\circle*{3}}
    \put(6.2,-4){\circle*{3}}
    \put(-.5,-2.7){\line(1,-2){6}}
    \put(6.2,-3){\line(0,-1){10}}
    \put(-.5,-2.7){\line(2,5){4}}
    \put(6.5,-2.7){\line(-2,5){4}}
    \put(3,6){\circle*{3}}
    \end{picture}
  & \hspace{-2ex}
    \begin{picture}(10,8)(-5,2)
    \put(-.2,6){\circle*{3}}
    \put(6.2,6){\circle*{3}}
    \put(6.2,-4){\circle*{3}}
    \put(6.2,-14){\circle*{3}}
    \put(-.3,5){\line(0,-1){20}}
    \put(6.2,5){\line(0,-1){20}}
    \put(-.2,-14){\circle*{3}}
    \end{picture}
  & \hspace{-2ex}
    \begin{picture}(10,8)(-5,2)
    \put(-.2,-13.7){\circle*{3}}
    \put(6.2,-13.7){\circle*{3}}
    \put(-.2,-4){\circle*{3}}
    \put(6.2,6){\circle*{3}}
    \put(0,-2.7){\line(1,-2){6}}
    \put(6.2,-13){\line(0,1){20}}
    \put(-.3,-2.7){\line(0,1){10}}
    \put(-.2,6){\circle*{3}}
    \end{picture}
  & \hspace{-2ex}
    \begin{picture}(10,8)(-3,2) 
    \put(-1,0){\circle*{3}}
    \put(4,0){\circle*{3}}
    \put(9,0){\circle*{3}}
    \put(-1,-10){\circle*{3}}
    \put(-1,0){\line(1,-1){10}}
    \put(9,0){\line(-1,-1){10}}
    \put(-1,0){\line(0,-1){10}}
    \put(9,0){\line(0,-1){10}}
    \put(4,0){\line(1,-2){5}}
    \put(4,0){\line(-1,-2){5}}
    \put(9,-10){\circle*{3}}
    \end{picture}
  & \hspace{-2ex}
    \begin{picture}(10,8)(-3,2) 
    \put(-1,0){\circle*{3}}
    \put(4,0){\circle*{3}}
    \put(9,0){\circle*{3}}
    \put(-1,-10){\circle*{3}}
    \put(-1,0){\line(1,-1){10}}
    \put(-1,0){\line(0,-1){10}}
    \put(9,0){\line(0,-1){10}}
    \put(4,0){\line(1,-2){5}}
    \put(4,0){\line(-1,-2){5}}
    \put(9,-10){\circle*{3}}
    \end{picture}
  & \hspace{-2ex}
    \begin{picture}(10,8)(-3,2) 
    \put(-1,0){\circle*{3}}
    \put(4,0){\circle*{3}}
    \put(9,0){\circle*{3}}
    \put(-1,-10){\circle*{3}}
    \put(-1,0){\line(0,-1){10}}
    \put(9,0){\line(0,-1){10}}
    \put(4,0){\line(1,-2){5}}
    \put(4,0){\line(-1,-2){5}}
    \put(9,-10){\circle*{3}}
    \end{picture}
  & \hspace{-2ex}
    \begin{picture}(10,8)(-3,2) 
    \put(-1,0){\circle*{3}}
    \put(4,0){\circle*{3}}
    \put(9,0){\circle*{3}}
    \put(-1,-10){\circle*{3}}
    \put(-1,0){\line(1,-1){10}}
    \put(-1,0){\line(0,-1){10}}
    \put(9,0){\line(0,-1){10}}
    \put(4,0){\line(1,-2){5}}
    \put(9,-10){\circle*{3}}
    \end{picture}
  & \hspace{-2ex}
    \begin{picture}(10,8)(-3,2) 
    \put(-1,0){\circle*{3}}
    \put(4,0){\circle*{3}}
    \put(9,0){\circle*{3}}
    \put(-1,-10){\circle*{3}}
    \put(-1,0){\line(0,-1){10}}
    \put(9,0){\line(0,-1){10}}
    \put(4,0){\line(1,-2){5}}
    \put(9,-10){\circle*{3}}
    \end{picture}
  \\[5ex]

  \hline  
  $k_n$ & 5 & & & & & & & & & & & & & & & & & & & &\\ 
  \hline 
  $m_n$ & 2 & 2 & 2 & 2 & 2 & 2 & 2 & 2 & 2 & 2 & 2 & 2 & 2 & 2 & 2 & 2 & 2 & 2 & 2 & 2 & 2\\
  \hline 
  $h_n$ & 3 & 3 & 3 & 3 & 3 & 3 & 3 & 3 & 3 & 3 & 3 & 3 & 3 & 3 & 3 & 3 & 2 & 2 & 2 & 2 & 2\\
  \hline
  $a_n$ & 4 & 2 & 1 & 2 & 2 & 2 & 1 & 1 & 1 & 1 & 1 & 1 & 2 & 2 & 1 & 1 &$\!$12$\!$& 2 & 2 & 2 & 2\\
  \hline
  $i_n$ &$\!$30$\!$&$\!$60$\!$&$\!\!\!$120$\!\!\!$&$\!$60$\!$&$\!$60$\!$&$\!\!\!$ 60$\!\!\!$&$\!\!\!$120$\!\!\!$&$\!\!\!$120$\!\!\!$&$\!\!\!$120$\!\!\!$&$\!\!\!$120$\!\!\!$&
  $\!\!\!$120$\!\!\!$&$\!\!\!$120$\!\!\!$&$\!$60$\!$&$\!$60$\!$&$\!\!\!$120$\!\!\!$&$\!\!\!$120$\!\!\!$&$\!$10$\!$&$\!$60$\!$&$\!$60$\!$&$\!$60$\!$&$\!$60$\hspace{.2ex}$\\
  \hline
  $d_n$ & 8 & 9 & 9 &$\!$10$\!$&$\!$11$\!$& 9 &$\!$10$\!$&$\!$10$\!$&$\!$10$\!$&$\!$11$\!$&$\!$11$\!$&$\!$12$\!$&$\!$12$\!$&$\!$12$\!$&$\!$12$\!$&$\!$14$\!$&$\!$11$\!$&$\!$12$\!$&$\!$13$\!$&$\!$14$\!$&$\!$15$\!$\\
  \hline
\end{tabular}


\hspace{-11.5ex}
\begin{tabular}{|c|ccccccccccccccccccccc|}
  \hline
$n$ &$\!$68$\!$&$\!$69$\!$&$\!$70$\!$&$\!$71$\!$&$\!$72$\!$&$\!$73$\!$&$\!$74$\!$&$\!$75$\!$&$\!$76$\!$&$\!$77$\!$&$\!$78$\!$&
$\!$79$\!$&$\!$80$\!$&$\!$81$\!$&$\!$82$\!$&$\!$83$\!$&$\!$84$\!$&$\!$85$\!$&$\!$86$\!$&$\!$87$\!$&$\!$88$\!$\\
  \hline
  \begin{picture}(10,8)(-2,1)
  \put(-3,-7){$P_n$}
  \end{picture}
  & \hspace{-2ex}
    \begin{picture}(10,8)(-3,2) 
    \put(-1,0){\circle*{3}}
    \put(4,0){\circle*{3}}
    \put(9,0){\circle*{3}}
    \put(-1,-10){\circle*{3}}
    \put(-1,0){\line(1,-1){10}}
    \put(9,0){\line(0,-1){10}}
    \put(4,0){\line(1,-2){5}}
    \put(9,-10){\circle*{3}}
    \end{picture}
  & \hspace{-2ex}
    \begin{picture}(10,8)(-3,2) 
    \put(-1,-14){\circle*{3}}
    \put(4,-14){\circle*{3}}
    \put(9,-14){\circle*{3}}
    \put(4,6){\circle*{3}}
    \put(-1,-14){\line(1,2){5}}
    \put(9,-14){\line(-1,2){5}}
    \put(4,-14){\line(0,1){20}}
    \put(4,-4){\circle*{3}}
    \end{picture}
  & \hspace{-2ex}
    \begin{picture}(10,8)(-3,2) 
    \put(-1,-14){\circle*{3}}
    \put(4,-14){\circle*{3}}
    \put(9,-14){\circle*{3}}
    \put(4,6){\circle*{3}}
    \put(-1,-14){\line(1,4){5}}
    \put(9,-14){\line(-1,4){5}}
    \put(4,-14){\line(1,4){2.5}}
    \put(6.5,-4){\circle*{3}}
    \end{picture}
  & \hspace{-2ex}
    \begin{picture}(10,8)(-3,2) 
    \put(-1,-14){\circle*{3}}
    \put(4,-14){\circle*{3}}
    \put(9,-14){\circle*{3}}
    \put(4,6){\circle*{3}}
    \put(9,-14){\line(-1,4){5}}
    \put(4,-14){\line(1,4){2.5}}
    \put(6.5,-4){\circle*{3}}
    \end{picture}
  & \hspace{-2ex}
    \begin{picture}(10,8)(-3,2) 
    \put(-1,-14){\circle*{3}}
    \put(4,-14){\circle*{3}}
    \put(9,-14){\circle*{3}}
    \put(4,6){\circle*{3}}
    \put(-1,-14){\line(1,4){5}}
    \put(9,-14){\line(-1,4){5}}
    \put(4,-14){\line(0,1){20}}
    \put(4,-4){\circle*{3}}
    \end{picture}
  & \hspace{-2ex}
    \begin{picture}(10,8)(-3,2) 
    \put(-1,-14){\circle*{3}}
    \put(4,-14){\circle*{3}}
    \put(9,-14){\circle*{3}}
    \put(4,6){\circle*{3}}
    \put(9,-14){\line(-1,4){5}}
    \put(4,-14){\line(0,1){20}}
    \put(4,-4){\circle*{3}}
    \end{picture}
  & \hspace{-2ex}
    \begin{picture}(10,8)(-3,2) 
    \put(-1,-14){\circle*{3}}
    \put(4,-14){\circle*{3}}
    \put(9,-14){\circle*{3}}
    \put(4,6){\circle*{3}}
    \put(4,-14){\line(0,1){20}}
    \put(4,-4){\circle*{3}}
    \end{picture}
  & \hspace{-2ex}
    \begin{picture}(10,8)(-3,2) 
    \put(-1,0){\circle*{3}}
    \put(4,-10){\circle*{3}}
    \put(9,0){\circle*{3}}
    \put(-1,-10){\circle*{3}}
    \put(-1,0){\line(1,-1){10}}
    \put(9,0){\line(-1,-1){10}}
    \put(-1,0){\line(0,-1){10}}
    \put(9,0){\line(0,-1){10}}
    \put(4,-10){\line(1,2){5}}
    \put(4,-10){\line(-1,2){5}}
    \put(9,-10){\circle*{3}}
    \end{picture}
  & \hspace{-2ex}
    \begin{picture}(10,8)(-3,2) 
    \put(-1,0){\circle*{3}}
    \put(4,-10){\circle*{3}}
    \put(9,0){\circle*{3}}
    \put(-1,-10){\circle*{3}}
    \put(-1,0){\line(1,-1){10}}
    \put(-1,0){\line(0,-1){10}}
    \put(9,0){\line(0,-1){10}}
    \put(4,-10){\line(1,2){5}}
    \put(4,-10){\line(-1,2){5}}
    \put(9,-10){\circle*{3}}
    \end{picture}
  & \hspace{-2ex}
    \begin{picture}(10,8)(-3,2) 
    \put(-1,0){\circle*{3}}
    \put(4,-10){\circle*{3}}
    \put(9,0){\circle*{3}}
    \put(-1,-10){\circle*{3}}
    \put(-1,0){\line(0,-1){10}}
    \put(9,0){\line(0,-1){10}}
    \put(4,-10){\line(1,2){5}}
    \put(4,-10){\line(-1,2){5}}
    \put(9,-10){\circle*{3}}
    \end{picture}
  & \hspace{-2ex}
    \begin{picture}(10,8)(-3,2) 
    \put(-1,0){\circle*{3}}
    \put(4,-10){\circle*{3}}
    \put(9,0){\circle*{3}}
    \put(-1,-10){\circle*{3}}
    \put(-1,0){\line(1,-1){10}}
    \put(9,0){\line(-1,-1){10}}
    \put(-1,0){\line(0,-1){10}}
    \put(9,0){\line(0,-1){10}}
    \put(9,-10){\circle*{3}}
    \end{picture}
  & \hspace{-2ex}
    \begin{picture}(10,8)(-3,2) 
    \put(-1,0){\circle*{3}}
    \put(4,-10){\circle*{3}}
    \put(9,0){\circle*{3}}
    \put(-1,-10){\circle*{3}}
    \put(-1,0){\line(1,-1){10}}
    \put(-1,0){\line(0,-1){10}}
    \put(9,0){\line(0,-1){10}}
    \put(4,-10){\line(-1,2){5}}
    \put(9,-10){\circle*{3}}
    \end{picture}
  & \hspace{-2ex}
    \begin{picture}(10,8)(-3,2) 
    \put(-1,0){\circle*{3}}
    \put(4,-10){\circle*{3}}
    \put(9,0){\circle*{3}}
    \put(-1,-10){\circle*{3}}
    \put(-1,0){\line(0,-1){10}}
    \put(9,0){\line(0,-1){10}}
    \put(4,-10){\line(-1,2){5}}
    \put(9,-10){\circle*{3}}
    \end{picture}
  & \hspace{-2ex}
    \begin{picture}(10,8)(-3,2) 
    \put(4,0){\circle*{3}}
    \put(4,-10){\circle*{3}}
    \put(9,0){\circle*{3}}
    \put(-1,-10){\circle*{3}}
    \put(4,0){\line(1,-2){5}}
    \put(4,0){\line(0,-1){10}}
    \put(9,0){\line(0,-1){10}}
    \put(9,-10){\circle*{3}}
    \end{picture}
  & \hspace{-2ex}
    \begin{picture}(10,8)(-3,2) 
    \put(-1,0){\circle*{3}}
    \put(4,-10){\circle*{3}}
    \put(9,0){\circle*{3}}
    \put(-1,-10){\circle*{3}}
    \put(-1,0){\line(0,-1){10}}
    \put(9,0){\line(0,-1){10}}
    \put(9,-10){\circle*{3}}
    \end{picture}
  & \hspace{-2ex}
    \begin{picture}(10,8)(-3,2) 
    \put(-1,0){\circle*{3}}
    \put(4,-10){\circle*{3}}
    \put(9,0){\circle*{3}}
    \put(-1,-10){\circle*{3}}
    \put(4,-10){\line(1,2){5}}
    \put(4,-10){\line(-1,2){5}}
    \put(9,-10){\circle*{3}}
    \end{picture}
  & \hspace{-2ex}
    \begin{picture}(10,8)(-3,1) 
    \put(-2,-10){\circle*{3}}
    \put(2,-10){\circle*{3}}
    \put(6,-10){\circle*{3}}
    \put(10,-10){\circle*{3}}
    \put(-2,-10){\line(1,2){5}}
    \put(2,-10){\line(1,6){1.7}}
    \put(6,-10){\line(-1,6){1.7}}
    \put(10,-10){\line(-1,2){5}}
     \put(4,0){\circle*{3}}
    \end{picture}
  & \hspace{-2ex}
    \begin{picture}(10,8)(-3,1) 
    \put(-2,-10){\circle*{3}}
    \put(2,-10){\circle*{3}}
    \put(6,-10){\circle*{3}}
    \put(10,-10){\circle*{3}}
    \put(2,-10){\line(1,6){1.7}}
    \put(6,-10){\line(-1,6){1.7}}
    \put(10,-10){\line(-1,2){5}}
     \put(4,0){\circle*{3}}
    \end{picture}
  & \hspace{-2ex}
    \begin{picture}(10,8)(-3,1) 
    \put(-2,-10){\circle*{3}}
    \put(2,-10){\circle*{3}}
    \put(6,-10){\circle*{3}}
    \put(10,-10){\circle*{3}}
    \put(2,-10){\line(1,6){1.7}}
    \put(6,-10){\line(-1,6){1.7}}
     \put(4,0){\circle*{3}}
    \end{picture}
  & \hspace{-2ex}
    \begin{picture}(10,8)(-3,1) 
    \put(-2,-10){\circle*{3}}
    \put(2,-10){\circle*{3}}
    \put(6,-10){\circle*{3}}
    \put(10,-10){\circle*{3}}
    \put(10,0){\line(0,-1){10}}
    \put(10,0){\circle*{3}}
    \end{picture}
  & \hspace{-2ex}
    \begin{picture}(10,8)(-3,-4) 
    \put(-4,-10){\circle*{3}}
    \put(0,-10){\circle*{3}}
    \put(4,-10){\circle*{3}}
    \put(8,-10){\circle*{3}}
    \put(12,-10){\circle*{3}}
    \end{picture}
  \\[5ex]
      
  \hline  
  $k_n$ & 5 & & & & & & & & & & & & & & & & & & & &\\ 
  \hline 
  $m_n$ & 2 & 3 & 3 & 3 & 3 & 3 & 3 & 3 & 3 & 3 & 3 & 3 & 3 & 3 & 3 & 3 & 4 & 4 & 4 & 4 & 5\\
  \hline 
  $h_n$ & 2 & 3 & 3 & 3 & 3 & 3 & 3 & 2 & 2 & 2 & 2 & 2 & 2 & 2 & 2 & 2 & 2 & 2 & 2 & 2 & 1\\
  
  \hline
  $a_n$ & 6 & 6 & 2 & 2 & 2 & 1 & 2 &$\!$12$\!$& 2 & 2 & 4 & 2 & 2 & 1 & 2 & 4 &$\!$24$\!$& 6 & 4 & 6 &$\!\!\!$120\\
  \hline
  $i_n$ &$\!$20$\!$&$\!$20$\!$&$\!$60$\!$&$\!$60$\!$&$\!$60$\!$&$\!\!\!$120$\!\!\!$&$\!$60$\!$&$\!$10$\!$&$\!$60$\!$&$\!$60$\!$&$\!$30$\!$&$\!$60$\!$&
  $\!$60$\!$&$\!\!\!$120$\!\!\!$&$\!$60$\!$&$\!$30$\!$&$\!$5$\!$&$\!$20$\!$&$\!$30$\!$&$\!$20$\!$&$\!$1$\!$\\
  \hline
  $d_n$ &$\!$18$\!$&$\!$10$\!$&$\!$11$\!$&$\!$12$\!$&$\!$13$\!$&$\!$14$\!$&$\!$16$\!$&$\!$11$\!$&$\!$12$\!$&$\!$13$\!$&$\!$14$\!$&$\!$14$\!$&$\!$15$\!$&$\!$16$\!$&$\!$18$\!$&$\!$20$\!$&$\!$17$\!$&$\!$18$\!$&$\!$20$\!$&$\!$24$\!$&$\!$32$\!$\\
  \hline
\end{tabular}

\end{footnotesize}
  
\newpage

\noindent {\bf \large The exponential functions}

\vspace{1ex}

\noindent $e_k(m)$ denotes the number of those labeled posets on $\underline{m+k}$ for which $\underline{m}$ is the set of minimal points. In particular, $e_k(0) = \delta_{k0}$.
By Theorem \ref{expo2},  

$\displaystyle{e_k(m) = \sum_{j=1}^{2^k} \sum_{i=0}^k (-1)^i b_{ijk}\, j^m \ \mbox{ with }}$

$\displaystyle{b_{ijk} = \sharp \{ (P,B)\mid P\in \mathcal{Q}_0(\underline{k}), \, B\subseteq {\rm Min}\, P, \, \sharp B = i, \, d(P \! - \! B) = j\}.}$

\vspace{1ex}

In particular, $b_{0jk} = p(j,k)$ and $b_{kjk} = \delta_{j1}$.

\vspace{1ex}

\noindent For $k\leq 5$, one obtains the following exponential sums (cf.\ Theorem \ref{asym}): 

\vspace{1ex}

$e_0(m) =1$

\vspace{1ex}

$e_1(m) = 2^m-1 = 2^m + O(1)$ 

\vspace{1ex}

$e_2(m) = 4^m+2\cdot 3^m-4\cdot 2^m+1 = 4^m + 2\,\binom{2}{2}\, 3^m + O(2^m)$ 

\vspace{1ex}

$e_3(m) = 8^m+6\cdot 6^m+6\cdot 5^m-6\cdot 4^m -18\cdot 3^m+12\cdot 2^m-1$

\hspace{6.7ex}$= 8^m + 2\,\binom{3}{2}\, 6^m + 6\,\binom{3}{3}\, 5^m + O(4^m)$ 

\vspace{1ex}

$e_4(m) = 16^m+12\cdot12^m+24\cdot 10^m+20\cdot 9^m+16\cdot 8^m+54\cdot 7^m-108\cdot 6^m$

$\hspace{9ex}-\,96\cdot 5^m+108\cdot 3^m-32\cdot 2^m+1$

\hspace{6.7ex}$= 16^m + 2\,\binom{4}{2}\, 12^m + 6\,\binom{4}{3}\, 10^m + 20\, \binom{4}{4}\, 9^m + O(8^m)$ 

\vspace{1ex}

$e_5(m) =  32^m+20\cdot 24^m+60\cdot 20^m+100\cdot 18^m+10\cdot 17^m+100\cdot16^m$

$\hspace{9ex}+120\cdot 15^m+390\cdot 14^m +240\cdot 13^m-180\cdot 12^m+500\cdot 11^m-540\cdot 10^m$

$\hspace{9ex}-300\cdot 9^m - 830\cdot 8^m -1650\cdot 7^m +1200\cdot 6^m +900  \cdot 5^m+320\cdot 4^m$

$\hspace{9ex}-540\cdot 3^m+80\cdot 2^m -1$

\hspace{6.7ex}$= 32^m + 2\,\binom{5}{2}\, 24^m + 6\,\binom{5}{3}\, 20^m + 20\, \binom{5}{4}\, 18^m + 10\, \binom{5}{5}\, 17^m  + O(16^m)$ 

\vspace{1ex}

\noindent These sums are obtained by summation of the exponential functions of all posets $P_n$ with $k_n=k$, taken with multiplicity $i_n = i(P_n)$, as listed on the following two pages.
The total number of posets with $k$ points is then

$\displaystyle{p(k) = e_k(1) = e_{k-1}(2) = \sum_{m=1}^k \binom{k}{m} e_{k-m}(m)}$

\vspace{-1ex}

$p(1) =1$

$p(2) = 3$

$p(3) = 19$

$p(4)= 219$

$p(5)= 4231$

$p(6) = 130023$

\newpage

\begin{small}

\begin{center}
\begin{tabular}{|r|r|l|}
  \hline
  $n$ & $i_n$ & $e(m,P_n)$\\  
  \hline
   1 &  1 & $\,\ 1^{m}$ \\
  \hline
   2 &  1 & $\,\ 2^{m}-1^{m}$ \\
  \hline
   3 &  2 & $\,\ 3^{m}-2^{m}$ \\
   4 &  1 & $\,\ 4^{m}-2\cdot 2^{m}+1^{m}$ \\
  \hline
   5 &  6 & $\,\ 4^{m}-3^{m}$ \\
   6 &  3 & $\,\ 5^{m}-4^{m}$ \\
   7 &  3 & $\,\ 5^{m}-2\cdot 3^{m}+2^{m}$ \\
   8 &  6 & $\,\ 6^{m}-4^{m}-3^{m}+2^{m}$ \\
   9 &  1 & $\,\ 8^{m}-3\cdot 4^{m}+3\cdot 2^{m}-1^{m}$ \\
  \hline
  10 & 24 & $\,\ 5^{m}-4^{m}$ \\
  11 & 12 & $\,\ 6^{m}-5^{m}$ \\
  12 & 12 & $\,\ 6^{m}-5^{m}$ \\
  13 & 24 & $\,\ 7^{m}-6^{m}$ \\
  14 &  4 & $\,\ 9^{m}-8^{m}$ \\
  15 & 12 & $\,\ 6^{m}-2\cdot 4^{m}+3^{m}$ \\
  16 & 24 & $\,\ 7^{m}-5^{m}-4^{m}+3^{m}$ \\
  17 & 24 & $\,\ 8^{m}-6^{m}-4^{m}+3^{m}$ \\
  18 &  6 & $\,\ 7^{m}-2\cdot 5^{m}+4^{m}$ \\
  19 & 24 & $\,\ 8^{m}-6^{m}-5^{m}+4^{m}$ \\
  20 & 12 & $\,\ 9^{m}-2 \cdot 6^{m}+4^{m}$ \\
  21 & 12 & $10^{m}-8^{m}-5^{m}+4^{m}$ \\
  22 &  4 & $\,\ 9^{m}-3 \cdot 5^{m}+3\cdot 3^{m}-2^{m}$ \\
  23 & 12 & $10^{m}-2 \cdot 6^{m}-5^{m}+4^{m}+2\cdot 3^{m}-2^{m}$ \\
  24 & 12 & $12^{m}-8^{m}-2\cdot 6^{m}+2\cdot 4^{m}+3^{m}-2^{m}$ \\
  25 &  1 & $16^{m}-4 \cdot 8^{m}+6\cdot 4^{m} -4\cdot 2^{m} +1^{m}$ \\
  \hline
    26 & 120 & $\,\ 6^{m}-5^{m}$ \\
  27 &  60 & $\,\ 7^{m}-6^{m}$ \\
  28 &  60 & $\,\ 7^{m}-6^{m}$ \\
  29 & 120 & $\,\ 8^{m}-7^{m}$ \\
  30 &  20 & $10^{m}-9^{m}$ \\
  31 &  60 & $\,\ 7^{m}-6^{m}$ \\
  32 & 120 & $\,\ 8^{m}-7^{m}$ \\
  33 & 120 & $\,\ 9^{m}-8^{m}$ \\
  34 &  30 & $\,\ 8^{m}-7^{m}$ \\
  35 & 120 & $\,\ 9^{m}-8^{m}$ \\
  36 &  60 & $10^{m}-9^{m}$ \\
  37 &  60 & $11^{m}-10^{m}$ \\
  38 &  20 & $10^{m}-9^{m}$ \\
  39 &  60 & $11^{m}-10^{m}$ \\
  40 &  60 & $13^{m}-12^{m}$ \\
  41 &   5 & $17^{m}-16^{m}$ \\
  42 &  60 & $\,\ 7^{m}-2 \cdot 5^{m}+4^{m}$ \\
  43 & 120 & $\,\ 8^{m}-6^{m}-5^{m}+4^{m}$ \\
  44 & 120 & $\,\ 9^{m}-7^{m}-5^{m}+4^{m} \ \hspace{47ex} $ \\
  \hline
\end{tabular} 
\end{center}

\begin{center}
\begin{tabular}{|r|r|l|}
  \hline
  $n$ & $i_n$ & $e(m,P_n)$\\  
  \hline

  45 & 120 & $10^{m}-8^{m}-5^{m}+4^{m}$ \\
  46 &  30 & $\,\ 8^{m}-2 \cdot 6^{m}+5^{m}$ \\
  47 &  30 & $\,\ 8^{m}-2 \cdot 6^{m}+5^{m}$ \\
  48 &  60 & $\,\ 9^{m}-2 \cdot 7^{m}+6^{m}$ \\
  49 & 120 & $\,\ 9^{m}-7^{m}-6^{m}+5^{m}$ \\
  50 &  60 & $10^{m}-2 \cdot 7^{m}+5^{m}$ \\
  51 &  60 & $11^{m}-9^{m}-6^{m}+5^{m}$ \\
  52 &  60 & $\,\ 9^{m}-7^{m}-6^{m}+5^{m}$ \\
  53 & 120 & $10^{m}-8^{m}-7^{m}+6^{m}$ \\
  54 & 120 & $10^{m}-8^{m}-6^{m}+5^{m}$ \\
  55 & 120 & $10^{m}-8^{m}-7^{m}+6^{m}$ \\
  56 & 120 & $11^{m}-2 \cdot 8^{m}+6^{m}$ \\
  57 & 120 & $11^{m}-9^{m}-7^{m}+6^{m}$ \\
  58 & 120 & $12^{m}-10^{m}-7^{m}+6^{m}$ \\
  59 &  60 & $12^{m}-10^{m}-6^{m}+5^{m}$ \\
  60 &  60 & $12^{m}-10^{m}-6^{m}+5^{m}$ \\
  61 & 120 & $12^{m}-9^{m}-8^{m}+6^{m}$ \\
  62 & 120 & $14^{m}-12^{m}-7^{m}+6^{m}$ \\
  63 &  10 & $11^{m}-2 \cdot 9^{m}+8^{m}$ \\
  64 &  60 & $12^{m}-10^{m}-9^{m}+8^{m}$ \\
  65 &  60 & $13^{m}-2 \cdot 10^{m}+8^{m}$ \\
  66 &  60 & $14^{m}-12^{m}-9^{m}+8^{m}$ \\
  67 &  60 & $15^{m}-12^{m}-10^{m}+8^{m}$ \\
  68 &  20 & $18^{m}-16^{m}-9^{m}+8^{m}$ \\
  69 &  20 & $10^{m}-3 \cdot 6^{m}+3 \cdot 4^{m}-3^{m}$ \\
  70 &  60 & $11^{m}-2 \cdot 7^{m}-6^{m}+5^{m}+2 \cdot 4^{m}-3^{m}$ \\
  71 &  60 & $12^{m}-2 \cdot 8^{m}+2 \cdot 4^{m}-3^{m}$ \\
  72 &  60 & $13^{m}-9^{m}-2 \cdot 7^{m}+2 \cdot 5^{m}+4^{m}-3^{m}$ \\
  73 & 120 & $14^{m}-10^{m}-8^{m}-7^{m}+6^{m}+5^{m}+4^{m}-3^{m}$ \\
  74 &  60 & $16^{m}-12^{m}-2 \cdot 8^{m}+2 \cdot 6^{m}+4^{m}-3^{m}$ \\
  75 &  10 & $11^{m}-3 \cdot 7^{m}+3 \cdot 5^{m}-4^{m}$ \\
  76 &  60 & $12^{m}-2 \cdot 8^{m}-7^{m}+6^{m}+2 \cdot 5^{m}-4^{m}$ \\
  77 &  60 & $13^{m}-9^{m}-2 \cdot 8^{m}+2 \cdot 6^{m}+5^{m}-4^{m}$ \\
  78 &  30 & $14^{m}-2 \cdot 10^{m}+8^{m}-7^{m}+2 \cdot 5^{m}-4^{m}$ \\
  79 &  60 & $14^{m}-10^{m}-2 \cdot 8^{m}+2 \cdot 6^{m}+5^{m}-4^{m}$ \\
  80 &  60 & $15^{m}-10^{m}-2 \cdot 9^{m}+3 \cdot 6^{m}-4^{m}$ \\
  81 & 120 & $16^{m}-12^{m}-10^{m}+6^{m}+5^{m}-4^{m}$ \\
  82 &  60 & $18^{m}-2 \cdot 12^{m}-9^{m}+8^{m}+2 \cdot 6^{m}-4^{m}$ \\
  83 &  30 & $20^{m}-16^{m}-2 \cdot 10^{m}+2 \cdot 8^{m}+5^{m}-4^{m}$ \\
  84 &   5 & $17^{m}-4 \cdot 9^{m}+6 \cdot 5^{m}-4 \cdot 3^{m}+2^{m}$ \\
  85 &  20 & $18^{m}-3 \cdot 10^{m}-9^{m}+3 \cdot 6^{m}+3 \cdot 5^{m}-4^{m}-3 \cdot 3^{m}+2^{m}$ \\
  86 &  30 & $20^{m}-2 \cdot 12^{m}-2 \cdot 10^{m}+8^{m}+4 \cdot 6^{m}+5^{m}-2 \cdot 4^{m}-2 \cdot 3^{m}+2^{m}$ \\
  87 &  20 & $24^{m} -16^{m} -3 \cdot 12^{m} +3 \cdot 8^{m} +3 \cdot 6^{m} -3 \cdot 4^{m} - 3^{m} + 2^{m}$ \\
  88 &   1 & $32^{m} -5 \cdot 16^{m} +10 \cdot 8^{m} -10 \cdot 4^{m} +5 \cdot 2^{m} -1$ \\  
  \hline
\end{tabular} 
\end{center}

\end{small}

\newpage

\begin{footnotesize}

\noindent In the same way one calculates the following two types of exponential sums. 
Let $e_{kn}(m)$ denote the exponential sum counting all posets $Q$ on $M\cup K$ such that \mbox{$\sharp\,{\rm Min}\,(Q|_K  ) = n$}  and ${\rm Min}\,Q = M$. 
Then, by Corollary \ref{kn} (with $n$ for $m$) and Theorem \ref{expo}, the leading term of $e_{kn}(m)$ is
$
\textstyle{\binom{k}{n}\cdot n \cdot (2^{\,k-1}\! + 2^{\,n-1})^m \ \mbox{ if } n < k, \mbox{ and } 2^{\,km} \mbox{ if } n = k.}
$
\begin{eqnarray*}
e_{11}(m) & \!\! = \!\! & 2^m\! - 1\\
e_{21}(m) & \!\! = \!\! & 2\cdot 3^m\! - 2\cdot 2^m\\
e_{22}(m) & \!\! = \!\! & 4^m\! - 2\cdot 2^m\! +1\\
e_{31}(m) & \!\! = \!\! & 3\cdot 5^m\! + 3\cdot 4^m\! - 6\cdot 3^m\\
e_{32}(m) & \!\! = \!\! & 6\cdot 6^m\! + 3\cdot 5^m\! - 6\cdot 4^m\! - 12\cdot 3^m\! + 9\cdot 2^m\\
e_{33}(m) & \!\! = \!\! & 8^m\! - 3\cdot 4^m\! + 3\cdot 2^m\! - 1\\
e_{41}(m) & \!\! = \!\! & 4\cdot 9^m\! - 4\cdot 8^m\! + 24\cdot 7^m\! - 24\cdot 4^m\\
e_{42}(m) & \!\! = \!\! & 12\cdot 10^m\! + 12\cdot 9^m\! + 36\cdot 8^m\! + 30\cdot 7^m\! - 60\cdot 6^m\! - 72\cdot 5^m\! - 18\cdot 4^m\! + 60\cdot 3^m\\
e_{43}(m) & \!\! = \!\! & 12\cdot 12^m\! + 12\cdot 10^m\! + 4\cdot 9^m\! - 12\cdot 8^m\! - 48\cdot 6^m\! - 24\cdot 5^m\! + 36\cdot 4^m\! + 48\cdot 3^m\! - 28\cdot 2^m\\
e_{44}(m) & \!\! = \!\! & 16^m\! - 4\cdot 8^m\! + 6\cdot 4^m\! - 4\cdot 2^m\! +1\\
e_{51}(m) & \!\! = \!\! & 5\cdot 17^m\! - 5\cdot 16^m\! + 60\cdot 13^m\! - 60\cdot 12^m\! + 120\cdot 11^m\! - 20\cdot 10^m\! + 140\cdot 9^m\! + 30\cdot 8^m \\
          &             & -90\cdot 7^m\! -60\cdot 6^m\! - 120\cdot 5^m\\
e_{52}(m) & \!\! = \!\! & 20\cdot 18 ^m\! - 20\cdot 16^m\! + 60\cdot 15^m\! + 180\cdot 14^m\! + 60\cdot 13^m\! + 180\cdot 12^m\! + 310\cdot 11^m \\
          &             & + 60\cdot 10^m\! - 100\cdot 9^m\! - 390\cdot 8^m\! - 1080\cdot 7^m\! + 180\cdot 6^m\! + 120\cdot 5^m\! + 420\cdot 4^m\\    
e_{53}(m) & \!\! = \!\! & 30\cdot 20^m\!\! + 60\cdot 18^m\!\! + 150\cdot 16^m\!\! + 60\cdot 15^m\!\! + 210\cdot 14^m\!\! + 120\cdot 13^m\!\! - 180\cdot 12^m\!\! + 70\cdot 11^m \\
          &             & - 460\cdot 10^m\! - 300\cdot 9^m\! - 570\cdot 8^m\! - 480\cdot 7^m\! + 840\cdot 6^m\! + 780\cdot 5^m\! + 50\cdot 4^m\! - 380\cdot 3^m\\
e_{54}(m) & \!\! = \!\! & 20\cdot 24^m\! + 30\cdot 20^m\! + 20\cdot 18^m\! + 5\cdot 17^m\! - 20\cdot 16^m\! - 120\cdot 12^m\! -120\cdot 10^m\! - 40\cdot 9^m \\
          &             & + 90\cdot 8^m\! + 240\cdot 6^m\! + 120\cdot 5^m\! - 140\cdot 4^m\! - 160\cdot 3^m\! + 75\cdot 2^m\\  
e_{55}(m) & \!\! = \!\! & 32^m - 5\cdot 16^m\! +10\cdot 8^m\! -10\cdot 4^m\! + 5\cdot 2^m\! - 1       
\end{eqnarray*}

\vspace{1ex}

\noindent Let $e_k^h(m)$ denote the exponential sum counting all posets $Q$ on $M \cup K$ of height $h\!+\!1$ with ${\rm Min}\,Q \!=\! M$. 
By Corollary \ref{kh} (with $n$ for $m$) and Corollary \ref{exposum2}, we have

\vspace{1ex}

\noindent \mbox{$e_k^1(m) = (2^m-1)^k,\ e_k^h(m) = (k)_h\,(2^{k-h}(h+1))^m + ... \ (1 \!<\! h \!<\! k),\ e_k^k(m) = k!((k+1)^m - k^m).$}

\vspace{-2ex}

\begin{eqnarray*}
e_3^1(m) & \!\! = \!\! & 8^m\! - 3\cdot 4^m\! + 3\cdot 2^m\! -1\\
e_3^2(m) & \!\! = \!\! & 6\cdot 6^m + 6\cdot 5^m\! - 9\cdot 4^m\! -12\cdot 3^m\! +9\cdot 2^m\\
e_3^3(m) & \!\! = \!\! & 6\cdot 4^m - 6\cdot 3^m \\
e_4^1(m) & \!\! = \!\! & 16^m\! - 4\cdot 8^m\! + 6\cdot 4^m\! - 4\cdot 2^m\! + 1\\
e_4^2(m) & \!\! = \!\! & 12\cdot 12^m\!\! + 24\cdot 10^m\!\! + 20\cdot 9^m\!\! - 4\cdot 8^m\!\! + 6\cdot 7^m\!\! -96\cdot 6^m\!\! - 72\cdot 5^m\!\! + 90\cdot 4^m\!\! +48\cdot 3^m\!\! - 28\cdot 2^m\\
e_4^3(m) & \!\! = \!\! & 24\cdot 8^m\! + 48\cdot 7^m\! - 12\cdot 6^m\! - 48\cdot 5^m\! - 72\cdot 4^m\! + 60\cdot 3^m\\
e_4^4(m) & \!\! = \!\! & 24\cdot 5^m - 24\cdot 4^m \\
e_5^1(m) & \!\! = \!\! & 32^m\! - 5\cdot 16^m\! +10\cdot 8^m\! -10\cdot 4^m\! + 5\cdot 2^m\! - 1\\
e_5^2(m) & \!\! = \!\! & 20\cdot 24^m\! + 60\cdot 20^m + 100\cdot 18^m\! + 10\cdot 17^m\! + 45\cdot 16^m + 120\cdot 15^m\! + 150\cdot 14^m\! + 120\cdot 13^m\!\\
         &             &  - 360\cdot 12^m\! + 20\cdot 11^m\! - 720\cdot 10^m\!\! - 440\cdot 9^m\!\! + 150\cdot 8^m\!\! - 120\cdot 7^m\!\! + 960\cdot 6^m\!\! + 600\cdot 5^m\\
         &             &  - 630\cdot 4^m\!\! - 160\cdot 3^m\!\! + 75\cdot 2^m\\
e_5^3(m) & \!\! = \!\! & 60\cdot 16^m\! + 240\cdot 14^m\! + 120\cdot 13^m\! + 180\cdot 12^m\!  + 480\cdot 11^m\! + 60 \cdot 10^m\! - 100 \cdot 9^m\! - 1110\cdot 8^m \\
         &             &  - 1410\cdot 7^m\! + 420\cdot 6^m\! + 900\cdot 5^m\! + 540\cdot 4^m\! - 380\cdot 3^m\\
e_5^4(m) & \!\! = \!\! & 120\cdot 10^m\! + 240 \cdot 9^m\! + 120\cdot 8^m\! - 120\cdot 7^m\! - 300\cdot 6^m\! - 480\cdot 5^m\! + 420\cdot 4^m\\
e_5^5(m) & \!\! = \!\! & 120\cdot 6^m - 120\cdot 5^m 
\end{eqnarray*}

\newpage
\noindent {\bf \large The representing matrices}

\hspace{-4ex}\begin{tabular}{|c|c|cc|ccccc|cccccccccccccccc|}
   \hline  
  \hspace{-2ex}
    \begin{picture}(3,8)(-1,-5) 
    \put(0,-2){$\emptyset$}
    \end{picture}
&  \hspace{-2ex}
    \begin{picture}(5,8)(-1,-4) 
    \put(3,3){\circle*{3}}
    \end{picture}
&  \hspace{-2ex}
    \begin{picture}(5,8)(1,-4) 
    \put(3,0){\circle*{3}}
    \put(3,6){\circle*{3}}
    \put(3,0){\line(0,1){6}}
    \end{picture}
&  \hspace{-2ex}
    \begin{picture}(4,8)(1,-7) 
    \put(0,0){\circle*{3}}
    \put(5,0){\circle*{3}}
    \end{picture}
&  \hspace{-2ex}
    \begin{picture}(5,8)(0,-5) 
    \put(3,-2.8){\circle*{3}}
    \put(3,2.2){\circle*{3}}
    \put(3,7.2){\circle*{3}}
    \put(3,-2.8){\line(0,1){10}} 
    \end{picture} 
&  \hspace{-2ex}
    \begin{picture}(5,8)(0,-5) 
    \put(.5,6.7){\circle*{3}}
    \put(6.5,6.7){\circle*{3}}
    \put(3.5,-1.7){\circle*{3}}
    \put(.5,7.7){\line(1,-3){3.5}}
    \put(6.5,7.7){\line(-1,-3){3.5}}
    \end{picture}
&  \hspace{-2ex}
    \begin{picture}(5,8)(0,-5) 
    \put(0,-1.7){\circle*{3}}
    \put(6,-1.7){\circle*{3}}
    \put(3,6.7){\circle*{3}}
    \put(0,-2.7){\line(1,3){3.5}}
    \put(6,-2.7){\line(-1,3){3.5}}
    \end{picture}
 &  \hspace{-2ex}
    \begin{picture}(5,8)(0,-5) 
    \put(.5,-1.7){\circle*{3}}
    \put(5.5,7){\circle*{3}}
    \put(5.5,-1.7){\circle*{3}}
    \put(5.5,-1.7){\line(0,1){10}}
    \end{picture}
 &  \hspace{-2ex}
    \begin{picture}(5,9)(-1,-8) 
    \put(-2,-1){\circle*{3}}
    \put(2,-1){\circle*{3}}
    \put(6,-1){\circle*{3}}
    \end{picture}   

  & \hspace{-2ex}
    \begin{picture}(7,18)(-2,-7) 
    \put(3,-7){\circle*{3}}
    \put(3,-2){\circle*{3}}
    \put(3,3){\circle*{3}}
    \put(3,8){\circle*{3}}
    \put(3,-6){\line(0,1){14}}
    \end{picture} 
  & \hspace{-3ex}
    \begin{picture}(10,18)(-3,-7) 
    \put(1,8.5){\circle*{3}}
    \put(7,8.5){\circle*{3}}
    \put(4,1){\circle*{3}}
    \put(4,-6.5){\circle*{3}}
    \put(.5,9.8){\line(2,-5){4}}
    \put(7.5,9.8){\line(-2,-5){4}}
    \put(4,-7.5){\line(0,1){10}}
    \end{picture}
  & \hspace{-4ex}
    \begin{picture}(10,18)(-5,-5) 
    \put(-.2,3.1){\circle*{3}}
    \put(6.2,3.1){\circle*{3}}
    \put(3,10.8){\circle*{3}}
    \put(3,-4.8){\circle*{3}}
    \put(-.5,4.4){\line(2,-5){4}}
    \put(6.5,4.4){\line(-2,-5){4}}
    \put(-.5,2){\line(2,5){4}}
    \put(6.5,2){\line(-2,5){4}}
    \end{picture}
  & \hspace{-4ex}
    \begin{picture}(10,18)(-5,-5) 
    \put(-.2,3.4){\circle*{3}}
    \put(6.2,10.8){\circle*{3}}
    \put(-.2,10.8){\circle*{3}}
    \put(3,-4.8){\circle*{3}}
    \put(-.2,4.6){\line(2,-5){4}}
    \put(6.6,10.8){\line(-1,-5){3}}
    \put(-.2,2){\line(0,1){10}}
    \end{picture}
  & \hspace{-4ex}
    \begin{picture}(10,18)(-6,-7) 
    \put(-3.2,6.3){\circle*{3}}
    \put(7.2,6.3){\circle*{3}}
    \put(2,6.3){\circle*{3}}
    \put(2,-4.5){\circle*{3}}
    \put(-3.5,7.5){\line(2,-5){5}}
    \put(7.5,7.5){\line(-2,-5){5}}
    \put(2,6){\line(0,-1){10}}
    \end{picture}        
  & \hspace{-4ex}
    \begin{picture}(10,18)(-5,-7) 
    \put(0,-6.5){\circle*{3}}
    \put(6,-6.5){\circle*{3}}
    \put(3,2){\circle*{3}}
    \put(3,9){\circle*{3}}
    \put(-.5,-7){\line(2,5){4}}
    \put(6.5,-7){\line(-2,5){4}}
    \put(3,.5){\line(0,1){10}}
    \end{picture}
  & \hspace{-4ex}
    \begin{picture}(10,18)(-6,-5) 
    \put(5.2,2.6){\circle*{3}}
    \put(5.2,-4.8){\circle*{3}}
    \put(-1.2,-4.8){\circle*{3}}
    \put(2.2,10.8){\circle*{3}}
    \put(6,1.4){\line(-2,5){4}}
    \put(-.2,-4.8){\line(1,5){3}}
    \put(5.2,4){\line(0,-1){10}}
    \end{picture}
  & \hspace{-4ex}
    \begin{picture}(10,18)(-4,-7) 
    \put(6,-7){\circle*{3}}
    \put(0,-7){\circle*{3}}
    \put(6,.5){\circle*{3}}
    \put(6,8){\circle*{3}}
    \put(6,-6){\line(0,1){14}}
    \end{picture} 
  & \hspace{-4ex}
    \begin{picture}(10,18)(-5,-7) 
    \put(-.2,6.3){\circle*{3}}
    \put(6.2,6.3){\circle*{3}}
    \put(-.2,-4.5){\circle*{3}}
    \put(6.2,-4.5){\circle*{3}}
    \put(0,-5){\line(1,2){6}}
    \put(6,-5){\line(-1,2){6}}
    \put(-.3,-5){\line(0,1){10}}
    \put(6.2,-5){\line(0,1){10}}
    \end{picture}
  & \hspace{-4ex}
    \begin{picture}(10,18)(-5,-7) 
    \put(-.2,6.3){\circle*{3}}
    \put(6.2,6.3){\circle*{3}}
    \put(-.2,-4.5){\circle*{3}}
    \put(6.2,-4.5){\circle*{3}}
    \put(6,-5){\line(-1,2){6}}
    \put(-.3,-5){\line(0,1){10}}
    \put(6.2,-5){\line(0,1){10}}
    \end{picture}
  & \hspace{-4ex}
    \begin{picture}(10,18)(-5,-7) 
    \put(-.2,6.3){\circle*{3}}
    \put(6.2,6.3){\circle*{3}}
    \put(-.2,-4.5){\circle*{3}}
    \put(6.2,-4.5){\circle*{3}}
    \put(-.3,-5){\line(0,1){10}}
    \put(6.2,-5){\line(0,1){10}}
    \end{picture}
  & \hspace{-4ex}
    \begin{picture}(10,18)(-5,-7) 
    \put(.8,6.3){\circle*{3}}
    \put(8.2,6.3){\circle*{3}}
    \put(-1,-4.5){\circle*{3}}
    \put(4.5,-4.5){\circle*{3}}
    \put(.5,7.5){\line(1,-3){4}}
    \put(8.5,7.5){\line(-1,-3){4}}
    \end{picture}
  & \hspace{-4ex}
    \begin{picture}(10,18)(-4,-7) 
    \put(-1.2,-4.5){\circle*{3}}
    \put(3,-4.5){\circle*{3}}
    \put(7.2,-4.5){\circle*{3}}
    \put(3,6.3){\circle*{3}}
    \put(-2,-5.5){\line(2,5){5}}
    \put(8,-5.5){\line(-2,5){5}}
    \put(3,6){\line(0,-1){10}}
    \end{picture}  
  & \hspace{-4ex}
    \begin{picture}(10,18)(-4,-7) 
    \put(-1.2,-4.5){\circle*{3}}
    \put(3,-4.5){\circle*{3}}
    \put(7.2,-4.5){\circle*{3}}
    \put(3,6.3){\circle*{3}}
    \put(-2,-5.5){\line(2,5){5}}
    \put(8,-5.5){\line(-2,5){5}}
    \end{picture}    
  & \hspace{-4ex}
    \begin{picture}(10,18)(-4,-7) 
    \put(-1.2,-4.5){\circle*{3}}
    \put(3,-4.5){\circle*{3}}
    \put(7.2,-4.5){\circle*{3}}
    \put(3,6.3){\circle*{3}}
    \put(3,6){\line(0,-1){10}}
    \end{picture}  
    
  & \hspace{-4ex}
    \begin{picture}(10,18)(-4,-12) 
    \put(-3,-4.5){\circle*{3}}
    \put(1,-4.5){\circle*{3}}
    \put(5,-4.5){\circle*{3}}
    \put(9,-4.5){\circle*{3}}
    \end{picture}
   \\
    
  \hline
    \hspace{-.5ex}1\hspace{-3ex}&\hspace{-.5ex}2\hspace{-2.5ex}&\hspace{-.5ex}3\hspace{-3ex}&\hspace{-2ex}4\hspace{-3ex}&\hspace{-1ex}5\hspace{-1ex}&\hspace{-.5ex}6\hspace{-1ex}&\hspace{-.5ex}7\hspace{-1ex}&\hspace{-.5ex}8\hspace{-1ex}&
    \hspace{-1ex}9\hspace{1ex}&\hspace{-1ex}10\hspace{-3ex}&\hspace{-2ex}11\hspace{-2ex}&\hspace{-2ex}12\hspace{-1ex}&\hspace{-1.5ex}13\hspace{-1ex}&\hspace{-2ex}14\hspace{-1ex}&\hspace{-2ex}15\hspace{-1ex}&
    \hspace{-2ex}16\hspace{-1ex}&\hspace{-2ex}17\hspace{-1ex}&\hspace{-2ex}18\hspace{-1ex}&\hspace{-2ex}19\hspace{-1ex}&\hspace{-2ex}20\hspace{-1ex}&\hspace{-2ex}21\hspace{-1ex}&\hspace{-2ex}22\hspace{-1ex}&
    \hspace{-2ex}23\hspace{-1ex}&\hspace{-1.5ex}24\hspace{-.5ex}&\hspace{-2ex}25\hspace{-1ex}\\
  \hline
\end{tabular}

\vspace{1ex}

\noindent {\normalsize $A: \ a_{mn} = \sum_{P_m \simeq U\in\, \mathcal{U}(P_n)} d(U^{\circ}) \ ( = \delta_{mn}\,d(P_n) \mbox{ for } m\geq 10 \mbox{ and } n\leq 25 )$} 


\hspace{-4ex}\begin{tabular}{|c|c|cc|ccccc|cccccccccccccccc|}
   \hline  
  1\hspace{-.7ex} & 1\hspace{-.3ex} & 1\hspace{-.7ex} & 1\hspace{-.7ex} & 1\hspace{-.7ex} & 1\hspace{-.7ex} & 1\hspace{-.7ex} & 1\hspace{-.7ex} & 1\hspace{-.7ex} & 1\hspace{-.7ex} & 
  1\hspace{-.7ex} & 1\hspace{-.7ex} & 1\hspace{-.7ex} & 1\hspace{-.7ex} & 1\hspace{-.7ex} & 1\hspace{-.7ex} & 1\hspace{-.7ex} & 1\hspace{-.7ex} & 1\hspace{-.7ex} & 1\hspace{-.7ex} & 
  1\hspace{-.7ex} & 1\hspace{-.7ex} & 1\hspace{-.7ex} & 1\hspace{-.7ex} & 1\hspace{.5ex}$ $\\
   \hline
 
   \hspace{-.7ex} & 2\hspace{-.7ex} & 1\hspace{-.7ex} & 4\hspace{-.7ex} & 1\hspace{-.7ex} & 2\hspace{-.7ex} & 1\hspace{-.7ex} & 3\hspace{-.7ex} & 6\hspace{-.7ex} & 1\hspace{-.7ex} & 
  2\hspace{-.7ex} & 1\hspace{-.7ex} & 2\hspace{-.7ex} & 3\hspace{-.7ex} & 1\hspace{-.7ex} & 1\hspace{-.7ex} & 3\hspace{-.7ex} & 2\hspace{-.7ex} & 2\hspace{-.7ex} & 2\hspace{-.7ex} & 
  4\hspace{-.7ex} & 1\hspace{-.7ex} & 3\hspace{-.7ex} & 5\hspace{-.7ex} & 8\hspace{.5ex}$ $\\
   \hline  
   
   \hspace{-.7ex} &  \hspace{-.7ex} & 3\hspace{-.7ex} & 0\hspace{-.7ex} & 1\hspace{-.7ex} & 0\hspace{-.7ex} & 4\hspace{-.7ex} & 3\hspace{-.7ex} & 0\hspace{-.7ex} & 1\hspace{-.7ex} & 
  0\hspace{-.7ex} & 2\hspace{-.7ex} & 1\hspace{-.7ex} & 0\hspace{-.7ex} & 1\hspace{-.7ex} & 3\hspace{-.7ex} & 1\hspace{-.7ex} & 0\hspace{-.7ex} & 2\hspace{-.7ex} & 6\hspace{-.7ex} & 
  0\hspace{-.7ex} & 6\hspace{-.7ex} & 4\hspace{-.7ex} & 3\hspace{-.7ex} & 0\hspace{.5ex}$ $\\
  
   \hspace{-.7ex} &  \hspace{-.7ex} &  \hspace{-.7ex} & 4\hspace{-.7ex} & 0\hspace{-.7ex} & 1\hspace{-.7ex} & 0\hspace{-.7ex} & 2\hspace{-.7ex} & \hspace{-1.5ex}12\hspace{-1ex} & 0\hspace{-.7ex} & 
  1\hspace{-.7ex} & 0\hspace{-.7ex} & 1\hspace{-.7ex} & 3\hspace{-.7ex} & 0\hspace{-.7ex} & 0\hspace{-.7ex} & 2\hspace{-.7ex} & 1\hspace{-.7ex} & 1\hspace{-.7ex} & 1\hspace{-.7ex} & 
  5\hspace{-.7ex} & 0\hspace{-.7ex} & 2\hspace{-.7ex} & 8 \hspace{-1.7ex} & \hspace{-1ex}24\\
   \hline 
   
   \hspace{-.7ex} &  \hspace{-.7ex} &  \hspace{-.7ex} &  \hspace{-.7ex} & 4\hspace{-.7ex} & 0\hspace{-.7ex} & 0\hspace{-.7ex} & 0\hspace{-.7ex} & 0\hspace{-.7ex} & 1\hspace{-.7ex} & 
  0\hspace{-.7ex} & 0\hspace{-.7ex} & 0\hspace{-.7ex} & 0\hspace{-.7ex} & 4\hspace{-.7ex} & 3\hspace{-.7ex} & 4\hspace{-.7ex} & 0\hspace{-.7ex} & 0\hspace{-.7ex} & 0\hspace{-.7ex} &   
  0\hspace{-.7ex} & 0\hspace{-.7ex} & 0\hspace{-.7ex} & 0\hspace{-.7ex} & 0\hspace{.5ex}$ $\\
  
   \hspace{-.7ex} &  \hspace{-.7ex} &  \hspace{-.7ex} &  \hspace{-.7ex} &  \hspace{-.7ex} & 5\hspace{-.7ex} & 0\hspace{-.7ex} & 0\hspace{-.7ex} & 0\hspace{-.7ex} & 0\hspace{-.7ex} & 
  1\hspace{-.7ex} & 0\hspace{-.7ex} & 0\hspace{-.7ex} & 0\hspace{-.7ex} & 0\hspace{-.7ex} & 0\hspace{-.7ex} & 0\hspace{-.7ex} & 4\hspace{-.7ex} & 3\hspace{-.7ex} & 0\hspace{-.7ex} & 
  5\hspace{-.7ex} & 0\hspace{-.7ex} & 0\hspace{-.7ex} & 0\hspace{-.7ex} & 0\hspace{.5ex}$ $\\
  
   \hspace{-.7ex} &  \hspace{-.7ex} &  \hspace{-.7ex} &  \hspace{-.7ex} &  \hspace{-.7ex} &  \hspace{-.7ex} & 5\hspace{-.7ex} & 0\hspace{-.7ex} & 0\hspace{-.7ex} & 0\hspace{-.7ex} & 
  0\hspace{-.7ex} & 1\hspace{-.7ex} & 0\hspace{-.7ex} & 0\hspace{-.7ex} & 0\hspace{-.7ex} & 2\hspace{-.7ex} & 0\hspace{-.7ex} & 0\hspace{-.7ex} & 0\hspace{-.7ex} & 0\hspace{-.7ex} & 
  0\hspace{-.7ex} & \hspace{-1ex}12\hspace{-1ex} & 5\hspace{-.7ex} & 0\hspace{-.7ex} & 0\hspace{.5ex}$ $\\
  
   \hspace{-.7ex} &  \hspace{-.7ex} &  \hspace{-.7ex} &  \hspace{-.7ex} &  \hspace{-.7ex} &  \hspace{-.7ex} &  \hspace{-.7ex} & 6\hspace{-.7ex} & 0\hspace{-.7ex} & 0\hspace{-.7ex} & 
  0\hspace{-.7ex} & 0\hspace{-.7ex} & 1\hspace{-.7ex} & 0\hspace{-.7ex} & 0\hspace{-.7ex} & 0\hspace{-.7ex} & 2\hspace{-.7ex} & 0\hspace{-.7ex} & 2\hspace{-.7ex} & 6\hspace{-.7ex} & 
  0\hspace{-.7ex} & 0\hspace{-.7ex} & 8\hspace{-.7ex} & \hspace{-1ex}12\hspace{-1ex} & 0\hspace{.5ex}$ $\\
  
   \hspace{-.7ex} &  \hspace{-.7ex} &  \hspace{-.7ex} &  \hspace{-.7ex} &  \hspace{-.7ex} &  \hspace{-.7ex} &  \hspace{-.7ex} &  \hspace{-.7ex} & 8\hspace{-.7ex} & 0\hspace{-.7ex} & 
  0\hspace{-.7ex} & 0\hspace{-.7ex} & 0\hspace{-.7ex} & 1\hspace{-.7ex} & 0\hspace{-.7ex} & 0\hspace{-.7ex} & 0\hspace{-.7ex} & 0\hspace{-.7ex} & 0\hspace{-.7ex} & 0\hspace{-.7ex} & 
  2\hspace{-.7ex} & 0\hspace{-.7ex} & 0\hspace{-.7ex} & 4\hspace{-.7ex} & \hspace{-1ex}32\\
  \hline   
       
\end{tabular}

\vspace{1ex}

\noindent {\normalsize B : \ $b_{mn} = \sharp \{U \in \mathcal{U}(P_n),\, U \simeq P_m\} \ ( = \delta_{mn} \mbox{ for } m\geq 10 \mbox{ and } n\leq 25 )$}


\hspace{-4ex}\begin{tabular}{|c|c|cc|ccccc|cccccccccccccccc|}
   \hline  
  1\hspace{-.7ex} & 1\hspace{-.3ex} & 1\hspace{-.7ex} & 1\hspace{-.7ex} & 1\hspace{-.7ex} & 1\hspace{-.7ex} & 1\hspace{-.7ex} & 1\hspace{-.7ex} & 1\hspace{-.7ex} & 1\hspace{-.7ex} & 
  1\hspace{-.7ex} & 1\hspace{-.7ex} & 1\hspace{-.7ex} & 1\hspace{-.7ex} & 1\hspace{-.7ex} & 1\hspace{-.7ex} & 1\hspace{-.7ex} & 1\hspace{-.7ex} & 1\hspace{-.7ex} & 1\hspace{-.7ex} & 
  1\hspace{-.7ex} & 1\hspace{-.7ex} & 1\hspace{-.7ex} & 1\hspace{-.7ex} & 1\hspace{.5ex}$ $\\
   \hline
 
   \hspace{-.7ex} & 1\hspace{-.7ex} & 1\hspace{-.7ex} & 2\hspace{-.7ex} & 1\hspace{-.7ex} & 2\hspace{-.7ex} & 1\hspace{-.7ex} & 2\hspace{-.7ex} & 3\hspace{-.7ex} & 1\hspace{-.7ex} & 
  2\hspace{-.7ex} & 1\hspace{-.7ex} & 2\hspace{-.7ex} & 3\hspace{-.7ex} & 1\hspace{-.7ex} & 1\hspace{-.7ex} & 2\hspace{-.7ex} & 2\hspace{-.7ex} & 2\hspace{-.7ex} & 2\hspace{-.7ex} & 
  3\hspace{-.7ex} & 1\hspace{-.7ex} & 2\hspace{-.7ex} & 3\hspace{-.7ex} & 4\hspace{.5ex}$ $\\
   \hline  
   
   \hspace{-.7ex} &  \hspace{-.7ex} & 1\hspace{-.7ex} & 0\hspace{-.7ex} & 1\hspace{-.7ex} & 0\hspace{-.7ex} & 2\hspace{-.7ex} & 1\hspace{-.7ex} & 0\hspace{-.7ex} & 1\hspace{-.7ex} & 
  0\hspace{-.7ex} & 2\hspace{-.7ex} & 1\hspace{-.7ex} & 0\hspace{-.7ex} & 1\hspace{-.7ex} & 2\hspace{-.7ex} & 1\hspace{-.7ex} & 0\hspace{-.7ex} & 1\hspace{-.7ex} & 2\hspace{-.7ex} & 
  0\hspace{-.7ex} & 3\hspace{-.7ex} & 2\hspace{-.7ex} & 1\hspace{-.7ex} & 0\hspace{.5ex}$ $\\
  
   \hspace{-.7ex} &  \hspace{-.7ex} &  \hspace{-.7ex} & 1\hspace{-.7ex} & 0\hspace{-.7ex} & 1\hspace{-.7ex} & 0\hspace{-.7ex} & 1\hspace{-.7ex} & 3\hspace{-.7ex} & 0\hspace{-.7ex} & 
  1\hspace{-.7ex} & 0\hspace{-.7ex} & 1\hspace{-.7ex} & 3\hspace{-.7ex} & 0\hspace{-.7ex} & 0\hspace{-.7ex} & 1\hspace{-.7ex} & 1\hspace{-.7ex} & 1\hspace{-.7ex} & 1\hspace{-.7ex} & 
  3\hspace{-.7ex} & 0\hspace{-.7ex} & 1\hspace{-.7ex} & 3\hspace{-.7ex} & 6\hspace{.5ex}$ $\\
   \hline 
   
   \hspace{-.7ex} &  \hspace{-.7ex} &  \hspace{-.7ex} &  \hspace{-.7ex} & 1\hspace{-.7ex} & 0\hspace{-.7ex} & 0\hspace{-.7ex} & 0\hspace{-.7ex} & 0\hspace{-.7ex} & 1\hspace{-.7ex} & 
  0\hspace{-.7ex} & 0\hspace{-.7ex} & 0\hspace{-.7ex} & 0\hspace{-.7ex} & 2\hspace{-.7ex} & 1\hspace{-.7ex} & 1\hspace{-.7ex} & 0\hspace{-.7ex} & 0\hspace{-.7ex} & 0\hspace{-.7ex} &   
  0\hspace{-.7ex} & 0\hspace{-.7ex} & 0\hspace{-.7ex} & 0\hspace{-.7ex} & 0\hspace{.5ex}$ $\\
  
   \hspace{-.7ex} &  \hspace{-.7ex} &  \hspace{-.7ex} &  \hspace{-.7ex} &  \hspace{-.7ex} & 1\hspace{-.7ex} & 0\hspace{-.7ex} & 0\hspace{-.7ex} & 0\hspace{-.7ex} & 0\hspace{-.7ex} & 
  1\hspace{-.7ex} & 0\hspace{-.7ex} & 0\hspace{-.7ex} & 0\hspace{-.7ex} & 0\hspace{-.7ex} & 0\hspace{-.7ex} & 0\hspace{-.7ex} & 2\hspace{-.7ex} & 1\hspace{-.7ex} & 0\hspace{-.7ex} & 
  1\hspace{-.7ex} & 0\hspace{-.7ex} & 0\hspace{-.7ex} & 0\hspace{-.7ex} & 0\hspace{.5ex}$ $\\
  
   \hspace{-.7ex} &  \hspace{-.7ex} &  \hspace{-.7ex} &  \hspace{-.7ex} &  \hspace{-.7ex} &  \hspace{-.7ex} & 1\hspace{-.7ex} & 0\hspace{-.7ex} & 0\hspace{-.7ex} & 0\hspace{-.7ex} & 
  0\hspace{-.7ex} & 1\hspace{-.7ex} & 0\hspace{-.7ex} & 0\hspace{-.7ex} & 0\hspace{-.7ex} & 1\hspace{-.7ex} & 0\hspace{-.7ex} & 0\hspace{-.7ex} & 0\hspace{-.7ex} & 0\hspace{-.7ex} & 
  0\hspace{-.7ex} & 3\hspace{-.7ex} & 1\hspace{-.7ex} & 0\hspace{-.7ex} & 0\hspace{.5ex}$ $\\
  
   \hspace{-.7ex} &  \hspace{-.7ex} &  \hspace{-.7ex} &  \hspace{-.7ex} &  \hspace{-.7ex} &  \hspace{-.7ex} &  \hspace{-.7ex} & 1\hspace{-.7ex} & 0\hspace{-.7ex} & 0\hspace{-.7ex} & 
  0\hspace{-.7ex} & 0\hspace{-.7ex} & 1\hspace{-.7ex} & 0\hspace{-.7ex} & 0\hspace{-.7ex} & 0\hspace{-.7ex} & 1\hspace{-.7ex} & 0\hspace{-.7ex} & 1\hspace{-.7ex} & 2\hspace{-.7ex} & 
  0\hspace{-.7ex} & 0\hspace{-.7ex} & 2\hspace{-.7ex} & 2\hspace{-.7ex} & 0\hspace{.5ex}$ $\\
  
   \hspace{-.7ex} &  \hspace{-.7ex} &  \hspace{-.7ex} &  \hspace{-.7ex} &  \hspace{-.7ex} &  \hspace{-.7ex} &  \hspace{-.7ex} &  \hspace{-.7ex} & 1\hspace{-.7ex} & 0\hspace{-.7ex} & 
  0\hspace{-.7ex} & 0\hspace{-.7ex} & 0\hspace{-.7ex} & 1\hspace{-.7ex} & 0\hspace{-.7ex} & 0\hspace{-.7ex} & 0\hspace{-.7ex} & 0\hspace{-.7ex} & 0\hspace{-.7ex} & 0\hspace{-.7ex} & 
  1\hspace{-.7ex} & 0\hspace{-.7ex} & 0\hspace{-.7ex} & 1\hspace{-.7ex} & 4\hspace{.5ex}$ $\\
  \hline   
       
\end{tabular}

\vspace{1ex}

\noindent {\normalsize $C: \ c_{mn} = \delta_{mn} - \sum_{j<m} b_{mj} c_{jn} \ ( = \delta_{mn} \mbox{ for } m\geq 10 \mbox{ and } n\leq 25 )$} 


\hspace{-4ex}\begin{tabular}{|c|c|cc|ccccc|cccccccccccccccc|}
   \hline  
  1\hspace{-.7ex} & \hspace{-.5ex}-1\hspace{-.3ex} & 0\hspace{-.7ex} & 1\hspace{-.7ex} & 0\hspace{-.7ex} & 0\hspace{-.7ex} & 0\hspace{-.7ex} & 0\hspace{-.7ex} & \hspace{-1ex}-1\hspace{-.7ex} & \hspace{-.5ex}0\hspace{-1ex} & 
  0\hspace{-.7ex} & 0\hspace{-.7ex} & 0\hspace{-.7ex} & 0\hspace{-.7ex} & 0\hspace{-.7ex} & 0\hspace{-.7ex} & 0\hspace{-.7ex} & 0\hspace{-.7ex} & 0\hspace{-.7ex} & 0\hspace{-.7ex} & 
  0\hspace{-.7ex} & 0\hspace{-.7ex} & 0\hspace{-.7ex} & 0\hspace{-.7ex} & 1\hspace{.5ex}$ $\\
   \hline
 
   \hspace{-.7ex} & 1\hspace{-.7ex} & \hspace{-1ex}-1\hspace{-.7ex} & \hspace{-1ex}-2\hspace{-.7ex} & 0\hspace{-.7ex} & 0\hspace{-.7ex} & 1\hspace{-.7ex} & 1\hspace{-.7ex} & 3\hspace{-.7ex} & 0\hspace{-.7ex} & 
  0\hspace{-.7ex} & 0\hspace{-.7ex} & 0\hspace{-.7ex} & 0\hspace{-.7ex} & 0\hspace{-.7ex} & 0\hspace{-.7ex} & 0\hspace{-.7ex} & 0\hspace{-.7ex} & 0\hspace{-.7ex} & 0\hspace{-.7ex} & 
  0\hspace{-.7ex} & \hspace{-1ex}-1\hspace{-.7ex} & \hspace{-1ex}-1\hspace{-.7ex} & \hspace{-1ex}-1\hspace{-.7ex} & \hspace{-1ex}-4\hspace{.5ex}$ $\\
   \hline  
   
   \hspace{-.7ex} &  \hspace{-.7ex} & 1\hspace{-.7ex} & 0\hspace{-.7ex} & \hspace{-1ex}-1\hspace{-.7ex} & 0\hspace{-.7ex} & \hspace{-1ex}-2\hspace{-.7ex} & \hspace{-1ex}-1\hspace{-.7ex} & 0\hspace{-.7ex} & 0\hspace{-.7ex} & 
  0\hspace{-.7ex} & 0\hspace{-.7ex} & 0\hspace{-.7ex} & 0\hspace{-.7ex} & 1\hspace{-.7ex} & 1\hspace{-.7ex} & 1\hspace{-.7ex} & 0\hspace{-.7ex} & 0\hspace{-.7ex} & 0\hspace{-.7ex} & 
  0\hspace{-.7ex} & 3\hspace{-.7ex} & 2\hspace{-.7ex} & 1\hspace{-.7ex} & 0\hspace{.5ex}$ $\\
  
   \hspace{-.7ex} &  \hspace{-.7ex} &  \hspace{-.7ex} & 1\hspace{-.7ex} & 0\hspace{-.7ex} & \hspace{-1ex}-1\hspace{-.7ex} & 0\hspace{-.7ex} & \hspace{-1ex}-1\hspace{-.7ex} & \hspace{-1ex}-3\hspace{-.7ex} & 0\hspace{-.7ex} & 
  0\hspace{-.7ex} & 0\hspace{-.7ex} & 0\hspace{-.7ex} & 0\hspace{-.7ex} & 0\hspace{-.7ex} & 0\hspace{-.7ex} & 0\hspace{-.7ex} & 1\hspace{-.7ex} & 1\hspace{-.7ex} & 1\hspace{-.7ex} & 
  1\hspace{-.7ex} & 0\hspace{-.7ex} & 1\hspace{-.7ex} & 2\hspace{-.7ex} & 6\hspace{.5ex}$ $\\
   \hline 
   
   \hspace{-.7ex} &  \hspace{-.7ex} &  \hspace{-.7ex} &  \hspace{-.7ex} & 1\hspace{-.7ex} & 0\hspace{-.7ex} & 0\hspace{-.7ex} & 0\hspace{-.7ex} & 0\hspace{-.7ex} & \hspace{-1ex}-1\hspace{-.7ex} & 
  0\hspace{-.7ex} & 0\hspace{-.7ex} & 0\hspace{-.7ex} & 0\hspace{-.7ex} & \hspace{-1ex}-2\hspace{-.7ex} & \hspace{-1ex}-1\hspace{-.7ex} & \hspace{-1ex}-1\hspace{-.7ex} & 0\hspace{-.7ex} & 0\hspace{-.7ex} & 0\hspace{-.7ex} &   
  0\hspace{-.7ex} & 0\hspace{-.7ex} & 0\hspace{-.7ex} & 0\hspace{-.7ex} & 0\hspace{.5ex}$ $\\
  
   \hspace{-.7ex} &  \hspace{-.7ex} &  \hspace{-.7ex} &  \hspace{-.7ex} &  \hspace{-.7ex} & 1\hspace{-.7ex} & 0\hspace{-.7ex} & 0\hspace{-.7ex} & 0\hspace{-.7ex} & 0\hspace{-.7ex} & 
  \hspace{-1ex}-1\hspace{-.7ex} & 0\hspace{-.7ex} & 0\hspace{-.7ex} & 0\hspace{-.7ex} & 0\hspace{-.7ex} & 0\hspace{-.7ex} & 0\hspace{-.7ex} & \hspace{-1ex}-2\hspace{-.7ex} & \hspace{-1ex}-1\hspace{-.7ex} & 0\hspace{-.7ex} & 
  \hspace{-1ex}-1\hspace{-.7ex} & 0\hspace{-.7ex} & 0\hspace{-.7ex} & 0\hspace{-.7ex} & 0\hspace{.5ex}$ $\\
  
   \hspace{-.7ex} &  \hspace{-.7ex} &  \hspace{-.7ex} &  \hspace{-.7ex} &  \hspace{-.7ex} &  \hspace{-.7ex} & 1\hspace{-.7ex} & 0\hspace{-.7ex} & 0\hspace{-.7ex} & 0\hspace{-.7ex} & 
  0\hspace{-.7ex} & \hspace{-1ex}-1\hspace{-.7ex} & 0\hspace{-.7ex} & 0\hspace{-.7ex} & 0\hspace{-.7ex} & \hspace{-1ex}-1\hspace{-.7ex} & 0\hspace{-.7ex} & 0\hspace{-.7ex} & 0\hspace{-.7ex} & 0\hspace{-.7ex} & 
  0\hspace{-.7ex} & \hspace{-1ex}-3\hspace{-.7ex} & \hspace{-1ex}-1\hspace{-.7ex} & 0\hspace{-.7ex} & 0\hspace{.5ex}$ $\\
  
   \hspace{-.7ex} &  \hspace{-.7ex} &  \hspace{-.7ex} &  \hspace{-.7ex} &  \hspace{-.7ex} &  \hspace{-.7ex} &  \hspace{-.7ex} & 1\hspace{-.7ex} & 0\hspace{-.7ex} & 0\hspace{-.7ex} & 
  0\hspace{-.7ex} & 0\hspace{-.7ex} & \hspace{-1ex}-1\hspace{-.7ex} & 0\hspace{-.7ex} & 0\hspace{-.7ex} & 0\hspace{-.7ex} & \hspace{-1ex}-1\hspace{-.7ex} & 0\hspace{-.7ex} & \hspace{-1ex}-1\hspace{-.7ex} & 
   \hspace{-1ex}-2\hspace{-.7ex} &   0\hspace{-.7ex} & 0\hspace{-.7ex} & \hspace{-1ex}-2\hspace{-.7ex} & \hspace{-1ex}-2\hspace{-.7ex} & 0\hspace{.5ex}$ $\\
  
   \hspace{-.7ex} &  \hspace{-.7ex} &  \hspace{-.7ex} &  \hspace{-.7ex} &  \hspace{-.7ex} &  \hspace{-.7ex} &  \hspace{-.7ex} &  \hspace{-.7ex} & 1\hspace{-.7ex} & 0\hspace{-.7ex} & 
  0\hspace{-.7ex} & 0\hspace{-.7ex} & 0\hspace{-.7ex} & \hspace{-1ex}-1\hspace{-.7ex} & 0\hspace{-.7ex} & 0\hspace{-.7ex} & 0\hspace{-.7ex} & 0\hspace{-.7ex} & 0\hspace{-.7ex} & 0\hspace{-.7ex} & 
  \hspace{-1ex}-1\hspace{-.7ex} & 0\hspace{-.7ex} & 0\hspace{-.7ex} & \hspace{-1ex}-1\hspace{-.7ex} & \hspace{-1ex}-4\hspace{.5ex}$ $\\
  \hline   
       
\end{tabular}

\vspace{1ex}

\noindent {\normalsize $D: d_{mn} = d(P_n)^m$ } 


\hspace{-4ex}\begin{tabular}{|ccccccccccccccccccccccccc|}
   \hline  
  \hspace{-1ex}$1^{m\!\!}$\hspace{-.7ex} & \hspace{-1ex}$2^{m\!\!}$\hspace{-.7ex} & \hspace{-1ex}$3^{m\!\!\!}$\hspace{-.7ex} & \hspace{-1ex}$4^{m\!\!\!}$\hspace{-.7ex} & \hspace{-1ex}$4^{m\!\!\!}$\hspace{-.7ex} & 
  \hspace{-1ex}$5^{m\!\!\!}$\hspace{-.7ex} & \hspace{-1ex}$5^{m\!\!\!}$\hspace{-.7ex} & \hspace{-1ex}$6^{m\!\!\!}$\hspace{-.7ex} &\hspace{-1ex}$8^{m\!\!\!}$\hspace{-.7ex}  & \hspace{-1ex}$5^{m\!\!\!}$\hspace{-.7ex} & 
  \hspace{-1ex}$6^{m\!\!\!}$\hspace{-.7ex} & \hspace{-1ex}$6^{m\!\!\!}$\hspace{-.7ex} & \hspace{-1ex}$7^{m\!\!\!}$\hspace{-.7ex} & \hspace{-1ex}$9^{m\!\!\!}$\hspace{-.7ex} & \hspace{-1ex}$6^{m\!\!\!}$\hspace{-.7ex} & 
  \hspace{-1ex}$7^{m\!\!\!}$\hspace{-.7ex} & \hspace{-1ex}$8^{m\!\!\!}$\hspace{-.7ex} & \hspace{-1ex}$7^{m\!\!\!}$\hspace{-.7ex} & \hspace{-1ex}$8^{m\!\!\!}$\hspace{-.7ex} & \hspace{-1ex}$9^{m\!\!\!}$\hspace{-.7ex} & 
  \hspace{-1.7ex}$10^{m\!\!\!}$\hspace{-.7ex} & \hspace{-1ex}$9^{m\!\!\!}$\hspace{-.7ex} & \hspace{-1.7ex}$10^{m\!\!\!}$\hspace{-.7ex} & \hspace{-1.7ex}$12^{m\!\!\!}$\hspace{-.7ex} & \hspace{-1.7ex}$16^{m\!\!\!}$\hspace{-.7ex}\\
  \hline
\end{tabular}
   
\vspace{1ex}

\noindent {\normalsize $E^T: e_{mn} = e(m,P_n) = \sharp\, \mathcal{E}(\underline{m},P_n)= \sharp\, \mathcal{F}(\underline{m},P_n)= \sharp\, \mathcal{G}(\underline{m},P_n)= \sharp\, \mathcal{H}(\underline{m},P_n)$ } 

\vspace{-2ex}

\hspace{-4ex}\begin{tabular}{|rrrrrrrrrrrr|}
   \hline  
  1\hspace{-.7ex} & \hspace{-.7ex}1\hspace{-.3ex} & 1\hspace{-.7ex} & 1\hspace{-1.7ex} & 1\hspace{-1.7ex} & 1\hspace{-1.7ex} & 1\hspace{-1.7ex} & 1\hspace{-1.7ex} & 1\hspace{-1.7ex} & \hspace{-1ex}$\hdots$\,& & \hspace{7.6ex}$1\!$\\
  1\hspace{-.7ex} & \hspace{-.7ex}3\hspace{-.3ex} & 7\hspace{-.7ex} & \hspace{-1ex}15\hspace{-1.7ex} & \hspace{-1ex}31\hspace{-1.7ex} & \hspace{-1ex}63\hspace{-1.7ex} & \hspace{-1ex}127\hspace{-1.7ex} \hspace{-1ex} & 
  \hspace{-1ex}255\hspace{-1.7ex} & \hspace{-1ex}511\hspace{-1.7ex} & \hspace{-1ex}$\hdots$\, &  & \hspace{-2.5ex}$2^m\!-1\!$\\
  1\hspace{-.7ex} &  \hspace{-.7ex}5\hspace{-.3ex} & \hspace{-1ex}19\hspace{-.7ex} & \hspace{-1ex} 65\hspace{-1.7ex} & \hspace{-1ex} 211\hspace{-1.7ex} & \hspace{-1ex} 665\hspace{-1.7ex} & \hspace{-1ex} 2059\hspace{-1.7ex} &
  \hspace{-1ex} 6305\hspace{-1.7ex} & \hspace{-1ex} 19171\hspace{-1.7ex} & \hspace{-1.5ex}$\hdots$\, &  & \hspace{-3.5ex}$3^m\!-2^m\!$\\
  1\hspace{-.7ex} &  \hspace{-.7ex}9\hspace{-.3ex} & \hspace{-1ex}49\hspace{-.7ex} & \hspace{-2ex}225\hspace{-1.7ex} & \hspace{-1ex}961\hspace{-1.7ex} & \hspace{-1ex}3969\hspace{-1.7ex} & \hspace{-1ex}16129\hspace{-1.7ex} & 
  \hspace{-1ex}65025\hspace{-1.7ex} & \hspace{-1ex}261121\hspace{-1.7ex} & \hspace{-1.5ex}$\hdots$\, &  & \hspace{-11.5ex}$4^m\!-2\cdot 2^m +1\!$\\
  1\hspace{-.7ex} &  \hspace{-.7ex}7\hspace{-.3ex} &  \hspace{-1ex}37\hspace{-.7ex} &  \hspace{-1ex}175\hspace{-1.7ex} &  \hspace{-1ex}781\hspace{-1.7ex} &  \hspace{-1ex}3367\hspace{-1.7ex} & \hspace{-1ex}14197\hspace{-1.7ex} &
  \hspace{-1ex}58975\hspace{-1.7ex} &  \hspace{-1ex}242461\hspace{-1.7ex} & \hspace{-1.5ex}$\hdots$\, &  & \hspace{-8.5ex}$4^m\!-3^m\!$\\
  1\hspace{-.7ex} &  \hspace{-.7ex}9\hspace{-.3ex} &  \hspace{-1ex}61\hspace{-.7ex} &  \hspace{-1ex}369\hspace{-1.7ex} &  \hspace{-1ex}2101\hspace{-1.7ex} &  \hspace{-1ex}11529\hspace{-1.7ex} & \hspace{-1ex}61741\hspace{-1.7ex} & 
  \hspace{-1ex}325089\hspace{-1.7ex} & \hspace{-1ex}1690981\hspace{-1.7ex} & \hspace{-1.5ex}$\hdots$\, &  & \hspace{-8.5ex}$5^m\!-4^m\!$\\
  1\hspace{-.7ex} & \hspace{-.7ex}11\hspace{-.3ex} &  \hspace{-1ex}79\hspace{-.7ex} &  \hspace{-1ex}479\hspace{-1.7ex} & \hspace{-1ex}2671\hspace{-1.7ex} & \hspace{-1ex}14231\hspace{-1.7ex} & \hspace{-1ex}73879\hspace{-1.7ex} & 
  \hspace{-1ex}377759\hspace{-1.7ex} & \hspace{-1ex}1914271\hspace{-1.7ex} & \hspace{-1.5ex}$\hdots$\, &  & \hspace{-14.5ex}$5^m\!-2\cdot 3^m + 2^m\!$\\
  1\hspace{-.7ex} & \hspace{-.7ex}15\hspace{-.3ex} & \hspace{-1ex}133\hspace{-.7ex} &  \hspace{-1ex}975\hspace{-1.7ex} & \hspace{-.1ex}6541\hspace{-1.7ex} & \hspace{-.1ex}41895\hspace{-1.7ex} & \hspace{-1ex}261493\hspace{-1.7ex} & 
  \hspace{-1ex}1607775\hspace{-1.7ex} & \hspace{-1ex}9796381\hspace{-1.7ex} & \hspace{-1.5ex}$\hdots$\, &  & \hspace{-5ex}$6^m\!- 4^m - 3^m + 2^m\!$\\
  1\hspace{-.7ex} & \hspace{-.7ex}27\hspace{-.3ex} & \hspace{-1ex}343\hspace{-.7ex} &3375\hspace{-1.7ex} & 29791\hspace{-1.7ex} & 250047\hspace{-1.7ex} & 2048383\hspace{-1.7ex} & 16581375\hspace{-1.7ex} &
  133432831\hspace{-1.7ex} &\hspace{-1.5ex}$\hdots$\, &  & \hspace{-22.5ex}$8^m\!-\! 3\cdot 4^m \!+\! 3\cdot 2^m \!-\! 1\!$\\ 
   \hline
\end{tabular}
   
\end{footnotesize}

\end{document}